\documentclass[11pt,twoside]{article}
\usepackage[latin1]{inputenc}
\usepackage{latexsym}
\usepackage[german]{babel}
\usepackage{amsmath}
\usepackage{amssymb,amsfonts}
\title{{\bf\hspace*{1.2cm} Kohomologie mit Schranken\\\hspace*{1.2cm} und \\
\hspace*{1.2cm} Fortsetzung
holomorpher Funktionen\\
\hspace*{1.2cm} durch\\
\hspace*{1.2cm} lineare stetige Operatoren}}
\author{\hspace*{0cm}Dissertation\\ zur Erlangung eines Doktors der \\
Naturwissenschaften\\
(doctor rerum naturalium)\\[4cm] Dem Fachbereich Mathematik der \\
Bergischen Universität - Gesamthochschule Wuppertal \\
vorgelegt von \\
{\sc Dipl.-Math. M. Matthias Schmitt}\\
aus Wuppertal\\[1cm]
Tag der mündlichen Prüfung: 29. Oktober 2001\\[1cm]
\begin{tabular}{lll}
Gutachter: & 1.& {\sc Prof. Dr. D. Vogt},
 Bergische Universität - Gesamthochschule Wuppertal\\
 & 2. & {\sc Prof. Dr. M. Langenbruch}, Carl von Ossietzky Universität Oldenburg
\end{tabular} }
\date{}

\newcommand{\bild}{\mbox{{\rm Bild}}\hspace{0.1em}}
\newcommand{\relkom}{\subset\subset}
\newcommand{\sumlim}{\sum\limits}
  
\newcommand{\C}{\mathbb{C}}
\newcommand{\R}{\mathbb{R}}

\newcommand{\N}{\mathbb{N}}

\newcommand{\komma}{\leer{0.2} , \leer{0.5}}
\newcommand{\doppelpunkt}{\leer{0.1} : \leer{0.3}}
\newcommand{\punkt}{\leer{0.2} . \leer{0.5}}
\newcommand{\abbkern}{\mbox{Kern\,}}

\newcommand{\eps}{\varepsilon}
\renewcommand{\phi}{\varphi}
\renewcommand{\rho}{\varrho}

\newcommand{\norm}[2]{\parallel #2 \parallel_{#1}}

\newcommand{\normoben}[3]{\parallel #3 \parallel_{#2}^{#1}}
\newcommand{\dreinorm}[3]{\parallel\hspace{-1pt}\mid #3 \parallel\hspace{-1pt}\mid_{#2}^{#1}}

\newcommand{\betrag}[1]{\mid #1\mid}
\renewcommand{\hat}[1]{\widehat{#1}}
\renewcommand{\tilde}[1]{\widetilde{#1}}

\renewcommand{\lim}{\mathop{\rm lim}}
\renewcommand{\sup}{\mathop{\rm sup}}

\newcommand{\F}{{\cal F}}
\renewcommand{\H}{{\cal H}}
\newcommand{\U}{{\cal U}}
\newcommand{\J}{{\cal J}}
\newcommand{\W}{{\cal W}}
\newcommand{\Rscript}{{\cal R}}
\newcommand{\oka}{{\cal O}}
\newcommand{\be}{\begin{enumerate}}
\newcommand{\ee}{\end{enumerate}}

\newcommand{\normleer}{\parallel \hspace{0.8em}\parallel}
\newcommand{\betragleer}{\mid \hspace{0.8em}\mid}

\newcommand{\sss}{{\sigma}}
\newcommand{\sminuseins}{{\sigma -1}}

\newcommand{\indnullsigma}{{{i_0,\ldots,i_{\sigma}}}}
\newcommand{\indeinssigma}{{{i_1,\ldots,i_{\sigma}}}}
\newcommand{\indzweisigma}{{{i_2,\ldots,i_{\sigma}}}}

\newcommand{\indnullsigmafehl}[1]{{{i_0,\ldots,\hat{i}_{#1},\ldots,i_{\sigma}}}}
\newcommand{\indeinssigmafehl}[1]{{{i_1,\ldots,\hat{i}_{#1},\ldots,i_{\sigma}}}}

\newcommand{\leer}[1]{\hspace{#1em}}

\newcommand{\PSH}{\mbox{PSH}}

\newcommand{\foralls}{\mbox{ für alle }}
\newcommand{\und}{\mbox{ und }}

\newcommand{\diam}{\mbox{diam}}

\newcommand{\einlabeln}[1]{{\label{{#1}}}
                   }
\newcommand{\auslabeln}[1]{(\ref{{#1}})}

\newcommand{\einlabelnsatz}[1]{\label{{#1}}
                   }
\newcommand{\auslabelnsatz}[1]{\ref{{#1}}}

\newcommand{\Ende}{$\Box$}

\newcommand{\dquer}{\overline{\partial}}

\newcommand{\ffa}[2]{\sum_{#1=1}^{#2}}

\newcommand{\ffb}[2]{\sup_{#1}\mid #2 \mid}

\newcommand{\ffc}[2]{#1_1,\ldots,#1_{#2}}

\newcommand{\ffd}[2]{l_1(H_{#1}(#2))}

\newcommand{\ffe}[4]{#1_{#2_{#3} + #4}}

\newcommand{\ffg}[2]{{i_{#1},\ldots,i_{\sigma  #2}}}

\newcommand{\zza}{H_U(\F)}

\newcommand{\zzb}{H_U(\oka)}

\newcommand{\zzc}{l_1(\zzb)}

\newcommand{\zzd}{{\cal{R}}}

\newenvironment{satz}[1]%
{
\refstepcounter{subsection} 
{\normalsize\bf #1
 \bf \arabic{section}.\arabic{subsection}. :} \hspace{\fill}\\
\leftmargin3em
\it }{

\sc Beweis:}
{
\refstepcounter{subsection} {\normalsize\bf #1
 \bf \arabic{section}.\arabic{subsection}. \bf #2:} \hspace{\fill}\\
\leftmargin3em
\it }{

\sc Beweis:}
\newenvironment{satzohnebeweis}[1]%
{
\refstepcounter{subsection} {\normalsize\bf #1
 \bf \arabic{section}.\arabic{subsection}. :} \hspace{\fill}\\
\leftmargin3em
\it }{

}
{
\refstepcounter{subsection} {\normalsize\bf #1
 \bf \arabic{section}.\arabic{subsection}. \bf #2:} \hspace{\fill}\\
\leftmargin3em
\it }{}

\newcommand{\formel}[1]{

\parbox{14.7cm}{\begin{eqnarray*}\refstepcounter{subsubsection} #1 \end{eqnarray*}}
\hfill \parbox{1.4cm}{\rm 
(\arabic{section}.\arabic{subsection}.\arabic{subsubsection})}

 }
\newcommand{\formelohne}[1]{

{\begin{eqnarray*} #1 \end{eqnarray*}}

 }
\begin{document}
\maketitle
\thispagestyle{empty}
\tableofcontents \newpage

\section[{\sc Einführung und Ergebnisse}]{Einführung und Ergebnisse}

Ein häufig von verschiedenen Autoren behandeltes Problem ist die
Frage nach Bedingungen oder Kriterien für die Existenz einer 
Rechtsinversen $r$ zu einer surjektiven Abbildung $q$. Diese Fragestellung
tritt oft in Form  einer
kurzen exakten Sequenz von Fr\'echeträumen
\[   0\to E \stackrel{i}{\to} G \stackrel{ q}{\to} F \to 0 \]
auf. Hier ist also $q\circ r=id_F$. $F$ ist in diesem Fall isomorph
zu einem projizierten Unterraum von $G$, nämlich Bild $r$. 
Man sagt auch, daß die Sequenz zerfällt.

Bekannte Beispiele für eine solche Fragestellung ergeben sich in
der Komplexen Analysis in den Situationen, in denen $F$ der Raum der
holomorphen Funktionen auf einer komplexen Untermannigfaltigkeit $V$
einer Steinschen Mannigfaltigkeit $X$ ist.

Ist $X$ eine Steinsche Mannigfaltigkeit, dann sei stets $H(X)$ der
Raum der holomorphen Funktionen auf $X$. Für Definitionen und
Sätze in diesem Zusammenhang sei auf das Buch von Hörmander \cite{Hoer},
Kapitel 5, verwiesen. Sei $\oka_X$ die Garbe
von Keimen holomorpher Funktionen auf $X$.
 Es sei mit $\Gamma(X,J_V)$ der Raum der
globalen Schnitte über $X$ mit Werten in $J_V$ bezeichnet, wobei
$J_V$ die Idealgarbe in der Garbe $\oka_X$ der Keime holomorpher
Funktionen auf $X$ ist, die auf
$V$ verschwinden (vgl. Kapitel 7 in \cite{Hoer}). $\Gamma(X,J_V)$ können
wir mit dem Teilraum der  holomorphen Funktionen auf $X$ 
identifizieren, die auf $V$ verschwinden.

Nach der Oka-Cartan-Theorie ist die Sequenz
\[  0\to \Gamma(X,J_V)\hookrightarrow H(X) \stackrel{R}{\to} H(V) \to 0 
            \leer{2} (\ast) \]
stets exakt. $R$ ist hierbei die Restriktionsabbildung auf $V$.
  Die Surjektivität von $R$ bedeutet gerade die globale
Fortsetzbarkeit 
 von holomorphen Funktionen auf $V$ zu solchen auf $X$ 
(Theorem 7.4.8 in \cite{Hoer}).

Mityagin und Khenkin haben eine Rechtsinverse zu $R$ in $(\ast)$ konstruiert 
für den Fall, daß $X$ ein Gebiet in einem $\C^n$ ist, das
sich darstellen läßt durch eine strikt plurisubharmonische
Funktion $u$ auf einer Umgebung $X'$ von $\overline{X}$ als
$X=\{z\in X'\doppelpunkt u(z)<0 \}$, wobei zusätzlich grad $u\not=0$
auf $\partial X$ ist. Es wird weiter vorausgesetzt, daß $V$ den Rand
von $X$ transversal schneidet
(Theorem 4.2 in \cite{M-H}).

Haupthilfsmittel in \cite{M-H} ist eine Abwandlung von Hörmanders
Theorem B mit Schranken (Theorem 7.6.10 in \cite{Hoer}), angewandt
auf die Garbe $J_V$, die in diesem Fall durch endlich viele globale 
Schnitte erzeugt wird. Im allgemeinen können wir nicht erwarten, daß
 die Idealgarbe einer analytischen Untervarietät des $\C^N$ durch endlich
viele globale Schnitte  erzeugt wird.

Wir werden in dieser Arbeit die Technik von Mityagin und Khenkin
auf Situationen übertragen, in denen keine endliche Menge von
erzeugenden Schnitten für $J_V$ existiert.

Von großer Bedeutung sind in diesem Kontext
 Potenzreihenräume, da für viele Räume holomorpher Funktionen Darstellungen
als Potenzreihenräume bekannt sind.
Ist $\alpha=(a_k)_{k\in \N_0}$ eine Folge positiver reeller Zahlen mit 
$a_k\nearrow \infty$ für $k\to \infty$ und $r\in \{0,\infty\}$, dann sei
\[   \Lambda_r(\alpha):= \{ x=(x_0,x_1,\ldots)\in \C^{\N_0}\doppelpunkt
     \norm{\rho}{x} := \sum_{k=0}^\infty \betrag{x_k} e^{\rho\, a_k} <\infty
     \mbox{ für alle } \rho<r \}    \]
Im Falle $r=0$ heißt $\Lambda_r(\alpha)$ Potenzreihenraum endlichen Typs,
im Falle $r=\infty$ heißt $\Lambda_r(\alpha)$ Potenzreihenraum 
unendlichen Typs.

Für die Charakterisierung von Fr\'echeträumen mit aufsteigendem 
Fundamentalsystem von Halbnormen $(\normleer_n)_n$, die
zu einem Potenzreihenraum isomorph sind, sind die folgenden von 
D. Vogt eingeführten linear topologischen Invarianten hilfreich
(vgl. Kapitel 29 in \cite{M-V}), wobei diese Eigenschaften auch in vielen anderen
Zusammenhängen bedeutsam sind. Sei hierfür $(E,\normleer_n)$ ein Fr\'echetraum:

\begin{enumerate}

\item[$(\underline{DN})$:] Es gibt ein $n_0\in \N$, so daß für jedes $k\in \N$ ein
$K\in \N$, ein $C>0$ und ein $0<\tau<1$ existiert mit
\[   \norm{k}{x} \le C\norm{n_0}{x}^{1-\tau} \norm{K}{x}^{\tau} 
                                       \mbox{ für jedes } x \in E \punkt \]
\item[$(DN)$:]  Es gibt ein $n_0\in \N$, so daß für jedes $k\in \N$ ein
$K\in \N$ und ein $C>0$  existiert mit
\[   \norm{k}{x}^2 \le C\norm{n_0}{x} \norm{K}{x} \mbox{ für jedes } x \in E \punkt \]
\item[$(\Omega)$:]  Für jedes  $n\in \N$ gibt es ein $k\in \N$, so daß zu jedem 
$K\in \N$ ein $C>0$ und ein $0<\tau<1$ existiert mit
\[   \norm{k}{y}^{\ast} \le C\norm{n}{y}^{\ast 1-\tau} \norm{K}{y}^{\ast \tau} 
                          \mbox{ für jedes } y \in E'\punkt \]
\item[$(\overline{\Omega})$:]  Für jedes  $n\in \N$ gibt es ein $k\in \N$, so daß zu jedem 
$K\in \N$ ein $C>0$  existiert mit
\[   \norm{k}{y}^{\ast 2} \le C{\norm{n}{y}^{\ast}} {\norm{K}{y}^{\ast}} 
                          \mbox{ für jedes } y \in E'\punkt \]

\end{enumerate}

Hierbei sei $\normleer_k^{\ast}$ die duale Norm zu der Halbnorm $\normleer_k$, 
d.h. für ein Dualelement $\Phi$ ist $\norm{k}{\Phi}^{\ast}
=\sup \{ \betrag{\Phi(x)} \doppelpunkt 
\norm{k}{x}\le 1\} \in \R\cup \{\infty\}$.
  
Wir unterscheiden die Fälle, in denen $H(X)$ isomorph zu einem 
Potenzreihenraum endlichen Typs ist, von den Fällen, in denen
$H(X)$ isomorph zu einem Potenzreihenraum unendlichen Typs ist.
Es muß nicht zwingend einer der beiden Fälle vorliegen.
Wir sprechen vom endlichen bzw. unendlichen Fall.

Mityagin und Khenkin haben gezeigt, daß die Isomorphie von $H(V)$
und $\Lambda_1(\alpha)$ in dem von Ihnen betrachteten endlichen Fall hinreichend
für die Existenz einer Rechtsinversen der Restriktionsabbildung ist
(Proposition 4.2 in \cite{M-H}). In allen endlichen 
Fällen ist die Isomorphie von $H(V)$
und $\Lambda_1(\alpha)$ nach einem Ergebnis von Mityagin aber auch notwendig, denn jeder
unendlichdimensionale, komplementierte Unterraum eines Potenzreihenraumes
endlichen Typs ist isomorph zu einem Potenzreihenraum
endlichen Typs (siehe z.B. Folgerung 29.20 in \cite{M-V}). Im unendlichen
Fall gibt es bislang keine analoge Aussage.

Sei jetzt der unendliche Fall betrachtet, daß $V$ eine abgeschlossene 
Untermannigfaltigkeit des $\C^N$ ist (vgl. \cite{Hoer}, Def. 5.1.4.).
 $V$ ist dann lokal endlich, d.h.
jeder Punkt des $\C^N$ besitzt eine Umgebung, die nur endlich
viele Zusammenhangskomponenten schneidet.

$H(\C^N)$ soll im folgenden stets die Fr\'echetraumtopologie tragen, die
durch die Suprema auf einer abzählbaren kompakten Ausschöpfung des 
$\C^N$ erzeugt wird. Die Fr\'echetraumtopologie von $H(V)$ soll für 
eine abgeschlossene Untermannigfaltigkeit $V$ stets durch die 
Spuren der kompakten Ausschöpfung des $\C^N$ gebildet werden.

Wir wenden uns also solchen Untermannigfaltigkeiten $V$ des $\C^N$ mit
Dimension $d$ zu, für die $(\ast)$ zerfällt, d.h. für die
es einen linearen stetigen Operator $E:H(V)\to H(\C^N)$ gibt mit
$(Ef)\mid_V=f$ für jedes $f\in H(V)$.

Für die Existenz eines solchen Ausdehnungsoperators
 haben verschiedene Autoren notwendige und hinreichende Bedingungen gezeigt.
Wir fassen diese in folgendem Satz zusammen und nennen die Zuordnung
zu den einzelnen Autoren im Beweis:

\begin{satz}{Satz}\einlabelnsatz{15.2}
Äquivalent sind:
\be
\item[1)] Es gibt einen linearen stetigen Operator $E:H(V)\to H(\C^N)$
mit $(Ef)\mid_V=f$ für jedes $f\in H(V)$. 
\item[2)] $H(V)$ besitzt (DN).
\item[3)] Jede beschränkte plurisubharmonische Funktion auf $V$
ist konstant (Stark-Liouvillesche Eigenschaft).
\item[4)] $H(V)$ ist isomorph zu $\Lambda_\infty(k^{\frac{1}{d}})$.
\item[5)] $H(V)$ ist isomorph zu $H(\C^d)$.
\item[6)] Es gibt eine Basisfolge $(f_k)_{k\in \N_0}$ in $H(V)$, so daß
für eine fest gewählte Folge $\rho_n\nearrow \infty$ gilt:
   \be
   \item[a)] Zu jedem $n\in \N$ gibt es ein $m\in \N$ und ein 
    $C>0$, so daß 
    \[   \norm{n}{f_k}\le C \rho_m^{(k^{\frac{1}{d}})} \foralls k\in \N_0\punkt \]
   \item[b)] Zu jedem $n\in \N$ gibt es ein $m\in \N$ und ein 
    $C>0$, so daß
    \[   \rho_n^{(k^{\frac{1}{d}})} \le C \norm{m}{f_k} \foralls k\in \N_0\punkt \]
   \ee
\ee
\end{satz}

Für 2)$\Leftrightarrow$1) ist der Splittingsatz von Vogt für kurze
exakte Sequenzen von Fr\'echet-Hilberträumen anwendbar (siehe z.B. \cite{M-V}, Satz 30.1).
Die hier vorliegenden Fr\'echeträume sind nuklear und daher Fr\'echet-Hilbert.

Der Vogtsche Splittingsatz besagt, daß $(\ast)$ (mit $X=\C^N$) spaltet, 
falls $H(V)$ die Eigenschaft (DN)
und $\Gamma(\C^N,J_V)$ die Eigenschaft $(\Omega)$ besitzt. Letzteres ist
für jede analytische Teilmenge des $\C^N$ der Fall (siehe z.B. \cite{V1}).

Wenn umgekehrt  $(\ast)$ spaltet,
dann ist $H(V)$ isomorph zu einem projizierten Teilraum von $H(\C^N)$.
Auf abgeschlossene Teilräume vererbt sich stets die Eigenschaft (DN).
Also folgt, daß $H(V)$ auch die Eigenschaft (DN) besitzt.

4)$\Leftrightarrow$5) ist elementar (siehe z.B. \cite{Ayt2}).

3)$\Leftrightarrow$2) findet man  in \cite{Ayt2},  Theorem I.12.
Auch Vogt und Zaharyuta sind unabhängig zu diesem Ergebnis
gekommen (siehe auch \cite{Ayt1};
bzw. \cite{Z}).

4)$\Rightarrow$6) ist elementar. Man nehme als Basisfolge die Urbilder des
Isomorphismus $H(V)\to \Lambda_\infty(k^{\frac{1}{d}})$.

Für  6)$\Rightarrow$4) zeigen wir, daß $T:\Lambda_\infty(k^{\frac{1}{d}})
\to H(V)$, $(\lambda_k)_{k\in \N_0} \mapsto \sum_{k\in \N_0} \lambda_k f_k$
ein Isomorphismus ist.

Einerseits gibt es wegen a) zu jedem $n\in \N$  ein $m\in \N$ und ein 
    $C>0$, so daß für  jedes $\lambda:=(\lambda_k)_{k\in \N_0}$ in 
$\Lambda_\infty(k^{\frac{1}{d}})$ gilt:
\[ \norm{n}{T\lambda} \le \sum_{k\in \N_0} \betrag{\lambda_k} \norm{n}{f_k}
   \le C  \sum_{k\in \N_0} \betrag{\lambda_k}\rho_m^{(k^{\frac{1}{d}})} 
   = C\norm{m}{\lambda} \punkt\]
Andererseits ist $H(V)$ nuklear, und daher ist $(f_k)_k$ in $H(V)$ absolut
nach dem Satz von Dynin-Mityagin (siehe z.B. \cite{M-V}, Satz 28.12). Mit 
b) folgt also, daß zu jedem $n\in \N$  ein $m\in \N$ und ein 
    $C>0$ existiert mit 
\[     \sum_{k\in \N_0} \betrag{\lambda_k}\rho_n^{(k^{\frac{1}{d}})} \le
       C\norm{m}{\sum_{k\in \N_0} \lambda_k f_k}  \foralls 
       \sum_{k\in \N_0} \lambda_k f_k \in H(V)  \punkt \]

5)$\Rightarrow$2) ist klar, da $H(\C^d)$ die Eigenschaft (DN) besitzt.

Es verbleibt 2)$\Rightarrow$5) zu zeigen.

Im unendlichen Fall können wir nicht wie im endlichen Fall ein allgemeines
Argument verwenden, denn es ist  offen, ob jeder unendlichdimensionale,
komplementierte Unterraum eines Potenzreihenraumes unendlichen Typs
wieder isomorph zu einem  Potenzreihenraum unendlichen Typs ist
(vgl. \cite{M-V} S. 354).
 
In ihrer Arbeit \cite{A-K-T} konnten  aber Aytuna, Krone und Terzio\v{g}lu unter
Verwendung eines Ergebnisses von Vogt zeigen, daß 5) aus 2)
folgt (\cite{A-K-T}, Proposition 2.1). Damit ist der Beweis vollständig \Ende 

Anmerkungen zum Beweis zu Satz \auslabelnsatz{15.2}:

Ist $V$ eine Steinsche Mannigfaltigkeit, so kann man nach dem bekannten Einbettungssatz
(siehe z.B. \cite{Hoer}, Theorem 5.3.9) $V$ als abgeschlossene Untermannigfaltigkeit
in den $\C^N$ einbetten. Die Existenz des Extensionsoperators hängt nach Satz
\auslabelnsatz{15.2} nicht von der Art der Einbettung ab. 

In Proposition 6.4 in \cite{M-H} stellten Mityagin und Khenkin die These auf,
daß 1) aus 5)  folgt, wenn man die aus der Isomorphie in 5) herrührende
Basisfolge in $H(V)$ gliedweise mit geeigneten Schranken zu einer Folge ganzer Funktionen
fortsetzen kann. Dies ist exakt die Technik in \cite{M-H} zur Konstruktion einer
 Rechtsinversen in der Situation,
wo $V$ eine 
abgeschlossene Untermannigfaltigkeit in einem strikt pseudokonvexen Gebiet ist, die
den Rand transversal schneidet. 

Will man der Argumentation in \cite{M-H} unter Verwendung von 
Hörmanders Theorem B mit Schranken folgen, so ergibt sich als
wesentliches Hindernis, daß die Idealgarbe $J_V$ im allgemeinen nicht
durch die Keime endlich vieler globaler Schnitte über dem $\C^N$ 
erzeugt werden kann, so wie das über beschränkten pseudokonvexen
Gebieten möglich ist (siehe z.B. \cite{Hoer} Theorem 7.2.1 in Verbindung mit 
Theorem 7.2.9).

Wir werden in dieser Arbeit zeigen, wie man dieses Hindernis überwinden
kann. Damit wird die Konstruktion eines Ausdehnungsoperators durch 
Fortsetzungen von Basisfunktionen in $H(V)$ auch im 
unendlichen Fall  möglich.

Insbesondere geht es um die Frage, ob dieser konstruktive Ansatz weiterführende
Informationen über die \glqq Qualität\grqq\  der Stetigkeit des
Ausdehnungsoperators liefert. Das bedeutet für jeweils gegebene Gradierungen
auf $H(V)$ und $H(X)$, zu fragen, wie groß man zu einer vorgegebenen Stufe
in $H(X)$ die Stufe in $H(V)$ wählen muß, damit sich die Fortsetzungen,
die der Ausdehnungsoperator liefert, gegen die Urbildfunktionen jeweils abschätzen.

Tatsächlich ermöglicht uns dieser Ansatz, eine Reihe von Kriterien
und Bedingungen analog zu denjenigen aus Satz 
\auslabelnsatz{15.2} verschiedener Autoren zu einem Satz über die 
Existenz eines linear zahmen Ausdehnungsoperators zusammenzufassen.

Wir werden also die Argumentation von Mityagin und Khenkin wieder
aufgreifen und ein Theorem B mit Schranken für lokal endlich erzeugte Untergarben von 
$\oka^p$, $p\in \N$, in einer Steinschen Mannigfaltigkeit $\Omega$ zeigen.
Hiermit ist eben auch die Idealgarbe $J_A$ einer beliebigen analytischen
Teilmenge $A$ von $\Omega$ erfaßt.

Der Beweis umfaßt die Kapitel 2 bis 5.
Die Grundstruktur ist ein Induktionsverfahren (vgl. Beweis zu Theorem 7.6.10 in \cite{Hoer}).

Wir betrachten dazu das folgende kommutative Diagramm von 
linearen  Räumen und linearen Abbildungen:

\formel{\begin{array}{lllllllll} 
     &     & 0       &     & 0        &    &   0     &      &                      \\
     &     & \uparrow&     & \uparrow &    &\uparrow &      &                      \\
\cdots & \stackrel{P_{1,2}}{\to} & F_{1,2} & \stackrel{P_{1,1}}{\to} & F_{1,1} 
       & \stackrel{P_{1,0}}{\to} & F_{1,0} & \to & 0                               \\
     &     & \uparrow d_{1,2} & & \uparrow d_{1,1} & & \uparrow d_{1,0} & &        \\
\cdots & \stackrel{P_{2,2}}{\to} & F_{2,2} & \stackrel{P_{2,1}}{\to} & F_{2,1} 
       & \stackrel{P_{2,0}}{\to} & F_{2,0} & \to & 0                               \\
     &     & \uparrow d_{2,2} & & \uparrow d_{2,1} & & \uparrow d_{2,0} & &        \\ 
\cdots & \stackrel{P_{3,2}}{\to} & F_{3,2} & \stackrel{P_{3,1}}{\to} & F_{3,1} 
       & \stackrel{P_{3,0}}{\to} & F_{3,0} & \to & 0                               \\
     &     & \uparrow d_{3,2} & & \uparrow d_{3,1} & & \uparrow d_{3,0} & &        \\
     &     &\,\vdots  &     &\,\vdots  &    &\,\vdots   &      &   
\end{array} \einlabeln{15.1.1}  }

In diesem Diagramm seien sämtliche Zeilen und Spalten exakt bis auf die
erste Spalte von rechts. Diese sei lediglich ein Komplex.

Die folgende Induktion zeigt nun, daß auch die erste Spalte exakt ist.
Dazu zeigen wir indukiv über aufsteigendes $n$, 
daß für jedes $k=0,1,\ldots$ gilt: 
Zu jedem $x_k\in F_{n,k}$ mit $ d_{n-1,k}\, x_k=0$ und $P_{n,k-1}\,x_k=0$
gibt es ein $y_k\in F_{n+1,k}$ mit $d_{n,k}\,y_k=x_k$ und $P_{n+1,k-1}\, y_k=0$
(hierbei setzen wir $P_{n,-1}:=0$).

Der Induktionsanfang ist klar, wenn wir uns vorstellen, man würde 
das Diagramm nach obenhin durch Nullräume  und Nullabbildungen erweitern.

Um von $n$ nach $n+1$ zu gelangen, sei $x\in F_{n+1,k}$ mit $d_{n,k}\,x=0$ und 
$P_{n+1,k-1}\,x=0$. Wegen der Exaktheit der Zeilen gibt es ein $y_1\in F_{n+1,k+1}$,
so daß $P_{n+1,k}\, y_1=x$. Sei $y_2=d_{n,k+1} \,y_1$, dann ist wegen
der Kommutativität des Diagramms $P_{n,k}\,y_2=d_{n,k}P_{n+1,k}\,y_1=d_{n,k}\,x=0$.

Da auch $d_{n-1,k+1}\,y_2=0$ ist, gibt es nach der Induktionsvoraussetzung
ein $y_3\in F_{n+1,k+1}$ mit $d_{n,k+1}\,y_3=y_2$ und $P_{n+1,k}\,y_3=0$.
Es liegt $y_4:= y_1-y_3$ im Kern von $d_{n,k+1}$. Da $k+1\ge 1$ ist, ist
die $k+1$-te Spalte exakt. Also gibt es $y_5\in F_{n+2,k+1}$ mit 
$d_{n+1,k+1}\, y_5 =y_4$. Man rechnet leicht nach, daß für 
$y:=P_{n+2,k}\,y_5\in F_{n+2,k}$ gilt $d_{n+1,k}\,y=x$ und $P_{n+2,k-1}\,y=0$.

Für die Definitionen der im folgenden vorkommenden Begriffe wie dem des Kokettenraumes
und des Korandoperators $\delta$ 
verweisen wir auf \cite{Hoer}, Chapter 7. Durch den Korandoperator werden Koketten
der Länge $\sss$ auf Koketten der Länge $\sss+1$ abgebildet.

Das klassische Beispiel für Koketten ist die Cousin-I-Verteilung 
bezüglich einer Überdeckung $(U_i)_{i\in \N}$ eines Gebietes in 
$\C$ bzw. allgemeiner
einer Steinschen Mannigfaltigkeit. Die Cousin Data sind nichts anderes
als eine Kokette $c=(c_{ij})_{ij}$ der Länge 1, bestehend aus 
holomorphen Funktionen $c_{ij}$ auf $U_i\cap U_j$, für die
$\delta c=0$ gilt (vgl. \cite{Hoer}, Theorem 5.5.1). 
Die Lösung $c'=(c'_i)_{i}$ des Cousin-I-Problems ist eine Kokette
der Länge 0, bestehend aus holomorphen Funktionen $c_i$ auf $U_i$,
für die $\delta c'=c$ gilt. Aus der Lösung des 
Cousin-I-Problems auf $\C$ kann man z.B. den Satz von Mittag-Leffler
folgern.

Die Aussage des Theorem B von Cartan, nämlich daß die
$p$-te Kohomologiegruppe $H^p(\Omega,\F)$ einer
Steinschen Mannigfaltigkeit $\Omega$ mit Werten in
einer kohärenten analytischen Garbe $\F$  für
jedes $p>0$ verschwindet (siehe z.B. \cite{Hoer},
Theorem 7.4.3), ist  wesentlich allgemeiner und  abstrakter, umfaßt aber die
Lösbarkeit von Cousin-I-Verteilungen auf Steinschen Mannigfaltigkeiten.   

 Sei jetzt $\U=(U_i)_{i\in \N}$ eine
Überdeckung einer Steinschen Mannigfaltigkeit
 $\Omega$ mit pseudokonvexen relativ kompakten Gebieten mit
der Eigenschaft, daß eine feste Überdeckungsmenge von höchstens $M$ 
paarweise verschiedenen Überdeckungsmengen geschnitten wird. Dann ist der
Schnitt von $M+1$ paarweise verschiedenen Überdeckungsmengen leer.

In das Diagramm \auslabeln{15.1.1}
setzen wir für $F_{n,0}$, $n=0,\ldots,M+1$ jeweils den Kokettenraum 
$C^{M+1-n}(\U,\F)$ von
Schnitten über $U_{i_0,\ldots ,i_{M+1-n}}:=
U_{i_0}\cap\ldots\cap U_{i_{M+1-n}}$ mit Werten in $\F$ ein.
$d_{n,0}$ ist dann der Korandoperator $\delta$ der Stufe $M-n$.
$F_{M+2,0}$ sei $\Gamma(\Omega,\F)$ und $d_{M+1,0}:\Gamma(\Omega,\F)\to
C^{0}(\U,\F)$ definieren wir dann durch $f\mapsto (f\mid_{U_i})_{i\in \N}$.
Für $n\ge M+3$ setzen wir $F_{n,0}:=\{0\}$, wodurch
die obige Induktion nach endlich vielen Schritten abgeschlossen ist.

Einen  Kokettenraum von Koketten der Länge $\sss+1$ 
 können wir als Teilraum von $\prod_{i\in \N}$
 $\Gamma(U_{i_0,\ldots,i_\sss},$ $\F)$ auffassen und, wie wir im Kapitel 2 
sehen werden, mit
einer Fr\'echetraumtopologie versehen, die sich aus der Fr\'echetraumtopologie
der Räume $\Gamma(U_{i_0,\ldots,i_\sss},\F)$ herleiten läßt. Da wir stets
voraussetzen, daß $\F$ eine Untergarbe einer $p$-fachen Kopie von
$\oka$ ist, können wir auf $\Gamma(U_{i_0,\ldots,i_\sss},\F)$ jeweils die
Fr\'echetraumtopologie der kompakt gleichmäßigen Konvergenz 
auf einer ausschöpfenden Folge von kompakten Teilmengen von 
$U_{i_0,\ldots,i_\sss}$ für
$p$-tupel holomorpher Funktionen auf $U_{i_0,\ldots,i_\sss}$ induzieren.

In jedem Schritt des oben beschriebenen Induktionsverfahrens werden wir die
Stetigkeits- bzw. Urbildabschätzungen nachhalten, um zu einer 
Lösungsabschätzung für den $\delta$-Operator zu gelangen.

Die Hauptarbeit besteht nun darin, das Diagramm
\auslabeln{15.1.1} mit Fr\'echeträumen  und stetigen linearen Abbildungen so
aufzufüllen, daß das Diagramm kommutiert und die erforderlichen
Exaktheitsbedingungen der Zeilen und Spalten erfüllt sind. 
Ferner sind Stetigkeits- bzw. 
Urbildabschätzungen zu zeigen. Dies alles geschieht in den Kapiteln 
3 bis 6. 

Die Exaktheit der ersten Spalte in dem Diagramm unter Berücksichtigung
von Lösungsabschät\-zungen ist nichts anderes als das
angestrebte Theorem B mit Schranken.

Wir wollen mit $\Phi_\gamma(\U)$ die Menge der stetigen plurisubharmonischen
Funktionen bezeichnen, für die $\phi\ge 0$ und $\sup_{U_i}\phi\le
\gamma^{\frac{1}{M}} \inf_{U_i}\phi$ für jedes $i\in \N$ gilt.

Das in Kapitel 6 gezeigte Theorem B mit Schranken ergibt dann auf der 
nullten Stufe den folgenden Satz, den wir hier ohne die Verwendung
der Begriffe Kokette und Korandoperator formulieren wollen:

\begin{satzohnebeweis}{Satz}\einlabelnsatz{15.2a}
Sei $(K_i)_{i\in \N}$ eine Folge von kompakten Teilmengen mit $K_i\subset U_i$
für $i\in \N$ und $\psi$ eine nicht negative plurisubharmonische
Funktion auf $\Omega$. Dann gibt es eine plurisubharmonische Funktion
$\chi$ auf $\Omega$, so daß folgendes gilt:

Zu jeder Doppelfolge $(c_{i,j})_{i,j\in \N}$ mit $c_{i,j} \in 
\Gamma(U_i\cap U_j, \F)$, $c_{i,j}=-c_{j,i} $ und $c_{i,j}+
c_{j,k}+c_{k,i}=0$ für alle $i,j,k\in \N$, sowie jedem $\phi\in \Phi_\gamma(\U)$
gibt es eine Folge $(c'_i)_{i\in \N}$ mit $c_i\in \Gamma(U_i,\F)$, 
$c'_i-c'_j=c_{i,j}$ für jedes $i,j\in \N$ und eine Abschätzung
\[   \sum_{i\in \N} \sup_{K_i} e^{-\gamma\phi - \psi -\chi}
                 \sup_{K_i} \betrag{c'_i} \le 
     \sum_{i,j\in \N} \sup_{U_i\cap U_j}e^{-\phi}\sup_{U_i\cap U_j} e^{-\psi}
        \sup_{U_i\cap U_j} \betrag{c_{i,j}} \punkt   \]

\end{satzohnebeweis}

Das entscheidende an der gewonnenen Lösungsabschätzung ist, daß
der Gewichtverlust $\chi$ nicht von $\phi\in \Phi_\gamma(\U)$ abhängt.
Durch direkte Anwendung des Satzes von der offenen Abbildung auf
den Korandoperator kann dieses Ergebnis nicht gezeigt werden.

Wenn wir für  $\Omega$ den $\C^N$ und für $\F$ die Idealgarbe $J_V$ der Keime der
auf $V$ verschwindenden holomorphen Funktionen einsetzen, können
wir mit Hilfe von Satz \auslabelnsatz{15.2a} lokale Fortsetzungen 
einer auf $V$ holomorphen Funktion zu einer
ganzen Funktion zusammensetzen. Hierzu muß   $V$ nicht singularitätenfrei
sein. Da wir aber für die lokalen Fortsetzungen die Existenz
eines holomorphen Retraktes für $V$ verwenden wollen und damit
gewissermaßen optimale Daten für die Anwendung von Satz \auslabelnsatz{15.2a}
erhalten, haben wir gefordert, daß $V$ eine abgeschlossene 
Untermannigfaltigkeit des $\C^N$, und damit Steinsch ist.

Wir erhalten aus Satz \auslabelnsatz{15.2a} den in Kapitel 8 gezeigten 
Fortsetzungssatz für holomorphe Funktionen
auf einer abgeschlossenen Untermannigfaltigkeit $V$ des $\C^N$ 
( siehe Satz \auslabelnsatz{12.4}):

\begin{satzohnebeweis}{Satz}\einlabelnsatz{15.3}
Seien $\phi_1>0$ und $\phi_2>0$ stetige plurisubharmonische Funktionen auf dem 
$\C^N$. Sei ferner $\gamma>1$ fest. Dann gibt es zu jedem Kompaktum $K$ eine Konstante
$B_K$, so daß zu jedem $f\in H(V)$ und jedem $\alpha >0$ mit 
\formel{
   \log \betrag{f(z)} \le \phi_1(z) +\alpha \phi_2(z) \mbox{ für jedes }z\in V   }

eine ganze Funktion $F$ existiert mit 
\be
\item[a)] $F(z)=f(z)$ für jedes $z\in V$,
\item[b)] $\sup_K \betrag{ F(z)} \le e^{B_K} e^{\gamma \sup_K (\phi_1+ \alpha\phi_2)}$ 
          für jedes Kompaktum $K$.
\ee
Hierbei hängt $B_K$ zwar von $\phi_1$ und $\phi_2$ ab, aber nicht von $\alpha$.
\end{satzohnebeweis}

Aus Satz \auslabelnsatz{15.3} folgt nun leicht in Anlehnung an die
Argumentation in \cite{M-H} die Existenz eines Ausdehnungsoperators, wenn
man eine Basisfolge in $H(V)$ zur Verfügung hat, die von  der
vorausgesetzten
Isomorphie von $H(V)$ zu $H(\C^d)$ induziert wird.

Wir wollen hier wie angekündigt als Anwendung von Satz \auslabelnsatz{15.3} die
Bedingungen untersuchen, unter denen der Ausdehnungsoperator darüberhinaus 
linear zahmen Abschätzungen genügt.

Eine stetig lineare Abbildung $T:(E,\normleer_{E,n})\to 
(F,\normleer_{F,n})$ zwischen den gradierten Fr\'echet\-räumen
$E$ und $F$ wird als linear zahm bezeichnet, falls es $a,b\in \N_0$ mit $a\not= 0$
gibt, so daß zu jedem $n$ ein $C_n$ existiert mit
\[   \norm{n}{Tx} \le C_n \norm{an+b}{x}  \foralls x\in E  \punkt\]

$T$ heißt zahm, falls $a=1$ gewählt werden kann. Gradierungen in 
einem Fr\'echetraum heißen linear zahm äquivalent, falls die
Identität in der jeweiligen Richtung linear zahm ist.

Eine linear zahme Isomorphie von $H(V)$ und $\Lambda_\infty(k^{\frac{1}{d}})$
ist leicht durch Eigenschaften einer Basisfolge $(f_k)_k$ in $H(V)$ zu
charakterisieren. Ist nämlich $H(V)$ zu $\Lambda_\infty(k^{\frac{1}{d}})$
linear zahm isomorph, so existiert eine Basisfolge $(f_k)_k$ in $H(V)$ mit
folgenden Eigenschaften:
\be
\item[B1)] Es existieren $a,b\in \N_0$, $a\ge 1$, so daß es zu jedem
$n\in \N$ ein $C\ge 1$ gibt mit
\[   \norm{n}{f_k}\le C e^{(an+b)\, k^{\frac{1}{d}}} \foralls k  \punkt\]
\item[B2)] Es existieren $a,b\in \N_0$, $a\ge 1$, so daß es zu jedem
$n\in \N$ ein $C\ge 1$ gibt mit
\[   \betrag{\lambda_k}\le \norm{an+b}{f} e^{-n\, k^{\frac{1}{d}}} \]
für alle $f=\sum_k \lambda_k f_k\in H(V)$ und alle $k$.
\ee

Besitzt andererseits $H(V)$ eine Basisfolge $(f_k)_k$ mit B1) und B2), dann
ist $H(V)$ zu $\Lambda_\infty(k^{\frac{1}{d}})$
linear zahm isomorph, denn wegen B1) ist 
\[   \norm{n}{\sum_k \lambda_k f_k}\le C\sum_k \betrag{\lambda_k} 
      e^{(an+b)\, k^{\frac{1}{d}}}  \] 
und wegen B2) ist 
\[  \sum_k \betrag{\lambda_k} e^{n\,k^{\frac{1}{d}}}=
     \sum_k \betrag{\lambda_k} e^{(n+1)\,k^{\frac{1}{d}}} e^{-k^{\frac{1}{d}}} 
    \le C \norm{a(n+1) +b}{f} \sum_k e^{-k^{\frac{1}{d}}} = C'\norm{a(n+1) +b}{f}\punkt \]

Für Potenzreihenräume ist der Begriff linear zahm  von besonderer
Bedeutung, da zueinander isomorphe Potenzreihenräume 
$\Lambda_\infty(a_k)$ bzw.  $\Lambda_\infty(b_k)$ stets
linear zahm isomorph sind, vermöge der Identität. Andererseits ist zum Beispiel
bei nuklearen Potenzreihenräumen eine zahme Isomorphie
sehr viel seltener zu erwarten, da dies $\lim \frac{a_k}{b_k}=
\lim \frac{b_k}{a_k}=1$ impliziert. Wir bemerken hierzu zunächst,
daß der nukleare Potenzreihenraum $\Lambda_\infty(c_k)$ stets per 
Identität zahm isomorph zu 
\[  \{ x\in \C^{\N_0}\doppelpunkt \norm{n}{x}:=\sum_{k=0}^\infty \betrag{x_k}^2
        e^{2n\,c_k} <\infty \mbox{ für alle } n\in \N \}  \]
ist. Dann erhalten wir dem Beweis zu Proposition 29.1 in \cite{M-Venglisch} folgend
zu jedem $\theta > 1$ eine Konstante $D\ge 1$, so daß 
$a_k\le \theta b_k + D$ für jedes $k\in \N$ ist und umgekehrt.  

Wir wollen durch die Suprema auf folgenden Polyzylindern des $\C^N$
\[  \Delta_n:= \{ z\in \C^N \doppelpunkt \betrag{z_j}\le e^n 
    \foralls 0\le j\le N \} \]
die Standardgradierung $(\normleer_n)_n$ für $H(\C^N)$ festlegen.

Auf $V$ sei die Standardgradierung $(\normleer_{H(V),n})_n$ gerade durch die Spuren 
$\Delta_n\cap V$ gegeben. Die Gradierung hängt also
von der Einbettung von $V$ in den $\C^N$ ab. Ohne Einschränkung
können wir $0\in V$ annehmen.

Wenn wir die Polyzylinder $\Delta_n$, $n\in \N$, durch
$B_n:=\{z\in \C^N\doppelpunkt \betrag{z}< e^n\}$ ersetzen, erhalten wir
auf $H(\C^N)$ bzw. auf $H(V)$ eine zahm äquivalente Gradierung. Solche
zahm äquivalente Gradierungen fassen wir auch als Standardgradierungen für
$H(\C^N)$ bzw.  $H(V)$ auf. 

Auf $H(V)$ ist wegen der Surjektivität der Einschränkungsabbildung
$H(\C^N)\to H(V)$ ein Quotientenhalbnormensystem durch
\[  \norm{n}{f}^q:=\inf\{ \norm{n}{F}\doppelpunkt F\in H(\C^N)\komma
    F\mid_V=f \}  \]
 gegeben. Es gilt natürlich $\normleer_n\le \normleer_n^q$ auf 
$H(V)$ für jedes $n$. Da der Satz von der offenen Abbildung gilt,
ist $(\normleer_n^q)_n$ ein zu $(\normleer_{H(V),n})_n$ äquivalentes
Halbnormensystem.

Die Monome in $H(\C^N)$ besitzen bezüglich der Standardgradierung
die Eigenschaften B1) und B2) und zwar jeweils mit $a=1$ und $b=0$.
Die Eigenschaft B2) ist äquivalent zu den Cauchyschen Ungleichungen
für ganze Funktionen auf den Polyzylindern $\Delta_n$ (siehe z.B.
\cite{Hoer}, Theorem 2.2.7). Wir haben also gezeigt, daß $H(\C^N)$ sogar
zahm isomorph zu einem Potenzreihenraum $\Lambda_\infty(a_k)$ ist, falls
$a_k$ gerade die Summe der Exponenten des $k$-ten Monoms ist.

$H(\C^N)$ ist daher auch linear zahm isomorph zu $\Lambda_\infty(k^{\frac{1}{N}})$,
aber nicht zahm isomorph, falls $N\ge 2$. Um das einzusehen, sei $(z_l)_k$ eine
Abzählung der Monome $z^{\alpha_1}_1\cdot\ldots\cdot z^{\alpha_N}_N$,
wobei $\alpha_1,\ldots,\alpha_N\in \N_0$, mit der Eigenschaft, daß für
$s(z^{\alpha_1}_1\cdot\ldots\cdot z^{\alpha_N}_N):=\alpha_1+\ldots+\alpha_N$ gilt
$s(z_l)\le s(z_{k+1})$.

Man zeigt leicht, daß die Anzahl der Monome $z_l$ mit $s(z_l)\le m$
jeweils $\binom{N+m}{m}$ beträgt. Es gilt daher
\[   \frac{m}{\binom{N+m}{m}^{\frac{1}{N}}} 
     \le \frac{s(z_l)}{k^{\frac{1}{N}}} 
     \le \frac{m}{\left(\binom{N+(m-1)}{m-1} +1\right)^{\frac{1}{N}}} \]
für alle $k$ mit $s(z_l)=m$. Die linke und die rechte Seite der
Ungleichung haben den gleichen Grenzwert  für $m\to \infty$. Wenn
man den linken Ausdruck ausmultipliziert, sieht man sofort, daß
dieser gegen die $N$-te Wurzel aus $N!$ konvergiert für $m\to \infty$.

Für den angekündigten Charakterisierungssatz benötigen wir
noch einige Begriffe und Definitionen.

Für einen gradierten Fr\'echetraum $(E,\normleer_n)$ wollen wir
die Eigenschaften (DN) und $(\Omega)$ verschär\-fen:

Wir sagen, daß $(E,\normleer_n)$ (DN) in Standardform besitzt, 
falls für alle $n\ge 2$  gilt
\[   \normleer_n^2 \le C_n \normleer_{n-1} \normleer_{n+1}  \punkt  \]

$(E,\normleer_n)$ besitze $(\Omega)$ in Standardform, falls für
alle $n\ge 2$ gilt
\[   \normleer_n^{\ast 2} \le D_n \normleer_{n-1}^{\ast} \normleer_{n+1}^{\ast}\]

Aufgrund der Cauchy-Schwarzschen Ungleichung besitzt $\Lambda_\infty(a_k)$
stets (DN) und $(\Omega)$ in Standardform. Die Konstanten sind alle gleich 1.

Von jetzt an wollen wir annehmen, daß $V$ als analytische 
Untervarietät des $\C^N$
nur aus einer irreduziblen Komponente besteht. Dies ist genau dann der Fall,
wenn $V$ zusammenhängend ist. Da nach Satz \auslabelnsatz{15.2a} die
Existenz eines Ausdehnungsoperators notwendig die Eigenschaft (DN) für
$H(V)$ voraussetzt,
besitzt $H(V)$ mindestens eine stetige Norm. Daher impliziert die 
Existenz eines Ausdehungsoperators, daß $V$ ohnehin  nur endlich viele
Zusammenhangskomponenten besitzen kann.

$V$ heißt algebraisch, falls es eine endliche Menge S  von Polynomen im $\C^N$
gibt, so daß $V=\{ z\in \C^N\doppelpunkt p(z)=0 \foralls p\in S\}$.
Die Anzahl der Polynome in $S$ kann man auf $N$ begrenzen (siehe z.B.
\cite{Kunz}, Chapter V.1).

Das nun folgende  Kriterium für eine Teilmenge einer
analytischen Untervarietät des $\C^N$, eine Teilmenge einer algebraischen
Varietät zu sein, stammt
von A. Sadullaev (\cite{S1}, Theorem 2.2). 

Sei dazu $L$ die Menge aller plurisubharmonischen Funktionen auf dem $\C^N$, 
für die es ein $\alpha(u)\in \R$ gibt, so daß $u(z)\le \log(1+\norm{}{z})
+\alpha(u)$ auf dem $\C^N$. Ist $K\subset \C^N$ ein Kompaktum, dann sei 
 \[  S(z,K):= \sup\{u(z)\doppelpunkt u\in L\komma u\mid_K\le 0\} \punkt \]
Die obere Regularisierte von $S(z,K)$ heißt Siciaksche 
Extremalfunktion bezüglich $K$.

Das Kriterium besagt nun folgendes:

Eine Teilmenge $A$ einer analytischen Menge im $\C^N$ ist genau dann eine Teilmenge
einer algebraischen Varietät, falls es ein Kompaktum in $A$ gibt,
so daß $S(z,K)$ lokal wesentlich beschränkt ist in $A$.

A. Aytuna verwendet in \cite{Ayt3} dieses Kriterium, um für eine irreduzible Komponente $V$
einer analytischen Varietät im $\C^N$ (die also auch singuläre Punkte
haben darf) zu zeigen, daß $V$ genau dann algebraisch ist, wenn
ein linear zahmer Ausdehnungsoperator $H(V)\to H(\C^N)$ bezüglich der
Standardgradierungen existiert (\cite{Ayt3}, Theorem 2).


Wir werden das Kriterium von Sadullaev in Kapitel \ref{Kap10} dazu verwenden,
ein Kriterium für die Existenz eines linear zahmen Ausdehnungsoperators abzuleiten,
das das Wachstum plurisubharmonischer Funktionen auf $V$ relativ zu einer 
Standardgradierung beschreibt.

Ist $\phi $ plurisubharmonisch auf $V$, so sei für $r\in \R^+_0$
\[   m_{\phi}(r):= \sup\{ \phi(z)\doppelpunkt z\in V, \betrag{z}<e^r \} \punkt \]
Das Kriterium ist erfüllt, falls es $a\in \N$ und $b\in \N_0$ gibt, so daß für
jedes auf $V$ plurisubharmonische $\phi $, für jedes $\lambda\in [0,1]$ und
alle $r_1,r_2\in \R^+_0$ gilt
\[   m_{\phi}(\lambda r_1 + (1-\lambda) r_2)\le \lambda m_{\phi}(ar_1 +b) + 
         (1-\lambda)m_{\phi}(ar_2 +b)  \punkt \]

Diese Eigenschaft ist stärker als die Stark-Liouvillesche Eigenschaft (Satz
\auslabelnsatz{15.2}, 3)), denn ist $\phi$ plurisubharmonisch und
 beschränkt auf $V$, dann gilt für $r>b$
\[   m_{\phi }(r)=m_{\phi }(\frac{n-1}{n}0 + \frac{1}{n} n r)\le 
     \frac{n-1}{n}m_{\phi }(b) + \frac{1}{n}m_{\phi }(anr+b) \stackrel{n\to\infty}{\to}
     m_{\phi }(b)\punkt  \]
Wegen des Maximumprinzips für  plurisubharmonische Funktionen auf $V$ (siehe
z.B. \cite{S1}, 1.3) ist $\phi$ auf jedem Kompaktum, das in $\{\betrag{z}< e^r\}$
enthalten ist, konstant. Damit ist $\phi$ auf ganz $V$ konstant.

Wir fassen nun alle Kriterien zusammen:

Für eine zusammenhängende abgeschlossene  Untermannigfaltigkeit $V$
des $\C^N$ der Dimension $d$ gilt der folgende Satz:

\begin{satz}{Satz}\einlabelnsatz{15.4}
Äquivalent sind:
\be
\item[1)] Es gibt einen linear zahmen Ausdehnungsoperator 
  \[    E: (H(V),\normleer_{H(V),n})\to (H(\C^N),\normleer_{H(\C^N),n})   \]
bezüglich der durch die Polyzylinder $(\Delta_n)_n$ gegebenen Standardgradierungen
auf $H(V)$ bzw. $H(\C^N)$.
\item[2)] $(H(V),\normleer_{H(V),n})$ ist linear zahm isomorph zu einem
nuklearen gradierten Fr\'echetraum $(F,\normleer_n)$, der (DN) in 
Standardform besitzt, und die Quotientenhalbnormen $(\normleer_n^q)_n$ auf
$H(V)$ sind
linear zahm äquivalent zu $(\normleer_{H(V),n})_n$.
\item[3)] Es gibt ein Kompaktum $K$ in $V$, so daß $S(z,K)$ in $V$
lokal wesentlich beschränkt ist.
\item[4)] Es gibt $a\in \N$ und $b\in \N_0$, so daß für 
jede auf $V$ plurisubharmonische Funktion $\phi$, jedes $\lambda\in [0,1]$ und 
alle $r_1,r_2\in \R_0^+$ gilt mit
$m_{\phi}(r):= \sup\{ \phi(z)\doppelpunkt z\in V, \betrag{z}<e^r \}$:
\[  m_{\phi}(\lambda r_1 + (1-\lambda) r_2)\le \lambda m_{\phi}(ar_1 +b) + 
         (1-\lambda)m_{\phi}(ar_2 +b)  \punkt \] 
 \item[5)] $(H(V),\normleer_{H(V),n})$ ist linear zahm isomorph zu 
$\Lambda_\infty(k^{\frac{1}{d}})$.
\item[6)] $(H(V),\normleer_{H(V),n})$ ist linear zahm isomorph zu
$(H(\C^d), \normleer_{H(\C^d),n})$.
\item[7)] In $H(V)$ gibt es eine Basisfolge $(f_k)_k$ mit B1) und B2).
\item[8)] $V$ ist algebraisch.
\ee
\end{satz}

1)$\Rightarrow$2) ist schlicht, denn das Bild von $E$ ist linear zahm 
isomorph zu einem projizierten Teilraum von $\Lambda_\infty(k^{\frac{1}{N}})$,
der darum (DN) in Standardform besitzt. Zum anderen ist für jedes
$f\in H(V)$ die Quotientenhalbnorm $\norm{n}{f}^q\le \norm{H(\C^N),n}{Ef}$.
Also folgt auch die zweite Bedingung von 2) aus 1).

2)$\Rightarrow$5) ist ein Ergebnis von D. Vogt (\cite{V2}, Theorem 2.3), das
wir auf die wegen 2) linear zahme Sequenz $\Lambda_\infty(k^{\frac{1}{N}})
\to F\to 0$ anwenden. Wir erhalten eine linear zahme Isomorphie von
$F$ zu einem Potenzreihenraum $\Lambda_\infty(a_k)$, woraus 5) folgt.

5)$\Rightarrow$1) ist der Spezialfall $\phi(z):=\log \betrag{z}$ von 
Satz \auslabelnsatz{13.2}. Wir zeigen in Kapitel 9, daß Satz
\auslabelnsatz{13.2} leicht aus dem Fortsetzungssatz \auslabelnsatz{15.3} folgt.

Die Äquivalenzen 3)$\Leftrightarrow$8) und 1)$\Leftrightarrow$8) sind
Spezialfälle der oben zitierten Ergebnisse von A. Sadullaev
bzw. A. Aytuna. Die Äquivalenz 5)$\Leftrightarrow$6)$\Leftrightarrow$7)
sind nach den obigen Ausführungen klar. 

1)$\Rightarrow$4) folgt aus Satz \auslabelnsatz{16.4}. 4)$\Rightarrow$8) folgt
aus Satz \auslabelnsatz{16.5} \Ende

Die Eigenschaft 4) aus Satz \auslabelnsatz{15.4} läßt sich in einer Weise 
verallgemeinern,
die unabhängig von der Einbettung der Steinschen Mannigfaltigkeit in einen
$\C^N$ ist:

Eine Steinsche Mannigfaltigkeit $V$ habe die Eigenschaft (LK), falls eine 
plurisubharmonische Ausschöpfungsfunktion $\psi$, sowie $a\in \N$, $b\in \N_0$
existieren, so daß für jede plurisubharmonische Funktion $\phi$ auf 
$V$ und alle $r_1,r_2\in \R_0^{+}$ und jedes $\lambda \in [0,1]$ gilt
\[   m_{\phi}^{\psi}(\lambda r_1 + (1-\lambda)r_2) \le 
       \lambda m_{\phi}^{\psi}(ar_1 + b) + (1-\lambda) m_{\phi}^{\psi}(ar_2 + b)\komma \]
wobei
\[   m_{\phi}^{\psi}(r):= \sup\{ \phi(z)\doppelpunkt z\in V\komma \psi(z)<r \} 
     \mbox{ ist für jeder }r\in \R \punkt  \]

Eine in den $\C^N$ eingebettete abgeschlossene zusammenhängende Untermannigfaltigkeit
ist also genau dann algebraisch nach Satz \auslabelnsatz{15.4}, falls sie (LK)
mit $\psi(z)=ln\betrag{z}$ für $z\in V$ besitzt.

Ferner impliziert (LK) die Stark-Liouvillesche Eigenschaft 
(Satz \auslabelnsatz{15.3}, 3)) für zusammenhängende Steinsche Mannigfaltigkeiten.

Hieraus leiten sich zwei Fragen ab:
\begin{enumerate}
\item[1)] Ist jede Steinsche Mannigfaltigkeit mit (LK) biholomorph äqui\-valent zu 
 einer algebraischen Varietät in einem $\C^N$?
\item[2)] Ist die Stark-Liouvillesche Eigenschaft auf einer zusammenhängenden
Mannigfaltigkeit äqui\-valent zur Eigenschaft (LK)?
\end{enumerate}

Wenn und nur wenn beide Fragen eine positive Antwort haben, wäre jede 
zusammenhängende Steinsche Mannigfaltigkeit $V$, für die $H(V)$ die Eigenschaft
(DN) besitzt, isomorph  zu einer algebraischen Varietät eines $\C^N$.

\vspace{4cm}

An dieser Stelle möchte ich mich bei den Herren Dr. L. Frerick und Dr. M. Tidten
sehr herzlich bedanken für die Bereitschaft zu vielen Anregungen und Gesprächen,
die das Entstehen dieser Arbeit begleitet haben. Mein ganz besonderer Dank gilt
meinem Doktorvater Herrn Prof. Dr. D. Vogt, der durch seine geduldige und 
stetige Anregung und Unterstützung entscheidend zum Gelingen dieser Arbeit 
beigetragen hat.

\newpage   
\section[{\sc Verallgemeinerter Lösungssatz für den Korandoperator}]{Der verallgemeinerte Lösungssatz für den Korandoperator $\delta$}
\label{Kap2} 

Wir führen in diesem Kapitel die in Kapitel 1 dargestellte Induktion
in dem Diagramm \auslabeln{15.1.1} durch, indem wir das Diagramm mit 
speziellen Fr\'echeträumen und linearen stetigen Abbildungen auffüllen.
Dabei ist jede Spaltenabbildung ein Korandoperator.
Wir behandeln die Kommutativität des Diagramms und die Exaktheitsbedingungen
an die Zeilen und Spalten  als Eigenschaften der Familie der vorkommenden
Fr\'echeträume. Wir zeigen, daß sich die Urbildabschätzungen
der exakten Spalten auf die erste Spalte überträgt.

In den folgenden Kapiteln muß dann gezeigt werden, daß es ein
solches Diagramm von Fr\'echet\-räumen gibt, wo in der ersten Spalte die Sequenz
der Korandabbildungen zwischen den
Kokettenräumen
mit Werten in einer lokal endlich erzeugten Untergarbe von
$\oka^p$ steht.

Sei $I$ eine abz"ahlbare Indexmenge. Sei f"ur jedes $k\in \N_0$, $\sss\in \N_0$ eine
Familie  $S_{k,\sss}$ von Fr\'echetr"aumen gegeben, so daß zu jedem Multiindex 
$\alpha=\ffg{0}{} \in I^{\sss +1}$, paarweise verschiedener Indizes aus $I$ ein 
Fr\'echetraum $F_\alpha^k$ mit einem jeweiligen aufsteigenden 
Halbnormensystem $(\normoben{k}{
\alpha,n}{\leer{0.5}})_{n\in \N}$ existiert. Es gelte $F_{\pi(\ffg{0}{})}^k=
F_{\ffg{0}{}}^k$ f"ur jede Permutation $\pi$ der Indizes $\ffg{0}{}$. Für die
lokalen Banachräume von $F_{\ffg{0}{}}^k$ bezüglich $\normoben{k}{
\ffg{0}{},n}{\leer{0.5}}$ sollen die  kanonischen Verbindungsabbildungen injektiv sein und
dichtes Bild besitzen.

Ferner sei f"ur jeden Multiindex $\ffg{1}{}\in I^\sss$ die Anzahl der Fr\'echetr"aume
$F_{i,\ffg{1}{}}^k \in S_{k,\sss}$, $i\in I$, die nicht der Nullraum sind,
gleichm"aßig f"ur jedes $k\in \N_0$ und $\sss\in \N_0$ durch eine Schranke $M\in \N$
beschr"ankt. Die Menge solcher Indizes $i$ zu einem Multiindex $\ffg{1}{}$ sei
f"ur jedes $k\in \N_0$ gleich. Wir bezeichnen sie mit $I_\ffg{1}{} $. Wir verlangen ferner, 
daß $I_{j,i_1,\ldots,i_\sss} \subset I_{i_1,\ldots,i_\sss}$ für jedes $j\in I$ und
jedes $i_1,\ldots,i_\sss\in I^\sss$ ist.

Ist $\alpha=\ffg{0}{}\in I^{\sss +1}$ ein Multiindex, dann sei mit $\alpha_j$ derjenige 
Multiindex in $I^\sss$ bezeichnet, der durch Streichung des Index $i_j$ entsteht.

Wir f"uhren jetzt den Begriff einer Kokette auf $S_{k,\sss}$, $k\in N_0$, $\sss\in \N_0$ ein:

\begin{satzohnebeweis}{Definition}

Eine Kokette $c=(c_\ffg{0}{})_{\ffg{0}{}\in I^{\sss +1}}$ der L"ange $\sss$ auf der Stufe $k$ 
sei ein Element
aus $\leer{-1.55}\prod\limits_{\ffg{0}{}\in I^{\sss +1}} \leer{-1}F_\ffg{0}{}^k$ 
mit der folgenden alternierenden Eigenschaft:
\[   c_\ffg{0}{}=(-1)^{\mbox{sign }\pi} c_{\pi(\ffg{0}{})} \komma \]
wobei $\pi$ eine Permutation der Indizes $\ffg{0}{}$ ist und sign $\pi$ die minimale Anzahl der 
Transpositionen angibt, in die $\pi$ zerlegt werden kann.
\end{satzohnebeweis}

Man beachte, daß $c_\ffg{0}{}$ und $c_{\pi(\ffg{0}{})}$ Elemente desselben Fr\'echetraumes sind.

Die Menge aller Koketten der L"ange $\sss$ sei mit ${\cal C}(S_{k,\sss})$ bezeichnet.

Wir wollen nun R"aume gewichteter Koketten einf"uhren und diese mit einer Fr\'echetraumtopologie
versehen.

Ist $J$ eine abz"ahlbare Indexmenge und $(F_j)_{j\in J}$ eine Familie von Fr\'echetr"aumen mit
aufsteigenden Halbnormensystemen $(\norm{j,n}{\leer{0.5}})_{n\in \N}$ jeweils in dem 
Fr\'echetraum $F_j$, dann sei 
\[   l_1((F_j)_j):=\{ x=(x_j)_{j\in J}\in \prod_{j\in J} F_j \doppelpunkt 
                \dreinorm{}{n}{x} := \sum_{j\in J}\norm{j,n}{x_j} <\infty\mbox{ f"ur jedes }n\}  \]
der Raum der bezüglich $(F_j)_j$ absolutsummierbaren Folgen.
$l_1((F_j)_j)$ ist also ein linearer Teilraum von $\prod\limits_j F_j$. 
Man zeigt wie im skalaren Fall,
daß $l_1((F_j)_j)$ versehen mit der durch $(\dreinorm{}{n}{\leer{0.5}})_n $
induzierten Topologie ein Fr\'echetraum ist.

Sei jetzt $C=(C_\alpha)_{\alpha\in I^{\sss +1}}$ eine Familie von  Konstanten aus dem halboffenen
Intervall $(0,1]$ mit 
$C_\alpha=C_{\pi(\alpha)}$ f"ur jede Permutation $\pi$ des Multiindex $\alpha$. Dann sei
\[    {\cal C}( S_{k,\sss},C) := \{ c=(c_\alpha)_{\alpha\in I^{\sss +1}} \in {\cal C}(S_{k,\sss})
         \doppelpunkt C\cdot c:=(C_\alpha\cdot c_\alpha )_{\alpha\in I^{\sss +1}} \in l_1(S_{k,\sss}) \}
         \punkt \]
F"ur $c\in {\cal C}( S_{k,\sss},C)$ gilt also
\[    \dreinorm{k}{n,C}{c} := \sum_{\alpha\in I^{\sss +1}} C_\alpha \normoben{k}{\alpha,n}{c_\alpha} <\infty \]
f"ur jedes $n\in \N$. Dasjenige Konstantensystem, das nur aus Einsen besteht, bezeichnen wir
mit $1$. 

${\cal C}( S_{k,\sss},C)$ ist versehen mit $(\dreinorm{k}{n,C}{\leer{0.5}})_n$ ein Fr\'echetraum,
denn ${\cal C}( S_{k,\sss},C)= l_1((F^k_\alpha)_{\alpha\in I^{\sss +1}})$, wobei die Halbnormen in
den $F^k_\alpha$ jeweils mit $C_\alpha$ skaliert sind. 

Ferner sei 
\[   {\cal C}_n(S_{k,\sss},C) := \{ c\in {\cal C}( S_{k,\sss})\doppelpunkt \dreinorm{k}{n,C}{c}
        < \infty \} \punkt  \]

Um einen Korandoperator f"ur Koketten auf der Stufe $k$ einf"uhren zu k"onnen, 
ben"otigen wir f"ur jedes 
$\alpha\in I^{\sss +1}$, $\sss\ge 1$ und $j=0,\ldots, \sss$ lineare stetige  Abbildungen
\[    \iota_{\alpha_j}^{k,\alpha} : F_{\alpha_j}^k \rightarrow F_\alpha^k   \] 
mit folgenden Eigenschaften:

\formel{{\einlabeln{10.1.1}}
         &\mbox{1)}& \normoben{k}{\alpha,n}{\iota_{\alpha_j}^{k,\alpha} x} \le \normoben{k}{\alpha_j,n}{x}
                       \mbox{ f"ur jedes } n\in \N  \\
         &\mbox{2)}&  \iota_{\alpha_j}^{k,\alpha}\circ \iota_{(\alpha_j)_i}^{k,\alpha_j} = 
                        \iota_{\alpha_i}^{k,\alpha}\circ \iota_{(\alpha_i)_{j-1}}^{k,\alpha_i} 
                         \mbox{ f"ur } 1\le i<j\le \sss+1  \\
         &\mbox{3)}&  \iota_{\alpha_j}^{k,\alpha}\circ \iota_{(\alpha_j)_i}^{k,\alpha_j} =
                        \iota_{\alpha_{i+1}}^{k,\alpha}\circ \iota_{(\alpha_{i+1})_j}^{k,\alpha_{i+1}} 
                         \mbox{ f"ur } 0\le j\le i \le \sss  \\
         &\mbox{4)}&   \iota_{\alpha_j}^{k,\alpha} \mbox{ ist injektiv, falls } F^k_\alpha \not=
                       \{0\} \mbox{ ist. }
}

Falls klar ist, auf welcher Stufe $k$ sich die Abbildungen befinden, lassen wir den Index $k$ fort.
Den Korandoperator auf der Stufe $k$ definieren wir nun wie folgt:

Für $\sss\ge 0$ sei
\formelohne{  \delta^\sss: {\cal C}( S_{k,\sss}) \ni (c_\alpha)_{\alpha\in I^{\sss +1}} 
            &\mapsto & ((\delta^\sss c)_\beta)_{\beta\in I^{\sss +2}} \in {\cal C}( S_{k,\sss+1}) \\
                     \mbox{ mit }  (\delta^\sss c)_\beta
            &:=&  \sum_{j=0}^{\sss+1} (-1)^j \iota_{\beta_j}^\beta c_{\beta_j}   }

Zwischen zwei Konstantensystemen $C=(C_\alpha)_{\alpha\in I^{\sss +1}}$ und $D=(D_\beta)_{\beta\in I^{\sss +2}}$
sei die Vergleichsrelation $C\succ D$ genau dann gegeben, falls 
\[   C_\alpha \ge \max_{i\in I_\alpha} D_{i,\alpha} \mbox{ f"ur alle } \alpha\in I^{\sss +1}\punkt \]

F"ur den Korandoperator gilt das folgende 

\begin{satz}{Lemma}\einlabelnsatz{10.2}
Falls $C\succ D$, dann ist $\delta^\sss$ ein stetiger linearer Operator
${\cal C}_n( S_{k,\sss},C) \rightarrow {\cal C}_n( S_{k,\sss+1},D)$. Ferner gilt 
\be
\item[1)] $ \delta^\sss \circ \delta^{\sss-1} =0$ f"ur $\sss=1,2,\ldots$.
\item[2)] $ \dreinorm{k}{n}{\delta^\sss c} \le M\cdot (\sss+2)\cdot \dreinorm{k}{n}{c}$ für
jedes $c\in{\cal C}_n( S_{k,\sss},C)$.
\ee
\end{satz}

Wir zeigen zun"achst $ \delta^\sss \circ \delta^{\sss-1}c =0$ f"ur ein $c\in {\cal C}( S_{k,\sss-1})$.
F"ur $\alpha\in I^{\sss +2}$ gilt:
\formelohne{  (\delta^\sss(\delta^{\sss-1} c))_\alpha 
     &=&  \sum_{j=0}^{\sss+1} (-1)^j \iota_{\alpha_j}^\alpha \left( 
                          \sum_{i=0}^\sss (-1)^i \iota_{(\alpha_j)_i}^{\alpha_j} c_{(\alpha_j)_i}\right) \\
     &=& \leer{-0.5}\sum_{0\le i<j\le \sss+1}\leer{-0.5} (-1)^{i+j} \iota_{\alpha_j}^\alpha\circ
                \iota_{(\alpha_j)_i}^{\alpha_j} c_{(\alpha_j)_i} + 
               \leer{-0.5}\sum_{0\le j\le i \le \sss}\leer{-0.5} (-1)^{i+j} \iota_{\alpha_j}^\alpha\circ
                \iota_{(\alpha_j)_i}^{\alpha_j} c_{(\alpha_j)_i}    }
Wir wenden die Kommutationseigenschaften \auslabeln{10.1.1} an und erhalten mit $n:=j$ und $m:=i+1$
\formelohne{   \leer{-0.5}\sum_{0\le j\le i \le \sss}\leer{-0.5} (-1)^{i+j} \iota_{\alpha_j}^\alpha\circ
                \iota_{(\alpha_j)_i}^{\alpha_j} c_{(\alpha_j)_i}
          &=&  \leer{-0.5}\sum_{0\le j\le i \le \sss}\leer{-0.5} (-1)^{i+j} \iota_{\alpha_{i+1}}^\alpha
                \circ \iota_{(\alpha_{i+1})_j}^{\alpha_{i+1}} c_{(\alpha_{i+1})_j} \\
          &=&  (-1)\cdot\leer{-1}\sum_{0\le n<m\le \sss+1}\leer{-0.5} (-1)^{m+n} \iota_{\alpha_m}^\alpha\circ
                \iota_{(\alpha_m)_n}^{\alpha_m} c_{(\alpha_m)_n} \punkt  }
Also ist $(\delta^\sss(\delta^{\sss-1} c))_\alpha = 0$.

Wir zeigen nun die Stetigkeitsabsch"atzung:

\formelohne{ \dreinorm{k}{n}{\delta^\sss c} 
          &=& \sum_{\alpha\in I^{\sss +2}} D_\alpha \normoben{k}{\alpha,n}{\sum_{j=0}^{\sss+1} (-1)^j 
                       \iota_{\alpha_j}^\alpha c_{\alpha_j} }           \\
          &\le& \sum_{\alpha\in I^{\sss +2}} D_\alpha \normoben{k}{\alpha,n}{\sum_{j=1}^{\sss+1} (-1)^j 
                       \iota_{\alpha_j}^\alpha c_{\alpha_j} }  + 
                \sum_{\alpha\in I^{\sss +2}} D_\alpha \normoben{k}{\alpha,n}{c_{\alpha_0}} \punkt   }
F"ur den zweiten Summanden der rechten Seite gilt:
\formelohne{ 
      \sum_{\alpha\in I^{\sss +2}} D_\alpha \normoben{k}{\alpha,n}{c_{\alpha_0}}  
      &=&  \sum_{\beta\in I^{\sss +1} } \sum_{i\in I_\beta^k} D_{i,\beta} \normoben{k}{\beta,n}{c_\beta} \\
      &\le& M  \sum_{\beta\in I^{\sss +1} } C_\beta \normoben{k}{\beta,n}{c_\beta} \\
      &=& M  \dreinorm{k}{n}{c} \punkt }
Durch iteratives Abspalten der weiteren Summanden erhalten wir schließlich 
\[   \dreinorm{k}{n}{\delta^\sss c} \le (\sss+2) M \dreinorm{k}{n}{c}  \punkt \]
\hfill $\Box$

Wir bemerken noch, daß f"ur jeden Multiindex $\alpha\in I^{M+2}$ mit paarweise verschiedenen
Indizes aus $I$ der Raum $F_\alpha^k$ stets der Nullraum ist f"ur alle $k$. Denn es gibt ein
$j\in \{1,\ldots,M+1\}$, so daß $F_{i_0,i_j}^k=\{ 0 \}$ ist, sonst w"are die Anzahl
der Indizes in $I_{i_0}$ gr"oßer als $M$, was der Definition von $M$
widerspräche. Ohne Einschr"ankung sei
$j=1$. Also ist $I_{i_0,i_1}=\emptyset$. Daher ist $I_{i_0,i_1,\ldots,i_M}=\emptyset$,
woraus  $F_{i_0,\ldots,i_{M+1}}^k=\{0\}$ folgt.

Im folgenden wollen wir unter einem System $(S_{k,\sss})_{k,\sss}$ stets ein System von
Fr\'echetr"aumen wie oben verstehen, das einen wohldefinierten Korandoperator
zul"aßt.

Ein Beispiel f"ur ein solches System $S_{k,\sss}$ mit festem $k$ und $\sss=0$  
wird durch die R"aume 
holomorpher Funktionen auf jeweils einem Element einer offenen zusammenhängenden
 "Uberdeckung $(U_i)_{i\in I}$
einer komplexen Mannigfaltigkeit gegeben, wobei sich nicht mehr als eine feste Anzahl
von Elementen der "Uberdeckung schneiden sollen. F"ur gr"oßere 
$\sss$ sind die R"aume holomorpher Funktionen auf den Schnittmengen von jeweils $\sss +1$
"Uberdeckungsmengen zu nehmen.

Die Halbnormensysteme der einzelnen Fr\'echetr"aume sind dann durch die Suprema auf kompakten
Aussch"opfungen der jeweiligen "Uberdeckungsmengen gegeben. Auf den Schnittmengen w"ahlt man
die Schnitte der jeweiligen kompakten Aussch"opfungsmengen.

Die Abbildungen $i_{\alpha_j}^\alpha$ f"ur ein $\alpha\in I^{\sss +1}$ und $ j\in \{0,\ldots,\sss\}$
kann man als Einschr"ankungen der holomorphen Funktionen jeweils auf durch den Schnitt mit
$U_{i_j}$ entstehende kleinen Menge $U_\alpha$ w"ahlen. Hierdurch sind die Eigenschaften
\auslabeln{10.1.1} sofort erf"ullt.

Die Koketten  werden   in diesem Beispiel zu den "ublichen Koketten
holomorpher Funktionen einer Cousin-I-Verteilung, $\delta^1 c=0$ ist dann exakt
die Zykelbedingung.

Wir fassen weitere Voraussetzungen, die wir benötigen, zu der folgenden Eigenschaft (E)
zusammen. Die Untereigenschaft (E.1) beschreibt die Bedingungen für die Kommutativität
des Diagramms \auslabeln{15.1.1}, die Exaktheit der Zeilen sowie Urbild- und 
Stetigkeitsabschätzungen. Die Untereigenschaft (E.2) beschreibt Urbildabschätzungen
der Spaltenabbildungen.

\begin{satzohnebeweis}{Definition}
Ein System $S=(S_{k,\sss})_{k,\sss}$ habe die Eigenschaft (E), falls die folgenden
Eigenschaften (E.1) und (E.2) erf"ullt sind:
\be
\item[(E.1):]  Ist $\sss\ge 0$ fest, so gibt es f"ur jedes $\alpha\in I^{\sss +1}$ eine 
exakte Sequenz von linearen stetigen Abbildungen $p^k_\alpha$, $k\ge 0$, die $\not=0$
sind, falls 
$F^{k+1}_\alpha\not=\{0\}$, und die 
invariant unter Permutationen der Indizes in $\alpha$ sind
\be
\item[(E.1.1)] \[    \ldots F_\alpha^3 \stackrel{p_\alpha^2}{\to} F_\alpha^2 \stackrel{p_\alpha^1}{\to}
                            F_\alpha^1 \stackrel{p_\alpha^0}{\to} F_\alpha^0 \to 0 \punkt \]
\ee
Es gelten ferner
\be
\item[(E.1.2)] \[     p_\alpha^k \circ \iota_{\alpha_j}^\alpha = \iota_{\alpha_j}^\alpha \circ p_{\alpha_j}^k
           \mbox{ f"ur } \alpha\in I^{\sss +1}\komma j=0,\ldots,\sss\komma k\ge 0 \komma \]
\ee
sowie die folgenden Eigenschaften:
\be
\item[(E.1.3)] Zu $\alpha\in I^{\sss +1}$, $k\ge 0$ und $n\in \N$ gibt es $C>0$, so daß es f"ur jedes 
$x\in F^k_\alpha$ mit $p_\alpha^{k-1} x=0$ und $\norm{\alpha,n+1}{x}\le 1 $ ein $y\in F^{k+1}_\alpha$
mit $p_\alpha^k y=x$ und $\norm{\alpha,n}{y}\le C$ gibt. Hierbei sei $p^{-1}_\alpha:F_\alpha^0 \to 0$
die Nullabbildung.
\item[(E.1.4)] F"ur alle $n\in \N$, $\alpha\in I^{\sss +1}$ und $k\ge 0$ gibt es ein Konstantensystem
$C_{n,k}=(C_{n,k,\alpha})_{\alpha \in I^{\sss +1}}$, so daß
\[   \dreinorm{k}{n,C_{n,k}}{p^k_\alpha x} \le \dreinorm{k+1}{n,1}{x} \mbox{ f"ur jedes } 
       x\in F_\alpha^{k+1}    \]
\ee
\item[(E.2):] F"ur $k\ge 1$ gelte:
Es gibt eine nicht leere Menge $M_\gamma$ zu einem festen $\gamma> 1$ mit
\be
\item[1.] $ M_\gamma \subset \{ [C:=(C_\alpha)_{\alpha\in I^{\sss +1}},C':=(C'_\alpha)_{\alpha\in I^{\sss}}]
           \doppelpunkt C\prec C' \}$,
\item[2.] aus $[C,C'] \in M_\gamma$ folgt $[C^\theta,{C'}^\theta]\in M_\gamma$ f"ur jedes $\theta\ge 1$,
\item[3.] ist $[C,C'] \in M_\gamma$ dann gibt es $C''\prec C$, so daß 
$[C'',C]\in M_\gamma$,
\ee
f"ur die gilt:
Zu jedem Konstantensystem $C_1=(C_{1,\alpha})_{\alpha\in I^{\sss +1}}$ und $n\in \N$ gibt es ein
Konstantensystem $C_2=(C_{2,\alpha})_{\alpha\in I^{\sss}}$, so daß f"ur alle $[C,C']\in M_\gamma$ gilt:

Ist $c\in {\cal C}_n( S_{k,\sss},C\cdot C_1)$ mit $\delta^\sss c=0$, dann existiert ein 
$c'\in {\cal C}_{n-1}( S_{k,\sss-1},{C'}^\gamma\cdot C_2)$ mit 
\formelohne{
      \delta^{\sss-1}c' &=& c \mbox{ und} \\
     \dreinorm{k}{n-1,{C'}^\gamma\cdot C_2}{c'} &\le& \dreinorm{k}{n,C\cdot C_1}{c}  } 
\ee
\end{satzohnebeweis}

Der nun folgende Satz zeigt, daß auch auf der Stufe $k=0$ L"osungen f"ur den verallgemeinerten
Korandoperator existieren, falls Eigenschaft (E) gilt.

Der Beweis verwendet das Induktionsverfahren, daß f"ur den Beweis des Theorem B von Cartan
auf relativ kompakten holomorph konvexen Stufen verwendet wird (vgl. Beweis von Theorem
7.4.3 bzw. 7.6.10 in \cite{Hoer}).

\begin{satz}{Satz}\einlabelnsatz{10.3}
 Sei $(S_{k,\sss})_{k,\sss}$ ein System von Fr\'echetr"aumen, f"ur das die 
Eigenschaft (E) gilt. Dann gilt f"ur $k=0$:

Sei $\gamma>1$ und $M_\gamma$ so gew"ahlt wie in (E.2), dann gibt es zu jedem Konstantensystem 
$C_1=(C_{1,\alpha})_{\alpha\in I^{\sss +1}}$ und $n\in \N$, $n>2(M+1)$ ein
Konstantensystem $C_2=(C_{2,\alpha})_{\alpha\in I^{\sss}}$, so daß f"ur alle $[C,C']\in M_\gamma$ gilt:

Ist $c\in {\cal C}_n( S_{0,\sss},C\cdot C_1)$ mit $\delta^\sss c=0$, dann existiert ein 
$c'\in {\cal C}_{n-2(M+1-\sss)}( S_{0,\sss-1},{C'}^{\gamma^{M-\sss}} \cdot C_2)$ mit 
\formelohne{
      \delta^{\sss-1}c' &=& c \mbox{ und} \\
     \dreinorm{0}{n-2(M+1-\sss),{C'}^{\gamma^{M-\sss}} \cdot C_2}{c'} 
                        &\le& \dreinorm{0}{n,C\cdot C_1}{c}  }
\end{satz}

Wir f"uhren eine Induktion "uber fallende $\sss$ durch.

Wir zeigen die folgende von $\sss$ abhängige Induktionsbehauptung für festes aber beliebiges
$k\ge 0$:

Zu jedem Konstantensystem 
$C_1=(C_{1,\alpha})_{\alpha\in I^{\sss +1}}$ und $n\in \N$, $n>2(M+1)$ gibt es ein
Konstantensystem $C_2=(C_{2,\alpha})_{\alpha\in I^{\sss}}$, so daß f"ur alle $[C,C']\in M_\gamma$ gilt:

Ist $c\in {\cal C}_n( S_{k,\sss},C\cdot C_1)$ mit $\delta^\sss c=0$ und
$p_\alpha^{k-1} c_\alpha=0$ für alle $\alpha\in I^{\sss+1}$, dann existiert ein 
$c'\in {\cal C}_{n-2(M+1-\sss)}( S_{k,\sss-1},{C'}^{\gamma^{M-\sss}} \cdot C_2)$ mit 
\formelohne{
      p_\alpha^{k-1} c'_\alpha &= & 0 \foralls \alpha\in I^\sss \\
      \delta^{\sss-1}c' &=& c \mbox{ und} \\
     \dreinorm{k}{n-2(M+1-\sss),{C'}^{\gamma^{M-\sss}} \cdot C_2}{c'} 
                        &\le& \dreinorm{k}{n,C\cdot C_1}{c}   \punkt}

Der Induktionsanfang $\sss=M+1$ ist klar, da nur Nullelemente auftauchen.

Zum Induktionschritt $\sss+1$ nach $\sss$:


 Sei $C_1=
(C_{1,\alpha})_{\alpha\in I^{\sss +1}}$ ein Konstantensystem und $n>2(M+1)$, $n\in \N$.
Sei ferner $[C=(C_{\alpha})_{\alpha\in I^{\sss +1}}, C'=(C'_{\alpha})_{\alpha\in I^{\sss }}]\in M_\gamma$.

Wir starten mit einem $c\in {\cal C}_{n}(S_{k,\sss},C\cdot C_1)$ mit $\delta^{\sss} c=0$
und $p_\alpha^{k-1}=0$ für alle $\alpha\in I^{\sss+1}$.
Nach (E.1.3) gibt es zu
$c_\alpha \in F_\alpha^k $ ein $d_\alpha \in F^{k+1}_\alpha$ mit $p_\alpha^k d_\alpha = c_\alpha$ und 
$\norm{\alpha,n-1}{d_\alpha} \le \tilde{C}_\alpha \norm{\alpha,n}{c_\alpha}$. Hierbei
hängt $\tilde{C}_\alpha$ nur von $\alpha$ und $n$ ab.

Wenn wir also $\tilde{C}:=(\tilde{C}_\alpha^{-1})_{\alpha\in I^{\sss +1}}$ setzen, so erhalten
wir ein $d=(d_\alpha)_{\alpha\in I^{\sss +1}} \in {\cal C}_{n-1}(S_{k+1,\sss},C\cdot\tilde{C}\cdot C_1)$
mit 
\[   \dreinorm{k+1}{n-1,C\cdot\tilde{C}\cdot C_1}{d} \le \dreinorm{k}{n,C\cdot C_1}{c}  \]

Wir wählen gemäß der 3. Eigenschaft von $M_\gamma$ ein $C''\prec C$, so daß
 $[C'',C]\in M_\gamma$. Nach Lemma \auslabelnsatz{10.2} folgt dann mit $D':=\tilde{C}\cdot C_1\cdot
(M(\sss+3))^{-1}$, daß $d':=\delta^{\sss+1} d \in{\cal C}_{n-1}(S_{k+1,\sss+1},C''\cdot D')$
und 
\[   \dreinorm{k+1}{n-1, C''\cdot D'}{d'} \le
                        \dreinorm{k+1}{n-1,C\cdot\tilde{C}\cdot C_1}{d} \]
Ist $\alpha \in I^{\sss+2}$, dann folgt wegen (E.1.2)
\formelohne{   p_\alpha^k d'_\alpha
                 &=& p_\alpha^k (\delta^{\sss+1} d)_\alpha  \\
                 &=& \sum_{j=0}^{\sss+2} (-1)^j p_\alpha^k i_{\alpha_j}^\alpha d_{\alpha_j} \\
                 &=& \sum_{j=0}^{\sss+2} (-1)^j i_{\alpha_j}^\alpha p_{\alpha_j}^k  d_{\alpha_j} \\           
                 &=& (\delta^{\sss+1} c)_\alpha \\
                 &=& 0 \punkt    }

Wir können also auf $d'$ die Induktionsvoraussetzung auf der Stufe $k+1$ anwenden. 
 Es gibt daher ein Konstantensystem $D''$, 
unabhängig von $C$, $C''$ und $d'$, so daß zu $d'$ eine 
Lösung der Gleichung $\delta^{\sss} d''= d'$ existiert mit $d''$ in 
${\cal C}_{n-1-2(M+1-\sss-1)}(S_{k+1,\sss},C^{\gamma^{M-\sss-1}}\cdot D'') $,
$p_\alpha^k d''_\alpha =0$ für alle $\alpha \in I^{\sss+1}$  und
\[   \dreinorm{k+1}{ n-1-2(M-\sss),C^{\gamma^{M-\sss-1}}\cdot D''}{d''} \le
      \dreinorm{k+1}{n-1,C''D'}{d'} \punkt   \]
Hierbei ist ohne Einschränkung $D''\le \tilde{C}\cdot C_1$.

Wir setzen nun $h:=d-d''\in {\cal C}( S_{k+1,\sss})$. Dann gilt 
\[   \delta^{\sss} h = \delta^{\sss} d - \delta^{\sss} d'' = d'-d' = 0  \mbox{ und} \]
\formelohne{
                 \dreinorm{k+1}{n-1-2(M-\sss),C^{\gamma^{M-\sss-1}} D''}{h}
     &\le&  \dreinorm{k+1}{n-1, C\cdot D''}{d} + \dreinorm{k+1}{n-1,C''D'}{d'}   \\
     &\le&  \dreinorm{k+1}{n-1, C\cdot\tilde{C}\cdot C_1}{d} + \dreinorm{k}{n,C\cdot C_1}{c} \\
     &\le&  2\dreinorm{k}{n,C\cdot C_1}{c}   }

Also ist $h\in  {\cal C}_{n-1-2(M-\sss)}(S_{k+1,\sss},C^{\gamma^{M-\sss-1}}\cdot D'')$

Wir wenden jetzt Eigenschaft (E.2) an: 

Es gibt demnach ein von $C$ und $C'$ unabhängiges $D'''=(D'''_\alpha)_{\alpha\in I^{\sss }}$,
 so daß für das Paar $[C^{\gamma^{M-\sss-1}},$ $
{C'}^{\gamma^{M-\sss-1}} ] \in M_\gamma$ gilt: Zu obigem $h\in {\cal C}_{n-1-2(M-\sss)}(S_{k+1 ,\sss},
C^{\gamma^{M-\sss-1}}\cdot D'')$ gibt es wegen $\delta^{\sss} h=0$ ein $h' \in
{\cal C}_{n-1-2(M-\sss)-1}$ $(S_{k+1,\sss-1},{C'}^{\gamma^{M-\sss}}\cdot D''')$, so daß 
\formelohne{
     \delta^{\sss-1} h'& = & h  \mbox{ und} \\
     \dreinorm{k+1}{n-1-2(M-\sss)-1,{C'}^{\gamma^{M-\sss}}\cdot D'''}{h'}
   &\le&  \dreinorm{k+1}{n-1-2(M-\sss),C^{\gamma^{M-\sss-1}}\cdot D''}{h}  }

Wir zeigen nun, daß $c':=(p^k_\alpha h'_\alpha)_{\alpha\in I^{\sss }}$ die gesuchte Lösung ist.

Es ist klar, daß $p_\alpha^{k-1}c_\alpha'=0$ für alle $\alpha\in I^{\sss}$ gilt. Ferner gilt für 
alle $\alpha\in I^{\sss+1}$, daß 
\formelohne{
             (\delta^{\sss-1} c')_\alpha 
     &=&  \sum_{j=0}^\sss (-1)^j \iota_{\alpha_j}^\alpha p^k_{\alpha_j} h'_{\alpha_j} \\
     &\stackrel{(E.1.2)}{=}&  \sum_{j=0}^\sss (-1)^j p^k_{\alpha} \iota_{\alpha_j}^\alpha  h'_{\alpha_j} \\
     &=&  p^k_\alpha(\delta^{\sss-1} h')_\alpha \\
     &=&  p^k_\alpha h_\alpha \\
     &=&  p^k_\alpha d_\alpha - p^k_\alpha d''_\alpha  \\
     &=&  c_\alpha \punkt   }

Wegen (E.1.4) gibt es zur Stufe $n-2(M+1-\sss)$ ein Konstantensystem $\tilde{\tilde{C}}=
(\tilde{\tilde{C}}_\alpha)_{\alpha\in I^{\sss +1}}$ mit 
\[  \dreinorm{k}{n-2(M+1-\sss),\tilde{\tilde{C}}_\alpha}{p^k_\alpha x} \le
    \dreinorm{k+1}{n-2(M+1-\sss),1}{ x}  \]
für jedes $x=(x_\alpha)_{\alpha\in I^{\sss +1}}$, $x_\alpha\in F^{k+1}_\alpha$.
Es folgt daher
\formelohne{
               \dreinorm{k}{n-2(M+1-\sss),{C'}^{\gamma^{M-\sss}}\cdot D'''\cdot \frac{
                  \tilde{\tilde{C}}}{2}}{c'}
    &\le&  \frac{1}{2} \dreinorm{k+1}{n-2(M+1-\sss),{C'}^{\gamma^{M-\sss}}\cdot D'''}{h'} \\
    &\le&  \frac{1}{2} \dreinorm{k+1}{n-1-2(M+1-\sss),C^{\gamma^{M-\sss-1}}\cdot D''}{h} \\
    &\le&  \dreinorm{k}{n,C\cdot C_1}{c}  }
Wir setzen also $C_2:= D'''\cdot \frac{\tilde{\tilde{C}}}{2}$ und erhalten damit
$c'\in {\cal C}_{n-2(M+1-\sss)}(S_{k,\sss-1},
{C'}^{\gamma^{M-\sss}}\cdot C_2)$. Damit ist der Induktionsschritt bewiesen \Ende

\newpage
\section[{\sc $\dquer$-Lösungen mit Abschätzungen}]{$\dquer$-Lösungen 
mit Abschätzungen auf Steinschen Mannigfaltig\-keiten }

Sei $\Omega$ eine komplexe Mannigfaltigkeit der Dimension $n$, die 
abzählbar gegen unendlich ist. $\oka$ sei die Garbe von Keimen analytischer Funktionen auf $\Omega$ und $\F$ eine kohärente analytische Garbe auf $\Omega$. Wir identifizieren für eine offene Menge $\Omega'\subset \Omega$ 
die holomorphen Funktionen auf $\Omega'$ mit den Schnitten $\Gamma(\Omega',\oka)$.

Eine $(p,q)$-Form $f$ auf $\Omega$ sei dadurch definiert, daß sie sich auf 
jeder Karte als $(p,q)$-Form in den lokalen Koordinaten schreiben läßt,
d.h. auf einer Karte $U$ existieren Koefizientenfunktionen aus 
$L^2(\Omega,loc)$ bezüglich eines festen Volumenelementes $dV$, so daß $f$
auf $U$ in den lokalen Koordinaten $z_1,\ldots,z_n$ die Form hat
\[   f= {\sum_{|I|=p}}'{\sum_{|J|=q}}' f_{I,J} dz^I\wedge d\overline{z}^J , \]
wobei hier die Notation aus [Hör] verwendet wird. Das bedeutet insbesondere,
daß die Summierung nur über aufsteigende Indizes der  Länge
p bzw. q erfolgt. 

Der Raum der $(p,q)$-Formen auf $\Omega$ heiße $L^2_{p,q}(\Omega,loc)$.

Wir wählen eine hermitsche Metrik $<\cdot,\cdot>$, die folgende 
Eigenschaft hat (siehe [Hör], S. 113): 

Es gibt eine Folge $(\eta_{\nu})_{\nu}$ in $C_{\circ}^{\infty}(\Omega)$,
so daß $0\le \eta_{\nu} \le 1$ und $\eta_{\nu}\equiv 1$ auf einem beliebigen
Kompaktum in $\Omega$, für $\nu$ genügend groß, und 
\[  \betrag{\dquer \eta_{\nu}}^2=<\dquer \eta_{\nu},\dquer \eta_{\nu}>\le 1
    \mbox{ auf } \Omega \leer{1} \nu=1,2,\ldots   \]

Mit einer Orthonormalisierung erreicht man, daß zu jedem $z\in \Omega$ eine
offene Umgebung $U$ existiert, sowie $(1,0)$-Formen $\omega_1,\ldots, \omega_n$
mit $C^{\infty}$-Koeffizienten, so daß 
\[   <\omega_j,\omega_k> = \delta_{jk} \mbox{ für } z\in U \]

Eine $(p,q)$-Form $f$ läßt sich in eindeutiger Weise in 
$\omega_1,\ldots, \omega_n$ entwickeln:
\[   f={\sum_{|I|=p}}'{\sum_{|J|=q}}' f_{I,J}\leer{0.2} \omega^I\wedge 
          \overline{\omega}^J      \mbox{ auf } U \]
und es gilt 
\[   <f,f>={\sum_{|I|=p}}'{\sum_{|J|=q}}' \betrag{f_{I,J}}^2 \mbox{ auf } U .\]
 
Sei $\phi\in C^2(\Omega)$ und $f\in L^2(\Omega,loc)$, so daß 
\formel{   \norm{\phi}{f}^2:= \int \betrag{f}^2 e^{-\phi} dV <\infty
\einlabeln{20}}
Der Raum dieser $(p,q)$-Formen sei mit $ L^2_{p,q}(\Omega,\phi)$ bezeichnet.

Sei von jetzt ab  $\Omega$ Steinsch.

Es gilt  das folgende Lemma über unbeschränkte
Operatoren zwischen Hilberträumen (siehe \cite{Hoer}, Lemma 4.1.1).

\begin{satzohnebeweis}{Lemma}\einlabelnsatz{3.1}

Sei $T:H_1\rightarrow H_2$ ein linearer dicht definierter abgeschlossener
Operator. Es gebe $C>0$, so daß
\formel{  
     \norm{H_2}{f}\le C \norm{H_1}{T^* f} \mbox{ für alle } f\in F\cap D_{T^*},
\einlabeln{3.0.5} }

wobei $F\subset H_2$ abgeschlossen ist und für das Bild von $T$ gilt 
$R_T\subset F$.

Dann existiert zu jedem $f\in F$ ein $g\in H_1$, so daß 
\formel{  Tg=f \mbox{ und } \norm{H_1}{g}\le C \norm{H_2}{f}\einlabeln{3.0.6}
 }
gilt. Hierbei kann in \auslabeln{3.0.6} die gleiche Konstante
$C$ gewählt werden wie in \auslabeln{3.0.5}.
\end{satzohnebeweis}

Tatsächlich beweist Hörmander die Äquivalenz von \auslabeln{3.0.5} und
\auslabeln{3.0.6}. Die Tatsache, daß in \auslabeln{3.0.6} die gleiche
Konstante gilt, wie in \auslabeln{3.0.5}, folgt unmittelbar aus der 
Anwendung des Satzes von Hahn-Banach  im Beweis.

$\dquer$ definiert folgende lineare abgeschlossene dicht definierte Operatoren
zwischen Hilberträumen:
\formel{   T:L_{p,q}^2(\Omega,\phi) & \rightarrow & L_{p,q+1}^2(\Omega,\phi)\\
           S:L_{p,q+1}^2(\Omega,\phi) & \rightarrow & L_{p,q+2}^2(\Omega,\phi)}

Den Raum der $(p,q)$-Formen auf $\Omega$ mit $C^{\infty}$-Koeffizienten 
und kompaktem Träger bezeichnen wir mit $D_{(p,q)}(\Omega)$.

Wir benötigen die folgende Aussage über die Lösungen von
$\dquer$ auf einer Steinschen Mannigfaltigkeit. Da die Beweisführung
sich weitestgehend auf den Beweis zu Theorem 5.3.4 in \cite{Hoer} stützt,
können wir uns hier kurz fassen.

\begin{satz}{Satz}\einlabelnsatz{3.3}
Sei $\Omega$ eine Steinsche Mannigfaltigkeit. Dann existiert eine 
streng plurisubharmonische Funktion $\psi\in C^2(\Omega)$, so daß
für jede plurisubharmonische Funktion $\phi$ auf $\Omega$ gilt:

Zu jedem $f\in L^2_{(p,q+1)}(\Omega,\phi)$ mit $\dquer f=0$ 
existiert eine Lösung $u\in L^2_{(p,q)}(\Omega,loc)$ der 
Gleichung $\dquer u=f$, so daß \formel{
     \int_{\Omega} \betrag{u}^2 e^{-\phi-\psi} dV \le
        \int_{\Omega} \betrag{f}^2 e^{-\phi} dV .\einlabeln{3.3.1} }
\end{satz}

Wir wählen eine streng plurisubharmonische Ausschöpfungsfunktion
$\kappa$ auf $\Omega$ (siehe z.B. \cite{Hoer}, Thm. 5.2.10). Dann
gibt es eine wachsende konvexe Funktion $\chi:\R\to \R$, so daß
für jede zweimal differenzierbare plurisubharmonische Funktion $\phi$
die sogenannte basic estimate 
\formel{  \norm{\chi(\kappa)+\phi}{f}^2 \le 
    \norm{\chi(\kappa)+\phi}{T^{\ast}f}^2 + \norm{\chi(\kappa)+\phi}{Sf}^2
    \foralls f\in D_{p,q+1}(\Omega)  
\einlabeln{3.3.2}} gilt.
Die basic estimate gilt auch für alle $f\in D_{T^{\ast}}\cap D_S$, da
$D_{p,q+1}(\Omega)$ bezüglich der Graphennorm 
\[   \norm{\chi(\kappa)+\phi}{f}^2 + 
    \norm{\chi(\kappa)+\phi}{T^{\ast}f}^2 + \norm{\chi(\kappa)+\phi}{Sf}^2 \]
dicht in $D_{T^{\ast}}\cap D_S$ liegt.

 Sei jetzt
$F$ der Kern von $S$, dann ist $F$ abgeschlossen in  $L^2_{(p,q+1)}(\Omega,
\chi(\kappa) + \phi)$,
da $S$ abgeschlossen ist. Wir erhalten also

\[  \norm{\chi(\kappa) + \phi}{f} \le \norm{\chi(\kappa) + \phi}{T^{\ast}f} 
  \leer{1}\mbox{ für alle }f\in F\cap D_{T^*} \punkt \]

Mit Lemma \auslabelnsatz{3.1} folgt hieraus, daß zu jedem $f\in 
L^2_{(p,q+1)}(\Omega,\chi(\kappa) + \phi)$ mit $\dquer f=0$ ein 
$u\in L^2_{(p,q)}(\Omega,\chi(\kappa) + \phi)$ existiert mit 
\[  \norm{\chi(\kappa)+\phi}{u} \le \norm{\chi(\kappa) + \phi}{f} 
    \le \norm{ \phi}{f} \punkt \]

Für den Übergang zu einer beliebigen  plurisubharmonischen
Funktion  $\phi$ verwenden wir 
die Standardargumentation wie z.B. im Beweis zu Theorem 4.4.2 in \cite{Hoer}. \hfill $\Box$

\begin{satz}{Korollar}\einlabelnsatz{3.4}
 Ist $\tilde{\Omega}\subset \Omega$ offen und
pseudokonvex, dann existiert zu jedem $f\in L^2_{(p,q+1)}(
\tilde{\Omega},\phi)$ mit $\dquer f=0$ eine Lösung $u\in L^2_{(p,q)}(
\tilde{\Omega},loc)$ der Gleichung $\dquer u=f$, so daß
\[ \int_{\tilde{\Omega}} \betrag{u}^2 e^{-\phi-\psi} dV \le
        \int_{\tilde{\Omega}} \betrag{f}^2 e^{-\phi} dV .  \]
Hier ist $\phi$ eine beliebige plurisubharmonische Funktion auf $\Omega$
 und $\psi$ aus obigem
Satz.
\end{satz}

Wir können auf $\tilde{\Omega}$ das eingeführte Volumenelement für
$\Omega$ eingeschränkt auf $\tilde{\Omega}$ verwenden.
Da \auslabeln{3.3.2} auch für alle $f\in D_{p,q+1}(\tilde{\Omega})$ gilt,
kann man wie im Beweis zu Satz \auslabelnsatz{3.3} argumentieren. \Ende

\begin{satz}{Lemma}\einlabelnsatz{3.4.1}
Seien $\Omega_1$, $\Omega_2$ und  $\Omega_3$ offene pseudokonvexe Teilmengen
einer komplexen Mannigfaltifkeit $\Omega$.
Sei $\Omega_1\subset \Omega_2$ relativ kompakt und $\Omega_2\subset \Omega_3$,
ferner sei $q\ge 1$ und $\phi \in \PSH(\Omega_2)$.

Dann gibt es zu jedem $f\in L^2_{(0,q)}(\Omega_2, \phi)$ mit
$\overline{\partial} f=0$ ein $F\in L^2_{(0,q)}(\Omega_3, loc) $ mit
\[ \overline{\partial} F=0 \mbox{ und } F=f \mbox{ auf } \Omega_1  \]
\end{satz}
Wir wenden Satz \auslabelnsatz{3.3} auf $f\in L^2_{(0,q)}(\Omega_2, \phi)$ an und
erhalten ein $u\in L^2_{(0,q-1)}(\Omega_2, loc) $ mit $\overline{\partial}u=f$.

Wir wählen eine $C^{\infty}$-Funktion $\chi$ mit Werten im Intervall $[0,1]$, 
so daß der Träger von $\chi $ in $\Omega_2$ liegt und $\Omega_1\subset \chi^{-1}(1)$.

Es sei $w:=\chi u$ und $F:=\overline{\partial} w$. Dann ist 
$F\in L^2_{(0,q)}(\Omega_3, loc)$ und $ \overline{\partial} F=0$.
Ferner gilt
\[   F\mid_{\Omega_1} = \left[ \overline{\partial} u\cdot \chi + u\cdot 
              \overline{\partial} \chi \right]\mid_{\Omega_1} = 
              \overline{\partial} u \mid_{\Omega_1} = f\mid_{\Omega_1} \punkt \]
\hfill $\Box$

\begin{satz}{Korollar}\einlabelnsatz{3.4.2}
Seien $\Omega_1\subset \Omega_2\subset\Omega_3$ pseudokonvexe Teilmengen von $\Omega$
aus Satz \auslabelnsatz{3.3} und sei $\Omega_1$ in $\Omega_2$ relativ kompakt.
 Sei ferner $\phi\in \PSH(\Omega)$ und $q\ge 1$.

Dann gibt es zu jedem $f\in L^2_{(0,q+1)}(\Omega_3, \phi)$ mit $\dquer f=0$
eine Lösung $u\in L^2_{(0,q)}(\Omega_3, loc) $ der Gleichung $\dquer u=f$, 
so daß die Abschätzung 
\[   \int_{\Omega_1} \betrag{u}^2 e^{-\phi-\psi} dV_{\Omega}
         \le \int_{\Omega_2} \betrag{f}^2 e^{-\phi} dV_{\Omega}  \]
mit $\psi$ aus Satz \auslabelnsatz{3.3} erfüllt wird.
\end{satz}
Wir wenden Korollar \auslabelnsatz{3.4} auf $f\in L^2_{(0,q+1)}(\Omega_3, \phi)$
mit $\dquer f=0$ zweifach an und erhalten $u_1\in L^2_{(0,q)}(\Omega_2, \phi + \psi)$
und $u_2\in L^2_{(0,q)}(\Omega_3, \phi + \psi)$ mit 
$f=\dquer u_1$ auf $\Omega_2$ und $f=\dquer u_2$ auf $\Omega_3$, sowie
\formelohne{  \int_{\Omega_2} \betrag{u_1}^2 e^{-\phi-\psi} dV_{\Omega}
             & \le& \int_{\Omega_2} \betrag{f}^2 e^{-\phi} dV_{\Omega}\und  \\
               \int_{\Omega_3} \betrag{u_2}^2 e^{-\phi-\psi} dV_{\Omega}
             & \le& \int_{\Omega_3} \betrag{f}^2 e^{-\phi} dV_{\Omega}\punkt  }

Dann erfüllt $w:=u_2|_{\Omega_2} -u_1$ die Voraussetzungen von Lemma \auslabelnsatz{3.4.1}.
Wir erhalten daher $W\in L^2_{(0,q)}(\Omega_3, loc)$ mit $\dquer W=0$ und 
$W=w$ auf $\Omega_1$.

Wir setzen $u:= u_2-W\in L^2_{(0,q)}(\Omega_3, loc)$. Da $\dquer u=\dquer u_2 - 
\dquer W = \dquer u_2=f$ gilt auf $\Omega_3$ und 
\formelohne{   \int_{\Omega_1} \betrag{u}^2 e^{-\phi-\psi} dV_{\Omega}
             & = & \int_{\Omega_1} \betrag{u_2-u_2+u_1}^2 e^{-\phi-\psi} dV_{\Omega}\\
             &\le& \int_{\Omega_2} \betrag{f}^2 e^{-\phi} dV_{\Omega}\komma }
ist das Korollar gezeigt \Ende

\newpage
\section[{\sc Beweis eines Theorem B mit Schranken für $\oka_{\Omega}$}]{Beweis eines Theorem B mit Schranken für Koketten mit Wer\-ten in 
der Garbe $\oka_{\Omega}$ von Keimen holomorpher Funktionen auf einer Steinschen
Mannigfaltigkeit}

Im nächsten Schritt werden wir ein Theorem B mit Lösungsabschätzungen
beweisen für die Garbe von Keimen holomorpher Funktionen auf
einer Steinschen Mannigfaltigkeit. Dies wird die Eigenschaft E.2  zeigen.

Ein ebensolches Lösungstheorem mit präzisen Angaben des Gewichtsverlustes
findet sich in \cite{Hoer}, Proposition 7.6.2. Dieses Ergebnis bezieht sich jedoch 
lediglich auf die Situation im $\C^n$. Deshalb streben wir hier
an, die Situation auf Steinschen Mannigfaltigkeiten zu untersuchen. Dabei
werden die Angaben des Gewichtsverlustes naturgemäß unpräziser. Der
Beweisgang wird im wesentlichen gleich verlaufen und zwar als Induktion
über $\sss$, wobei in den Beweis des Induktionsanfangs $\sss=1$ der
$\dquer$-Lösungssatz \auslabelnsatz{3.3} eingehen wird.

Sei also $\Omega$ eine Steinsche Mannigfaltigkeit der Dimension $n$ und
$\U:=(U_i)_{i\in I}$ eine abzählbare Überdeckung mit Steinschen Gebieten und
der Eigenschaft, daß eine feste Überdeckungsmenge
von höchstens $M$ verschiedenen Überdeckungsmengen geschnitten wird.
$\W:=(W_i)_{i\in I}$ sei eine Verfeinerung zu $\U$. Insbesondere gelte 
$W_i\subset U_i$ relativ kompakt und $W_i$ Steinsches Gebiet für jedes $i\in I$.
 $(\chi_i)_i$ sei eine Teilung der Eins
bezüglich $(W_i)_i$. Wir wählen ferner eine plurisubharmonische
Funktion $\phi$, die 
\[      \left( \sum_{i\in I}\betrag{\dquer \chi_i(z)}\right)^2 e^{-\phi}\le 1, 
           \leer{1} z\in \Omega \komma  \]
erfüllt. Es sei $\H_q$ die Garbe von Keimen von $(0,q)$-Formen.
Auf den Stufenindex des Korandoperators $\delta$ wird verzichtet,
wenn sich die entsprechende Stufe aus dem Zusammenhang ergibt.
Für einen Multiindex $i_0,\ldots ,i_\sss$ bezeichnen wir mit
$U_{i_0,\ldots ,i_\sss}$ den Schnitt $U_{i_0}\cap \ldots \cap U_{i_\sss}$.
Dies gelte für $W_{i_0,\ldots ,i_\sss}$ entsprechend.
Die Menge aller auf $\Omega$ plurisubharmonischen Funktionen bezeichnen
wie im folgenden mit PSH$(\Omega)$.

\begin{satz}{Poposition}\einlabelnsatz{3.5}
Sei $(\tilde{W}_i)_{i\in I}$ eine Familie von pseudokonvexen Gebieten mit 
$W_i\subset \tilde{W}_i$ relativ kompakt und $\tilde{W}_i\subset U_i$ für
jedes $i\in I$.
Ist $c\in C^\sss(\U,\H_q)$, $\delta c=0$, $\dquer c=0$ und 
\formel{ {\einlabeln{3.5.1} }  \sum_\indnullsigma \int_{\tilde{W}_\indnullsigma} \betrag{
             c_\indnullsigma}^2 e^{-\kappa} dV < \infty    }

für eine plurisubharmonische Funktion $\kappa$ auf $\Omega$, dann
existiert $c'\in C^{\sss-1}(\U,\H_q)$  mit $\dquer c'=0$, so
daß $\delta c'=c$ und 
\formel{    \sum_\indeinssigma \int_{W_\indeinssigma} \betrag{
             c'_\indeinssigma}^2 e^{-\kappa-(\phi+\psi)\cdot\sss} dV 
         \le C_\sss \sum_\indnullsigma \int_{\tilde{W}_\indnullsigma} \betrag{
             c_\indnullsigma}^2 e^{-\kappa} dV ,}
wobei $\psi\in \PSH(\Omega)$ unabhängig von $\kappa$ und $c$ gewählt
werden kann.
\end{satz}

{\sc Bemerkungen:} Das $\psi$ hier ist das gleiche wie in Satz \auslabelnsatz{3.3}.
Die Konstanten $C_\sss$ hängen nicht von den Gewichtsfunktionen ab.

Sei $c\in C^\sss(\U,\H_q)$, $\delta c=0$, $\dquer c=0$. Setze 
\[   b_\indeinssigma:=\sum_i\chi_i c_{i,\indeinssigma} \]
und $b:=(b_\indeinssigma)_\indeinssigma \in C^{\sss-1}(\U,\H_q)$.
Man beachte, daß die Summe jeweils endlich ist.
Es gilt:
\formel{   (\delta^\sminuseins b)_\indnullsigma &=& 
         \sum_{k=0}^\sss (-1)^k \sum_i \chi_i c_{i,\indnullsigmafehl{k}} \\
     &=& \sum_i \chi_i \sum_{k=0}^\sss (-1)^k c_{i,\indnullsigmafehl{k}} \\
     &=& \sum_i \chi_i c_\indnullsigma \\
     &=& c_\indnullsigma ,  }
denn $\delta^\sss c=0$.

Ist $U\subset U_\indeinssigma$ offen, dann gilt für $\beta\in PSH(\Omega)$
wegen der Cauchy-Schwarzschen Ungleichung:
\formel{{\einlabeln{3.5.3} } 
   \int_U \betrag{b_\indeinssigma}^2 e^{-\beta} dV &\le& 
   \int_U \sum_i \chi_i \betrag{c_{i,\indeinssigma}}^2
                                    e^{-\beta} dV  \\
   &\le&  \sum_i \int_{U\cap W_i} \betrag{c_{i,\indeinssigma}}^2 e^{-\beta} dV 
 \punkt}

Setze $f_\indeinssigma = (\dquer b)_\indeinssigma := \dquer 
b_\indeinssigma$. Hierdurch wird eine Kokette $f:=(f_\indeinssigma
)_\indeinssigma \in C^\sminuseins(\U,\-\H_{q+1})$ definiert. Es gilt:
\formelohne{  f_\indeinssigma =  \dquer b_\indeinssigma
         &=& \sum_i \dquer(\chi_i c_{i,\indeinssigma})  \\
         &=&   \sum_i \dquer\chi_i\wedge c_{i,\indeinssigma} , }
denn $\dquer c=0$. Da $\dquer$ auf jedes einzelne lokale Element
einer Kokette wirkt, vertauscht $\dquer$ mit dem Korandoperator $\delta$
auf jeder Stufe. Daher gilt
\[ \delta^\sminuseins f= \delta^\sminuseins \dquer b =  \dquer \delta^\sminuseins b
         = \dquer c=0  . \]
Für jedes $U\subset U_\indeinssigma$ und jedes $\beta\in PSH(\Omega)$ gilt
wegen der Wachstumseigenschaft von $\phi$
\formel{ {\einlabeln{3.5.4}} \int_U \betrag{f_\indeinssigma}^2 e^{-\phi-\beta} dV
      &\le&
          \int_U\left(\sum_i\betrag{\dquer \chi_i}\right)
             \left(\sum_i\betrag{\dquer \chi_i}\betrag{c_{i,\indeinssigma}}^2\right)
                                             e^{-\phi-\beta} dV \\
      &\le& \int_U \sum_i\betrag{\dquer \chi_i}\betrag{c_{i,\indeinssigma}}^2
                                             e^{-\frac{\phi}{2}-\beta} dV \\
      &\le& \sum_i \int_{U\cap W_i} \betrag{\dquer \chi_i}
                       \betrag{c_{i,\indeinssigma}}^2 e^{-\frac{\phi}{2}-\beta} dV \\ 
      &\le& \sum_i \int_{U\cap W_i} \betrag{c_{i,\indeinssigma}}^2 e^{-\beta} dV.  }

Wir führen nun eine Induktion über $\sigma$ durch:

Für den Induktionsanfang $\sss=1$ sei $c\in C^1(\U,\H_q)$, $\delta^1 c=0$,
$\dquer c=0$ und 
\[   \sum_{i,j} \int_{U_{i,j}}\betrag{c_{i,j}}^2 e^{-\kappa} dV<\infty \punkt \]
 $b$ und $f$ seien wie oben konstruiert, dann gilt 
$\delta^0 f=0$ und $\dquer f=\dquer\dquer b=0$. $f$ definiert also eine $(0,q+1)$-Form auf
$\Omega$, die $\dquer$-geschlossen ist.  Wir erhalten wegen \auslabeln{3.5.1} 
und \auslabeln{3.5.4} 
\formel{{\einlabeln{3.5.5} }  \int_{\Omega} \betrag{f}^2 e^{-\phi-\kappa} dV
    &\le&
           \sum_j \int_{W_j} \betrag{f_j}^2 e^{-\phi-\kappa} dV  \\
    &\le&
           \sum_j\sum_i \int_{W_i\cap W_j} \betrag{c_{i,j}}^2 e^{-\kappa} dV \\
    &<&   \infty }
Es gibt daher nach Satz \auslabelnsatz{3.3} ein $g\in L^2_{(0,q)}(\Omega,
loc)$ mit $\dquer g=f$ und 
\formel{ \int_{\Omega} \betrag{g}^2 e^{-\phi-\psi-\kappa}dV \le
         \int_{\Omega} \betrag{f}^2 e^{-\phi-\kappa}dV ,
\einlabeln{3.5.6}}
wobei $\psi$ unabhängig von $f$ und $\kappa+\phi$  gewählt werden kann.

Setze $c'_i := b_i-g$ für jedes $i\in I$. Dann wird durch $c':=(c'_i)_i$ eine Kokette in 
$C^0(\U,\H_q)$ definiert, für die 
\[   \delta^0 c'=\delta^0 b -\delta^0 g = \delta^0 b=c \] und
\[   \dquer c'=\dquer b-\dquer g= f-f=0 \mbox{   gilt. } \]
Wir können wegen \auslabeln{3.5.5} und \auslabeln{3.5.6} abschätzen
\formelohne{\sum_j \int_{W_j} \betrag{c'_j}^2 e^{-\phi-\psi-\kappa} dV
    &\le&
       2  \sum_j \int_{W_j} \betrag{b_j}^2 e^{-\phi-\psi-\kappa} dV +
       2  \sum_j \int_{W_j} \betrag{g}^2 e^{-\phi-\psi-\kappa} dV \\
    &\le&
       2  \sum_j \sum_i \int_{W_i\cap W_j} \betrag{c_{i,j}}^2 e^{-\phi-\psi-\kappa} dV +
       2M  \int_{\Omega} \betrag{g}^2 e^{-\phi-\psi-\kappa} dV \\
    &\le&
       2 \sum_{i,j} \int_{W_{i,j}} \betrag{c_{i,j}}^2 e^{-\kappa} dV +
       2M \sum_{i,j} \int_{W_{i,j}} \betrag{c_{i,j}}^2 e^{-\kappa} dV .  }
Mit $C_1:=2(M+1)$ ist der Induktionsanfang gezeigt.

Für den Induktionsschritt $\sminuseins\rightarrow \sss$ sei
$c\in C^\sss(\U,\H_q)$, $\delta^\sss c=0$, $\dquer c=0$ und
\[   \sum_{\indnullsigma} \int_{\tilde{W}_{\indnullsigma}}\betrag{c_{\indnullsigma}}^2
          e^{-\kappa} dV<\infty \punkt \]
 Wir bilden $b$ und $f$ wie oben. Dann ist $f\in C^\sminuseins(\U,\H_{q+1})$ und
$\delta^\sminuseins f=0$, $\dquer f=0$.

Wir können die Induktionsvorraussetzung auf $f$ und $\kappa + \phi$ anwenden,
denn wegen \auslabeln{3.5.4} gilt
\formel{{\einlabeln{3.5.7}}   \sum_\indeinssigma \int_{\tilde{W}_\indeinssigma} 
                                   \betrag{f_\indeinssigma}^2 e^{-\phi-\kappa} dV
    &\le&
           \sum_{i_0,\indeinssigma} \int_{\tilde{W}_{i_0,\indeinssigma}} 
                                   \betrag{c_\indnullsigma}^2 e^{-\kappa} dV \\
    &<&    \infty }
Es existiert also $b'\in C^{\sss-2}(\U,\H_{q+1})$, $\dquer b'=0$ mit 
$\delta^{\sss-2} b'=f$ und 
\formel{{\einlabeln{3.5.8}} \sum_\indzweisigma \int_{\tilde{\tilde{W}}_\indzweisigma} 
                            \betrag{b'_\indzweisigma}^2
                      & &\hspace{-1cm} e^{-(\kappa+\phi) - (\phi+\psi)(\sss -1)} dV \\
                       &\le&
         C_\sminuseins \sum_\indeinssigma \int_{\tilde{W}_\indeinssigma} \betrag{f_\indeinssigma}^2
                        e^{-\kappa-\phi } dV \komma}
wobei $(\tilde{\tilde{W}}_i)_i$ eine Verfeinerung zu $(\tilde{W}_i)_i$ ist mit 
$W_i\subset\tilde{\tilde{W}}_i$ und $\tilde{\tilde{W}}_i\subset \tilde{W}_i$ jeweils 
relativ kompakt für jedes $i\in I$.

Wir wenden Korollar \auslabelnsatz{3.4.2}  an auf jedes $b'_\indzweisigma$, $U_\indzweisigma$ und 
$\beta=(\kappa + \phi)+(\phi+\psi)(\sss-1)\in PSH(\Omega)$, sowie
$W_\indzweisigma\subset \tilde{\tilde{W}}_\indzweisigma \subset U_\indzweisigma$.

Es existiert also jeweils $b''_\indzweisigma\in L_{(0,q)}^2(U_\indzweisigma,loc)$
mit \formel{{\einlabeln{3.5.9}}   \dquer b''_\indzweisigma
     &=&
          b'_\indzweisigma    \mbox{ und } \\
          \int_{W_\indzweisigma} \betrag{b''_\indzweisigma}^2 e^{-\beta-\psi} dV
     &\le& 
          \int_{\tilde{\tilde{W}}_\indzweisigma} \betrag{b'_\indzweisigma}^2 e^{-\beta} dV   }
Sei $b'':= (b''_\indzweisigma)_\indzweisigma \in C^{\sigma-2}(\U,\H_q)$.
Mit \auslabeln{3.5.7}, \auslabeln{3.5.8} und \auslabeln{3.5.9} erhalten 
wir \formel{{\einlabeln{3.5.10}}  \sum_\indzweisigma \int_{W_\indzweisigma}
                      \betrag{b''_\indzweisigma}^2 e^{-\beta-\psi} dV 
   &\le&
         \sum_\indzweisigma \int_{\tilde{\tilde{W}}_\indzweisigma}
                      \betrag{b'_\indzweisigma}^2 e^{-(\kappa+\phi)-(\phi+\psi)(\sigma-1)} dV \\ 
   &\le&
         C_{\sigma-1} \sum_\indeinssigma \int_{\tilde{W}_\indeinssigma} 
           \betrag{f_\indeinssigma}^2 e^{-\phi-\kappa} dV   \\
   &\le&
         C_{\sigma-1} \sum_\indnullsigma \int_{\tilde{W}_\indnullsigma} 
           \betrag{c_\indnullsigma}^2 e^{-\kappa} dV.   }

Wir setzen $b''':= \delta^{\sigma-2} b''$. Dann ist $b'''\in C^{\sigma-1}(\U,\H_q)$.
Schließlich setzen wir $c':= b-b''' \in C^{\sigma-1}(\U,\H_q)$. Dann gilt
\formelohne{   \delta^{\sigma-1} c' = \delta^{\sigma-1}b - \delta^{\sigma-1}b''' 
         &=&  \delta^{\sigma-1}b - \delta^{\sigma-1}\delta^{\sss-2} b'' \\
         &=&  \delta^{\sigma-1}b \\
         &=&  c .  }
und
\formelohne{   \dquer c' = \dquer b - \dquer b''' 
         &=&
              f- \dquer\delta^{\sss-2} b'' \\
         &=&  f-\delta^{\sss-2} \dquer b'' \\
         &=&  f-\delta^{\sss-2} b'  \\
         &=&  f-f\\
         &=&  0 }

Es gilt $\beta + \psi= \kappa + (\phi+\psi)\sigma$. Daher erhalten wir
aus \auslabeln{3.5.10} und \auslabeln{3.5.4}:
\formelohne{  \sum_\indeinssigma  \int_{W_\indeinssigma}& &\leer{-1.5}\betrag{c'_\indeinssigma}^2
                   e^{-\kappa-(\phi+\psi)\sigma} dV \\ 
   &\le&
           2\sum_\indeinssigma \sum_i \int_{W_i\cap W_\indeinssigma}\betrag{c_{i,\indeinssigma}}^2
                   e^{-\kappa} dV  \\
   &+ & \leer{0}  2\sum_\indeinssigma \int_{W_\indeinssigma}\betrag{\sum_{j=1}^\sss 
                 (-1)^{j-1} b''_\indeinssigmafehl{j}}^2  e^{-\beta-\psi} dV \\
   &\le&  2   \sum_\indnullsigma \int_{W_\indnullsigma}\betrag{c_\indnullsigma}^2
                   e^{-\kappa} dV  \\
   &+ & \leer{0}   2\sum_\indeinssigma \sum_{j=1}^\sss \int_{W_\indeinssigma}\betrag{
                  b''_\indeinssigmafehl{j}}^2  e^{-\beta-\psi} dV   }
Es gilt für den letzten Summanden
\formelohne{  \sum_\indeinssigma \sum_{j=1}^\sss \int_{W_\indeinssigma}& &\leer{-1.5}\betrag{
                  b''_\indeinssigmafehl{j}}^2  e^{-\beta-\psi} dV \\
    &=&  \sum_\indzweisigma \sum_i \int_{W_i\cap W_\indzweisigma} \betrag{
                  b''_\indzweisigma}^2  e^{-\beta-\psi} dV \\
    &+&   \sum_\indeinssigma \sum_{j=2}^\sss \int_{W_\indeinssigma}\betrag{
                  b''_\indeinssigmafehl{j}}^2  e^{-\beta-\psi} dV \\
    &\le&  \sum_\indzweisigma M \int_{W_\indzweisigma} \betrag{b''_\indzweisigma}^2
                  e^{-\beta-\psi} dV \\
    &+&   \sum_\indeinssigma \sum_{j=2}^\sss \int_{W_\indeinssigma} \betrag{
                  b''_\indeinssigmafehl{j}}^2  e^{-\beta-\psi} dV, }
denn $W_\indzweisigma$ wird von höchstens $M$ Mengen aus $\W$ geschnitten.
Wenn wir diese Abschätzung für $j=2,\ldots,\sss$ iterativ fortsetzen erhalten wir
\formelohne{   \sum_\indeinssigma \sum_{j=1}^\sss \int_{W_\indeinssigma}& &\leer{-1.5}\betrag{
                  b''_\indeinssigmafehl{j}}^2  e^{-\beta-\psi} dV \\
   &\le&
         \sss M  \sum_\indzweisigma  \int_{W_\indzweisigma}\betrag{
                  b''_\indzweisigma}^2  e^{-\beta-\psi} dV \\
   &\le& \sss M C_{\sss-1} \sum_\indnullsigma  \int_{\tilde{W}_\indnullsigma}\betrag{c_\indnullsigma}^2
                   e^{-\kappa} dV. }
Wir haben also mit $C_\sss := 2(1+\sss^2M C_{\sss-1})$ den Induktionsschritt gezeigt $\Box$

{\sc Bemerkung:} Der Beweis zu Proposition \auslabelnsatz{3.5} ist auch mit 
$W_i=U_i$ für jedes $i\in I$ durchzuführen. Anstelle von Korollar \auslabelnsatz{3.4.2} 
ist dann Korollar \auslabelnsatz{3.4} anzuwenden. Daher gilt Proposition \auslabelnsatz{3.5}
sinngemäß auch für $W_i=\tilde{W}_i=U_i$ für jedes $i\in I$.

Wir schöpfen jedes $U_i$, $i\in I$ durch pseudokonvexe Gebiete $A_{i,n}$, $n\in \N$ aus, für die
für jedes $n\in \N$ jeweils $A_{i,n}\subset A_{i,n+1}$ relativ kompakt gilt.
Ohne Einschränkung sei $\bigcup_{i\in I} A_{i,1}=\Omega$. Für 
$(i_0,\ldots, i_\sss)\in I^{\sss +1}$ sei $A_{i_0,\ldots, i_\sss,n}:=A_{i_0,n}\cap\ldots\cap
A_{i_\sss,n}$.

Für eine Kokette der Länge $\sss$, die aus Schnitten über den Überdeckungsmengen aus der
Über\-deckung $\U$ besteht
mit Werten in der Garbe von Keimen von $L^2$-Formen sei für $\phi\in PSH(\Omega)$
\[   \norm{\phi,n}{c}^2 := \sum_\indnullsigma \int_{A_{i_0,\ldots, i_\sss,n}}
            \betrag{c_\indnullsigma}^2   e^{-\phi} dV  \]
bzw.
\[   \norm{\phi}{c}^2 := \sum_\indnullsigma \int_{U_\indnullsigma}
            \betrag{c_\indnullsigma}^2   e^{-\phi} dV  \]

eingeführt. Mit $C_n^\sss(\U,\H_q,\phi)$ bzw. $C^\sss(\U,\H_q,\phi)$ 
 bezeichnen wir die Koketten aus $C^\sss(\U,\H_q)$, 
für die 
\[  \norm{\phi,n}{c}^2 < \infty \mbox{ bzw. }       \norm{\phi}{c}^2 < \infty \mbox{ gilt.}   \]
Proposition \auslabelnsatz{3.5} können wir also auch wie folgt formulieren:

\begin{satzohnebeweis}{Satz}\einlabelnsatz{3.6}
Es gibt eine plurisubharmonische Funktion $\psi$, so daß für jedes $n\in \N$,
 eine beliebige plurisubharmonische
Funktion $\kappa$ und $c\in C_{n+1}^\sss(\U,\H_q,\kappa)$ mit $\delta c=0$ und $\dquer c=0$
ein $c'\in C_{n}^{\sss-1}(\U,\H_q)$ mit $\delta c'=c$ und $\dquer c'=0$ existiert, so daß
\[   \norm{\kappa+\psi,n}{c'}^2 \le \norm{\kappa,n+1}{c}^2  \]
erfüllt ist.
\end{satzohnebeweis}

Ein wichtiger Spezialfall von Satz \auslabelnsatz{3.6} ist $q=0$,
 denn die $\dquer$-geschlossenen Schnitte über $\Omega$ mit Werten in
 $\H_0$ sind die  holomorphen Funktionen auf $\Omega$.

\newpage
\section[{\sc Konstruktion einer exakten Sequenz}]{Konstruktion einer exakten Sequenz} \label{Kap5}
Es sei $\Omega$ eine Steinsche Mannigfaltigkeit und $U\subset\Omega$ offen
holomorph konvex. Sei ${\cal{F}}$ eine kohärente analytische Garbe auf $\Omega$.
Ist $U\subset\Omega$ relativ kompakt, dann kann auf den Schnitten $\Gamma(U,{\cal{F}})$
in der folgenden Weise eine Fr\'echetraumtopologie eingeführt werden:

Nach Theorem A von Cartan gibt es $F_1,\ldots,F_q\in \Gamma(\Omega,{\cal{F}})$, die
${\cal{F}}$ in $U$ erzeugen. In der Umgebung eines jeden Kompaktums $K_p$, $p\in \N$
einer holomorph konvexen kompakten Ausschöpfung $(K_p)_{p\in \N}$ von $U$ wird
${\cal{F}}$ von $F_1,\ldots,F_q$ erzeugt, d.h. für jeden Schnitt $x$ von einer 
Umgebung von $K_p$ in ${\cal{F}}$ gibt es auf einer Umgebung von $K_p$ holomorphe 
Funktionen $c_1,\ldots,c_q$, so daß in einer Umgebung von $K_p$ die Darstellung
\[    x=\sum_{j=1}^q c_j F_j   \]
existiert.

Es ist für jedes $p\in N$ 
\formel{    \norm{p}{x}:= \inf\limits_{x=\sum_j c_j F_j}\sum_{j=1}^q\ffb{K_p}{c_j} 
        \einlabeln{7.0.0} }
eine Halbnorm und $(\Gamma(U,{\cal{F}}),(\norm{p}{\leer{1}})_p)$ ein Fr\'echetraum
(siehe z.B. \cite{Hoer},Cor. 7.2.6 ).

Der Übergang zu einem anderen Erzeugendensystem liefert ein äquivalentes Halbnormensystem.
Wir bezeichnen mit $H_U({\cal{F}})$ den Raum $\Gamma(U,{\cal{F}})$ zusammen mit 
erzeugenden Schnitten $F_1,\ldots,F_q\in \Gamma(\Omega,{\cal{F}})$ und einer 
holomorph konvexen kompakten Ausschöpfung von $U$. Auf dem Fr\'echetraum $H_U({\cal{F}})$
ist damit eine Fr\'echetraumgradierung festgelegt. 

$H_U(\oka)$ fällt mit dem Fr\'echetraum der holomorphen Funktionen in $U$ (mit
dem durch Suprema auf der kompakten Ausschöpfung gebildeten Normensystem) zusammen,
falls als Erzeugendensystem die konstante Funktion $1$ gewählt wird. Im folgenden
werden wir $H_U(\oka)$ stets in dieser Weise auffassen, falls nichts anderes vermerkt ist.

Sei nun ${\cal{F}}\subset \oka^r$ eine lokal endlich erzeugte Untergarbe, dann ist 
${\cal{F}}$ kohärent nach dem Satz von Oka (siehe z.B. \cite{Hoer},Theorem 7.1.5). Sei ferner 
$\ffc{F}{q}$ ein Erzeugendensystem von ${\cal{F}}$ in $U$ und $(K_p)_p$ eine
holomorph konvexe kompakte Ausschöpfung.
 Wir wollen nun die entsprechende Gradierung $(\norm{p}{\leer{1}}
)_p$ auf $H_U(\F)$ mit den Suprema auf den $K_p$ vergleichen.
Für $f\in \Gamma(U,\oka^r)$ sei $[f]_j$ die $j$-te Komponente von $f$. 
$\Gamma(U,{\cal{F}})$ ist als Teilraum
des Raumes der $r$-tupel holomorpher Funktionen in $U$ in der kompakt 
gleichmäßigen Topologie abgeschlossen (siehe \cite{Hoer}, Theorem 7.2.12), also induzieren
die Supremumsnormen 
\[   \betrag{x}_p := \ffa{j}{r} \ffb{K_p}{[x]_j} \komma x\in \Gamma(U,\oka^r) \]

eine Fr\'echetraumgradierung auf $\Gamma(U,{\cal{F}})$. Es gilt das folgende

\begin{satz}{Lemma}\einlabelnsatz{7.0.1}
Für jedes $p\in \N$ gilt \be
\item $\betrag{x}_p \le C_p \norm{p}{x}$ und 
\item $\norm{p}{x}\le C'_p \betrag{x}_{p+1}$
\ee
für jedes $x\in \Gamma(U,{\cal{F}})$ mit unabhängigen Konstanten $C_p$ und $C'_p$.
\end{satz}

Sei $g\in \Gamma(U,\F)$, $g=\ffa{i}{q} d_j F_j$, dann ist 
\formelohne{
    \betrag{g}_p &=& \ffa{j}{r} \ffb{K_p}{[g]_j}    \\
        &=& \ffa{j}{r} \ffb{K_p}{[\ffa{i}{q} d_iF_i]_j}  \\
        &\le& \ffa{j}{r}\ffa{i}{q} \ffb{K_p}{d_j} \ffb{K_p}{[F_i]_j}  \\
        &\le&  C_p \ffa{j}{r} \ffb{K_p}{d_j}  }

Durch Übergang zum Infimum über alle Darstellungen $g=\ffa{i}{q} d_i F_i$
erhält man $\betrag{g}_p\le C_p\norm{p}{g}$.

Sei $p\in \N$ fest, dann wähle eine holomorph konvexe offene Menge $U'$
und $K_{p+1}\supset U' \supset K_p$. Sei $(A_n)_{n\in \N}$ eine 
Ausschöpfung von $U'$ mit holomorph konvexen kompakten Mengen mit 
$A_1=K_p$. Die Abbildung

\formel{ {\einlabeln{7.1.1}  }     \phi: H_{U'}(\oka^q) &\rightarrow &\Gamma(U',\F) \\
            \ffc{c}{q} &\mapsto & \ffa{j}{q} c_jF_j         
}
ist stetig linear zwischen Fr\'echeträumen, surjektiv und offen (Theorem 
7.2.12 \cite{Hoer}).

Da $\phi$ offen ist, gibt es zu jedem $n\in \N$ ein $m\in \N$ und $C_{n,m}$,
so daß zu jedem $g\in \Gamma(U',\F)$ ein Element $c=(\ffc{c}{q})\in
H_{U'}(\oka^q)$ existiert mit $g=\ffa{j}{q} c_jF_j$ und

\[   \ffa{j}{q} \ffb{A_n}{c_j} \le C_{n,m} \ffa{i}{r}\ffb{A_m}{[g]_i}\punkt \]

Hieraus folgt 
\[   \inf_{g=\sum\limits_{j=1}^{q}c_jF_j} \ffa{j}{q} \ffb{K_p}{c_j} \le
               C\ffa{i}{r} \ffb{K_{p+1}}{[g]_j} ,  \]
für $g\in \Gamma(U,\F)$.
Also gilt $\norm{p}{ g } \le C'_p \betrag{g}_{p+1}$  $\Box$

Wir wollen für Abbildungen wie in \auslabeln{7.1.1}  Urbildabschätzungen
zeigen, falls der Bildraum die Gradierung \auslabeln{7.0.0} trägt:

\begin{satz}{Lemma}\einlabeln{7.2a}
Sei $U$ eine Steinsche Mannigfaltigkeit und sei $(K_p)_p$ eine Ausschöpfung von 
$U$ mit holomorph konvexen kompakten Mengen. Seien ferner $\ffc{F}{q}\in 
\Gamma(U,\oka^r)$ und ${\cal{F}}\subset \oka^r$ die von $\ffc{F}{q}$ erzeugte
Untergarbe. Dann ist 
\[   \phi: H_U(\oka^q) \rightarrow H_U({\cal{F}})\komma \ffc{c}{q}\mapsto 
     \ffa{j}{q} c_jF_j   \]
eine stetige, lineare und surjektive Abbildung und es gilt für jedes $\theta >1$:
Zu jedem $f\in H_U({\cal{F}})$ gibt es ein $g\in H_U(\oka^q)$ mit $\phi(g)=f$
und $\betrag{g}_p \le \theta\norm{p}{f} $.
\end{satz}

Nach Theorem 7.2.12 aus \cite{Hoer} ist $\phi$, aufgefaßt als Abbildung in $\Gamma(U,{\cal{F}})$,
linear, stetig und surjektiv, wobei $\Gamma(U,{\cal{F}})$ die kompakt gleichmäßige Topologie 
von $H_U(\oka^q)$ trägt. Wegen Lemma \auslabelnsatz{7.0.1} ist daher $\phi$,
aufgefaßt als Abbildung in $H_U({\cal{F}})$ ebenfalls stetig.

Sei $\theta>1$. Sei nun $f\in H_U({\cal{F}})$ mit $\norm{p}{f}\le 1$. Dann gibt es  
nach der Definition von $\norm{p}{f}$
in einer Umgebung von $K_p$ holomorphe Funktionen $c=\ffc{c}{q}$, so 
daß dort $f=\ffa{j}{q}c_j F_j$ gilt und $\ffa{j}{q} \ffb{K_p}{c_j} < \sqrt{\theta}\norm{
p}{f}\le \sqrt{\theta}$.

Sei $\zzd:=\zzd(\ffc{F}{q})$ die Relationsgarbe von $\ffc{F}{q}$, 
d.h. $\Gamma(U,\zzd)$ ist der Kern von $\phi$. $\zzd$ ist 
nach dem Satz von Oka kohärent in $U$. Ist $d=\ffc{d}{q}\in H_U(\oka^q)$
und $\ffa{j}{q} d_j F_j=f$, dann ist $c-d$ ein Schnitt in die Relationsgarbe über
einer Umgebung von $K_p$.

Wir können eine Folge $(h_n)_n$ in $\Gamma(U,\zzd)$ finden, so daß $\norm{p}{
h_n-(c-d)}\rightarrow 0$ für $n\rightarrow\infty$ (siehe \cite{Hoer}, Theorem 7.2.7). Nach
Lemma \auslabelnsatz{7.0.1} folgt 
\[   \ffa{j}{q} \ffb{K_p}{[h_n]_j-(c_j-d_j)}\rightarrow 0 \mbox{ für }
n\rightarrow \infty           \]
Zu jedem $\eps>0$ können wir also $h=\ffc{h}{q}\in \Gamma(U,\zzd)$ finden,
so daß $\ffa{j}{q}\ffb{K_p}{(h_j+d_j) - c_j}<\eps$. Also gilt
$h+d\in \Gamma(U,\oka^q)$, $\phi(h+d)=\ffa{j}{q} (h_j+d_j)F_j=\ffa{j}{q} d_jF_j=f$ und
$\ffa{j}{q} \ffb{K_p}{h_j+d_j}<\sqrt{\theta }+ \eps$. 
Mit $\eps:=\theta-\sqrt{\theta}$ folgt daher $\betrag{h+d}_p\le \theta$. Hieraus
folgt die Behauptung $\Box$

Sei $\F$ wieder allgemein eine kohärente analytische Garbe über $\Omega$ 
und $\ffc{F}{q}\in \Gamma(\Omega,\F)$  erzeugende Schnitte von $\F$ über $U\subset\Omega$ 
offen, holomorph konvex und relativ kompakt. Es gilt das folgende 

\begin{satz}{Lemma}\einlabelnsatz{7.2}
Es gibt eine exakte Sequenz 
\[  \ldots  H_U(\oka^{r_3})\stackrel{p_2}{\rightarrow}
            H_U(\oka^{r_2})\stackrel{p_1}{\rightarrow}
            H_U(\oka^{r_1})\stackrel{p_0}{\rightarrow} \zza \rightarrow 0  \]
Ferner gilt für $i\ge 1$: Es gibt ein $C>0$, so daß 
$B_{p+1}^i\cap \abbkern p_{i-1}\subset p_i(C B^{i+1}_p)$, wobei
$B^i_p:=\{f\in H_U(\oka^{r_i}):\norm{p}{f}\le 1\}$, $i\ge 1$ und 
$B^0_p:=\{ f\in \zza : \norm{p}{f}\le 1\}$.

Für $i=0$ gilt sogar $B^0_p\subset p_0(\theta B_p^1)$ für jedes $\theta>1$.
\end{satz}
1) $p_0: H_U(\oka^q)\rightarrow \zza$, $\ffc{c}{q}\mapsto \ffa{j}{q} c_j F_j $ ist eine
lineare, stetige und surjektive Abbildung nach Lemma \auslabelnsatz{7.2a}. Setze also $r_1=q$.
Es folgt sofort $B^0_p\subset p_0(\theta B_p^1)$ für jedes $\theta>1$.

2) Sei $\Rscript$ die Relationsgarbe von $\ffc{F}{q}$. 
Es gibt  Schnitte $\ffc{G}{r_2}\in \Gamma(\Omega,\zzd)$, die $\Rscript$ über $U$ erzeugen.
Wie unter 1) ist $\tilde{p}_1:H_U(\oka^{r_2})$ $\rightarrow H_U(\zzd)$, 
$\ffc{c}{r_2}\mapsto \ffa{j}{r_2} c_j G_j$ linear, stetig und surjektiv und es gilt
$\{f\in H_U(\zzd):\norm{p}{f}\le 1\} \subset \tilde{p}_1(\theta B_p^2)$ für $\theta>1$.

Wenn wir mit $p_1$ die Abbildung $\tilde{p}_1$, aufgefaßt als Abbildung in $H_U(\oka^{r_1})$,
bezeichnen, dann gibt es nach Lemma \auslabelnsatz{7.0.1}, jetzt angewandt auf
$\Rscript$, ein $C'>0$, so daß 
\[ \{ f\in H_U(\oka^{r_1}): p_0 f=0 \mbox{ und } \betrag{f}_{p+1} \le \frac{1}{C'} \}
    \subset \{ f\in H_U(\zzd): \norm{p}{f}\le 1\} \subset \tilde{p}_1(\theta B_p^2) = 
      p_1(\theta B_p^2)      \]

Mit $C:=C'\cdot \theta$ folgt die Behauptung für $i=1$, denn $B_{p+1}^1=
C'\cdot\{ f\in H_U(\oka^{r_1})\doppelpunkt \betrag{f}_{p+1} \le \frac{1}{C'} \}$.

Durch Iteration des zweiten Schrittes mit den jeweiligen Relationsgarben folgt die
Behauptung für alle $i$ $\Box$

Wir betrachten für $U\subset \Omega$ offen holomorph konvex den Raum $l_1(\zza)
:= \{ x=(x_j)_j \in H_U({\cal{F}})^{\N}, \sum_{j=1}^{\infty} \norm{p}{x_j}<\infty
\forall p\in \N \}$.
Dies ist wiederum ein Fr\'echetraum mit dem Halbnormensystem 
\[   \norm{p}{(x_j)_j} := \ffa{j}{\infty} \norm{p}{x_j},\leer{1} p\in \N,  \]
wobei $\norm{p}{\cdot}$ die $p$-te Halbnorm auf $\zza$ ist bezüglich einer 
holomorph konvexen kompakten Ausschöpfung $(K_p)_p$ von $U$.

Wir wollen eine exakte Sequenz 
\[ \cdots\to \zzc\stackrel{p_2^U}{\rightarrow}
          \zzc\stackrel{p_1^U}{\rightarrow}
          \zzc\stackrel{p_0^U}{\rightarrow} \zza\rightarrow 0       \]
mit stetig linearen Abbildungen konstruieren. Hierbei soll für jedes 
$U'\subset U\subset\subset \Omega$ jeweils
\[   p_i^{U'}\mid_{{\displaystyle \zzc}} = p_i^U \mbox{ für } i\ge 0  \]
gelten. Hierzu benötigen wir zunächst die folgende



\begin{satz}{Proposition}\einlabelnsatz{7.3}
 Sei $\F$ ein kohärente analytische Garbe auf einer Steinschen
Mannigfaltigkeit $\Omega$, $(K_t)_{t\ge 0}$ eine holomorph konvexe kompakte
Ausschöpfung von $\Omega$. Sei $(f_j)$ eine Folge in $\bigoplus_{l\in \N} H_{\Omega}(\F)$.

Sei $p_U:\zzc \rightarrow \ffd{U}{\F}$, $(\lambda_j)_j \mapsto \sum_j \lambda_j f_j$
für jedes $U\subset\subset \Omega$ offen, holomorph konvex, eine wohldefinierte stetige lineare 
Abbildung. 

Es gebe Folgen natürlicher Zahlen $(a_t)_{t\ge 0}$, $(b_t)_{t\ge 0}$ mit
\[   a_t\le b_t \le a_{t+1},\, t\ge 0 \mbox{ und } b_t<b_{t+1},\, t\ge 0, \]
so daß auf einer Umgebung einer Stufe $K_s$ es Gleichungen
\[    \ffe{f}{b}{t}{i} = \ffa{j}{b_t} c_{t,i,j} f_j  \]
für alle $t\ge s$ und $b_t+i\le a_{t+1}$ gibt mit Funktionen $c_{t,i,j}$,
die auf einer Umgebung von $K_s$ holomorph sind und
\formel{{\einlabeln{7.0.2}}  \sup_{{ t\ge s \atop 1\le i\le a_{t+1}-b_t } }
           \ffa{j}{b_t} \ffb{K_s}{c_{t,i,j}} <1   }
erfüllen.

Ferner gelte für eine Folge $\lambda=(\lambda_j)_j\in \abbkern p_U$, $U\relkom K_s$,
$r$ beliebig, daß aus $\lambda_{b_r+i}=0$ für alle $r\ge t\ge s$, $i=1,\ldots,
a_{r+1}-b_r$ schon 
\formel{    \lambda_j=0 \leer{2}\mbox{ für alle $j\ge a_t + 1$ folgt.}  
\einlabeln{7.0.3}}

Sind die vorstehenden Voraussetzungen erfüllt, dann gibt es eine Folge 
$(g_j)_j\in \ffd{\Omega}{\oka}$, so daß für ein $U\relkom \Omega$ offen, holomorph konvex
die Abbildung 
\formelohne{  q_U: l_1(H_U(\oka)) &\rightarrow& l_1(H_U(\oka))  \\
              \lambda=(\lambda_j)_j & \mapsto & \sum_j \lambda_j g_j  }
stetig linear ist und $\bild q_U = \abbkern p_U$ ist.
Ferner gelten für $(g_j)_j$ die Eigenschaften \auslabeln{7.0.2} und \auslabeln{7.0.3}
entsprechend.

\end{satz}

Für jedes offene holomorph konvexe $U\relkom\Omega$ sei eine holomorph konvexe kompakte
Ausschöpfung $(A_p)_p$ fest gewählt. Wir bezeichnen mit $\normleer_p$, 
$p\in \N$ die $p$-te Halbnorm des Fr\'echetraumes $l_1(H_U(\oka))$.

Sei $U\subset K_s$ offen holomorph konvex. Definiere für $t\ge s+1$ eine Abbildung $P_{t}^{t-1}:
l_1(H_U(\oka)) \rightarrow l_1(H_U(\oka))$, $\lambda:=(\lambda_l)_l \mapsto
\mu:=(\mu_l)_l$ durch 
\[   \mu_l:= \left\{  \begin{array}{l@{\leer{0.5} ; \leer{0.5}}l} 
             \lambda_l + \ffa{i}{a_t-b_{t-1}} \lambda_{b_{t-1}+i} c_{t-1,i,l} & l\le b_{t-1} \\
             0 & b_{t-1}+1\le l \le a_t \\
             \lambda_l & \mbox{sonst}  
          \end{array}  \right.       \punkt\]
Eigenschaften von $P_{t}^{t-1}$: \begin{enumerate}
\item $P_{t}^{t-1}(\abbkern p_U) \subset \abbkern p_U$.
\item $P_{t}^{t-1}\circ P_{t}^{t-1} = P_{t}^{t-1}$.
\item $P_{t}^{t-1}$ ist stetig mit $\norm{p}{P_{t}^{t-1} \lambda} \le \norm{p}{\lambda}$, $p\in \N$.
\end{enumerate}

ad 1): Sei $\lambda=(\lambda_l)_l\in \abbkern p_U$, das heißt $\ffa{l}{\infty} \lambda_l f_l =0$.
Es ist 
\formelohne{  \ffa{l}{a_t} \lambda_l f_l &=& \ffa{l}{b_{t-1}} \lambda_l f_l + \ffa{i}{a_t - b_{t-1}} 
              \lambda_{b_{t-1}+i} f_{b_{t-1}+i} \\
             &=&   \ffa{l}{b_{t-1}} \lambda_l f_l  + \ffa{i}{a_t - b_{t-1}} \lambda_{b_{t-1}+i}
                   \left( \ffa{l}{b_{t-1}} c_{t-1,i,l } f_l\right) \\
             &=&  \ffa{l}{b_{t-1}} \left( \lambda_l + \ffa{i}{a_t - b_{t-1}} \lambda_{b_{t-1}+i}
                   c_{t-1,i,l} \right) f_l \\
             &=&   \ffa{l}{b_{t-1}} \mu_l f_l     }
Andererseits ist 
\[  \ffa{l}{a_t} \lambda_l f_l = - \sum_{l=a_t+1}^{\infty} \lambda_l f_l =   - \sum_{l=a_t+1}^{\infty}
       \mu_l f_l = - \sum_{l=b_{t-1}+1}^{\infty} \mu_l f_l   ,    \]
also ist $P_{t}^{t-1}(\lambda) \in \abbkern p_U $.

Die 2. Aussage ist klar.

ad 3): \formelohne{ \norm{p}{P_{t}^{t-1}(\lambda)} &=& \ffa{l}{b_{t-1}} \ffb{A_p}{\lambda_l + 
                       \ffa{i}{a_t - b_{t-1}} \lambda_{b_{t-1}+i} c_{t-1,i,l} }  + 
                             \sum_{l=a_t+1}^{\infty} \ffb{A_p}{\lambda_l}    \\
                    &\le&    \ffa{l}{b_{t-1}} \ffb{A_p}{\lambda_l} +  \sum_{l=b_{t-1}+1}^{a_t}
                             \ffb{A_p}{\lambda_l} +  \sum_{l=a_t+1}^{\infty}\ffb{A_p}{\lambda_l}  \\
                    &=&    \norm{p}{\lambda},    }
da $\ffb{A_p}{c_{t-1,i,l}} <1$ für $t\ge s+1$, $i=1,\ldots,a_t - b_{t-1}$, $l=1,
\ldots, b_{t-1}$ .

Definiere $C_i^t \in l_1(H_U(\oka)) $ für $i=1,\ldots,a_{t+1} - b_{t}$, $t\ge s+1$,
durch die Komponenten 
\[   \left[ C_i^t\right]_l := \left\{ \begin{array}{l@{\leer{0.5} ; \leer{0.5}}l} 
                 c_{t,i,l} & l=1,\ldots, b_t \\
                 -1 & l=b_t +i \\
                 0  & \mbox{sonst} \end{array}  \right.   \]
{\sc Bemerkung:} Es gilt
\[   \ffa{l}{b_{t+1}} [ C_i^t]_l f_l=0 \]
auf einer Umgebung von $K_t$, denn nach Voraussetzung gilt:
\[   \ffa{l}{b_{t+1}} [ C_i^t]_l f_l=\ffa{l}{b_{t}} c_{t,i,l} f_l + (-1)f_{b_t+i} =0\punkt\]
Soweit die Bemerkung.

Es gilt 
\[   \abbkern P_{t}^{t-1} =\left\{ \sum_{i=1}^{a_t-b_{t-1}} \alpha_i C_i^{t-1}\doppelpunkt \alpha_i\in 
        H_U(\oka) \right\}  ,\]
denn ist $\alpha\in H_U(\oka)$, dann gilt für $i=1,\ldots,a_t - b_{t-1}$
\formelohne{   
          \left[ P_{t}^{t-1}(\alpha C_i^{t-1})\right]_l 
               &=& \left\{
                \begin{array}{l@{\leer{0.5} ; \leer{0.5}}l}  
                   \alpha c_{t-1,i,l} + \alpha(-c_{t-1,i,l}) & l\le b_{t-1}  \\
                    0 & \mbox{sonst} \end{array} \right.  \\
               &=&  0  \punkt         }
Sei andererseits $ P_{t}^{t-1} \lambda =0$, dann folgt $\lambda_l=0 $ für $l\ge a_t+1$.
Für $l\le b_{t-1}$ gilt 
\formelohne{ 
            \lambda_l &=& - \sum_{i=1}^{a_t-b_{t-1}} \lambda_{b_{t-1}+i} c_{t-1,i,l} \\
                  &=& \sum_{i=1}^{a_t-b_{t-1}} (-\lambda_{b_{t-1}+i})[ C_i^{t-1}]_l \komma  }
und da für $l=b_{t-1} +i$, $ i=1,\ldots,a_t - b_{t-1}$ 
\[        \lambda_{b_{t-1}+i} = (- \lambda_{b_{t-1}+i} )[ C_i^{t-1}]_l \mbox{ ist}\komma \]
folgt 
\[        \lambda = (\ffc{\lambda}{a_t}) = \sum_{i=1}^{a_t-b_{t-1}}(- \lambda_{b_{t-1}+i} )
                 C_i^{t-1} \punkt \]

Setze  $ Q_{t+n}^t := P_{t}^{t-1}\circ \ldots\circ P^{t+n-1}_{t+n} $. Ist $\lambda\in l_1(H_U(\oka))$ 
fest, dann ist  $( Q_{t+n}^t \lambda)_n$ Cauchyfolge in $l_1(H_U(\oka))$:

Sei dazu $x_n:= Q_{t+n}^t \lambda$, dann ist für $p\in \N$ 
\formelohne{   \norm{p}{x_{n+1} - x_n} 
            &=& \norm{p}{
                 P_{t}^{t-1}\circ \ldots \circ P^{t+n-1}_{t+n} \circ P_{t+n+1}^{t+n} \lambda -
                 P_{t}^{t-1}\circ \ldots \circ P_{t+n}^{t+n-1} \lambda }       \\
            &\le& \norm{p}{P_{t+n+1}^{t+n}\lambda - \lambda }  \\
            &=&  \ffa{l}{b_{t+n}} \ffb{K_p}{\lambda_l +\leer{-0.5} \sum_{i=1}^{a_{t+n+1}-b_{t+n}}
                   \leer{-0.5}\lambda_{b_{t+n}+i} c_{t+n,i,l} -\lambda_l }
                  + \sum_{b_{t+n}+1}^{a_{t+n+1}} \ffb{K_p}{\lambda_l}  \\
            &\le&  \sum_{i=1}^{a_{t+n+1}-b_{t+n}} \ffb{K_p}{\lambda_{b_{t+n}+i}} 
                    \ffa{l}{b_{t+n}} \ffb{K_p}{c_{t+n,i,l}} +  \sum_{b_{t+n}+1}^{a_{t+n+1}}
                           \ffb{K_p}{\lambda_l}  \\
            &<& 2 \sum_{b_{t+n}+1}^{a_{t+n+1}}  \ffb{K_p}{\lambda_l}       \punkt    }

Es folgt 
\[    \ffa{n}{\infty} \norm{p}{x_{n+1} -x_n} \le 2 \norm{p}{\lambda} \punkt \]

Wir setzen also $Q_t(\lambda):= \lim_n Q_{t+n}^t(\lambda)$ für $\lambda\in \zzc$. 
Dann ist $Q_t$ eine lineare Abbildung $Q_t:\zzc\rightarrow \zzc$. Es sei
\[      F_t:=\{ \lambda\in \zzc \doppelpunkt \lambda_l=0 \mbox{ für } l\ge b_{t-1}  +1 \}   \]

Eigenschaften von $Q_t$:\be
\item $Q_t$ ist stetig und $\norm{p}{Q_t(\lambda)} \le\norm{p}{\lambda} $ für $\lambda\in \zzc$,
      $p\in \N$.
\item $Q_r(\lambda)=\lambda$ für $\lambda\in F_t$, $r\ge t$.
\item $Q_t(\lambda)\rightarrow \lambda$ für $t\rightarrow \infty$ in $\zzc$.
\item $Q_t(\abbkern p_U)\subset \abbkern p_U$.
\item $Q_t(\lambda) \in F_t$ für alle $\lambda\in \abbkern p_U$.
\ee

\be
\item[ad 1)] $ \norm{p}{Q_t\lambda} = \lim_n\norm{p}{Q^t_{t+n}(\lambda)} \le \norm{p}{\lambda}$.
\item[ad 2)] $P_{t}^{t-1} \lambda = \lambda$, falls $\lambda_l=0$ für $l\ge b_{t-1}+1$. Daher
ist $Q_{r+n}^r \lambda =\lambda$, falls $\lambda\in F_t$, $r\ge t$ und $n\in \N_0$ .
\item[ad 3)] Sei $p\in \N$ und $\eps>0$, dann gibt es ein $t$, so daß sich $\lambda$ für 
$r\ge t$ schreiben läßt als $\lambda=\lambda_1+ \lambda_2$ mit $\lambda_1\in F_r$ und $\norm{p}{\lambda_2}
<\frac{\eps}{2}$. 

Es folgt für jedes $r\ge t$ 
\formelohne{
   \norm{p}{Q_r(\lambda)-\lambda} &=& \norm{p}{Q_r(\lambda_1 + \lambda_2) - \lambda}  \\
                          &=& \norm{p}{\lambda_1 + Q_r(\lambda_2) - \lambda}  \\
                          &=& \norm{p}{Q_r(\lambda_2) - \lambda_2}        \\
                          &\le& \norm{p}{\lambda_2} + \norm{p}{\lambda_2}\\
                          & <& \eps   \punkt }
{\sc Bemerkung:} Da $F:=\{ \lambda\in \zzc:\lambda\in F_r \mbox{ für ein }r\}$
eine in $\zzc$ totale Menge ist, folgt 3) direkt aus dem Satz von Banach-Steinhaus 
[Köthe Bd 2, S.142].
\item[ad 4)] Ist $\lambda\in \abbkern p_U$, dann ist $Q_{t+n}^t(\lambda)\in \abbkern p_U$. Da
$\abbkern p_u$ in $\zzc$ abgeschlossen ist, folgt $Q_t(\lambda)\in \abbkern p_U$.
\item[ad 5)] Sei $r\ge t$, $t+n>r$, $\mu^n:= Q_{t+n}^t \lambda$, $\lambda\in\abbkern p_U$ . Setze 
\[   \mu' := P^{r}_{r+1} (P_{r+2}^{r+1}\circ \ldots \circ P_{t+n}^{t+n-1}\lambda) \komma \]
dann ist $[\mu']_{b_r+i}=0 $ für $i=1,\ldots,a_{r+1} - b_{r}$. Andererseits ist 
$\mu^n = P_{t}^{t-1}\circ \ldots \circ P_{r}^{r-1}\mu'$ und damit $\mu_l^n=\mu'_l$ für
$l\ge a_r +1$, also folgt $\mu_{b_r+i}^n=0$ für $i=1,\ldots,a_{r+1} - b_{r}$.

Sei $\mu:= Q_t(\lambda) = \lim_n  Q^t_{t+n}(\lambda) = \lim_n \mu^n$, dann gilt 
$\mu_l^n\rightarrow \mu_l$ für $n\rightarrow \infty$ in $H_U(\oka)$, also
ist $\mu_{b_r+i}=0$ für $r\ge t$, $i=1,\ldots,a_{r+1} - b_{r}$. 
Nach Voraussetzung \auslabeln{7.0.3} gilt dann aber bereits $\mu_l=0$ für $l\ge a_t+1$,
daher ist $\mu \in F_t$. 
\ee

Wir behalten  die Situation $U\relkom K_s$ bei und setzen für $\lambda\in \abbkern p_U$
\[    x_t:= \left\{  \begin{array}{l@{\leer{0.5} ; \leer{0.5}}l}  
                   Q_s(\lambda) & t=s-1 \\
                   Q_{t+1}(\lambda) - Q_t(\lambda) & t\ge s   \end{array}\right. \punkt  \]
Es gilt $\sum\limits_{t=s-1}^{\infty} x_t= \lambda$ da $Q_t(\lambda)\rightarrow \lambda$ 
für $t\rightarrow \infty$ in $\zzc$.

Für $t\ge s$ gilt 
\formelohne{   P_{t}^{t-1} x_t 
                   &=& P_{t}^{t-1} (\lim_n Q^{t+1}_{t+n}  (\lambda) - 
                          \lim_n Q^{t}_{t+n}(\lambda)) \\
                   &=& \lim_n( P_{t}^{t-1}\circ Q^{t+1}_{t+n}(\lambda) -
                               P_{t}^{t-1}\circ Q^{t}_{t+n}(\lambda))  \\
                   &=& 0 \punkt          }
Also ist $x_t\in \abbkern P_{t}^{t-1} = \{ \sum\limits_{i=1}^{a_{t}-b_{t-1}} \alpha_i C_i^{t-1}
\komma \alpha_i\in H_U(\oka)\}$. 

Es gibt also $\alpha_{i,t}\in H_U(\oka)$, $i=1,\ldots,a_{t+1} - b_{t} $, so daß 
\[      x_t= \sum_{i=1}^{a_{t}-b_{t-1}} \alpha_{i,t}  C_i^{t-1}  \punkt \] 

Daher ist $\lambda-Q_s(\lambda) = \sum\limits_{t=s}^{\infty} \sum\limits_{i=1}^{a_{t}-b_{t-1}}
\alpha_{i,t}  C_i^{t-1} $.

Zeige nun, daß $ \sum\limits_{t=s}^{\infty} \sum_{i=1}^{a_{t+1}-b_{t}} \ffb{K_p}{
\alpha_{i,t} } <\infty$ für $p\in \N$ ist, und damit $(\alpha_{i,t})_{i,t\ge s}
\in \zzc$ gilt:

Es gibt $D_p\in \R^{+}$, so daß $\norm{p}{\sum\limits_{t=s}^n x_t }< D_p$ gleichmäßig
in $n$ ist. Sei also $n\ge s$ beliebig, dann gilt 
\formelohne{   D_p &>& \norm{p}{\sum_{t=s}^n x_t}  \\
                   &=& \norm{p}{\sum_{t=s}^{n} \sum_{i=1}^{a_{t+1}-b_{t}}
                        \alpha_{i,t}  C_i^t  }   \\
                   &=& \ffa{l}{\infty} \ffb{K_p}{\sum_{t=s}^{n} \sum_{i=1}^{a_{t+1}-b_{t}}
                        \alpha_{i,t}  d_{t,i,l} }     \mbox{ mit}        }  
\[      d_{t,i,l}:= [C_i^t]_l =\left\{ \begin{array}{l@{\leer{0.5} ; \leer{0.5}}l}
                       c_{t,i,l} & l=1,\ldots b_t  \\
                       -1 & l=b_t + i  \\
                       0  &  \mbox{sonst}   \end{array}\right. \punkt   \]

Setze 
\[      \beta_l := \left\{ \begin{array}{l@{\leer{0.5} ; \leer{0.5}}l}
                      \alpha_{i,t} & l=b_t +i\komma 1\le i \le a_{t+1} - b_{t}\komma
                                         s\le t\le n  \\
                      0            &  \mbox{sonst}   \end{array}\right. \punkt   \]
Dann ist 
\[    \ffa{l}{\infty} \ffb{K_p}{\beta_l} = \sum_{t=s}^{\infty} \left(
            \sum_{l=b_t +1}^{a_{t+1}}\ffb{K_p}{\beta_l} + 
            \sum_{l=a_{t +1}+1}^{b_{t+1}}\ffb{K_p}{\beta_l} \right) = 
            \sum_{t=s}^{n} \sum_{i=1}^{a_{t+1}-b_{t}} \ffb{K_p}{\alpha_{i,t}} \punkt   \]
Also ist 
\formelohne{D_p &>& \leer{0}\sum_{l=1}^{\infty} \ffb{K_p}{\sum_{t=s}^{n}
                 \sum_{i=1}^{a_{t+1}-b_{t}} \alpha_{i,t} d_{t,i,l}  }   \\
           &=& \sum_{l=1}^{\infty} \ffb{K_p}{ \beta_l - \sum_{t=s}^{n} 
                     \sum_{{i=1\atop b_t+i\not=l}}^{a_{t+1}-b_{t}} \alpha_{i,t} d_{t,i,l}  }   \\
           &\ge&  \sum_{l=1}^{\infty}\ffb{K_p}{\beta_l} - 
                    \sum_{l=1}^{\infty}\sum_{t=s}^{n} \sum_{{i=1\atop b_t+i\not=l}}^{a_{t+1}-b_{t}}
                    \ffb{K_p}{ \alpha_{i,t}} \ffb{K_p}{  d_{t,i,l}}  \\
           &=& \sum_{l=1}^{\infty}\ffb{K_p}{\beta_l} -  \sum_{t=s}^{n} \sum_{i=1}^{a_{t+1}-b_{t}}
                    \ffb{K_p}{ \alpha_{i,t}} \ffa{l}{b_t} \ffb{K_p}{ c_{t,i,l}} \\
           &\ge&  \sum_{t=s}^{n} \sum_{i=1}^{a_{t+1}-b_{t}}
                    \ffb{K_p}{ \alpha_{i,t}}  -  \rho \sum_{t=s}^{n} \sum_{i=1}^{a_{t+1}-b_{t}}
                    \ffb{K_p}{ \alpha_{i,t}} \mbox{ mit $\rho<1 $ geeignet} \\
           &=& (1-\rho) \sum_{t=s}^{n} \sum_{i=1}^{a_{t+1}-b_{t}}
                    \ffb{K_p}{ \alpha_{i,t}}   \punkt    }


Oben wurde gezeigt, daß $Q_s(\lambda)\in F_s$, da $\lambda\in $ Kern $p_U$. Das bedeutet, daß 
$\ffa{l}{b_s}[Q_s(\lambda)]_l\cdot f_l =0$.

Da $f_j \in \bigoplus_{l\in \N} H_{\Omega}(\F)$ ist für $j=1,2,\ldots$, gibt es ein 
$N_s$, so daß $(f_1)_z,\ldots,(f_{b_s})_z$ einen Untermodul $(\F_s)_z$ von $\F_z^{N_s}$
erzeugen für 
$z\in \Omega$. Es ist $\bigcup\limits_{z\in \Omega} (\F_{s})_z=:\F_{s}$ eine kohärente 
analytische Garbe (\cite{G-R},Chap IV, Sec. B,  Proposition 12). $\F_{s}$ wird also
von $f_1,\ldots,f_{b_s}$ in $\Omega$ global erzeugt (\cite{Hoer}, Theorem 7.2.9).

Ist $\Rscript_{s}:=\Rscript(f_1,\ldots,f_{b_s})$ die Relationsgarbe von
$f_1,\ldots,f_{b_s}$, 
dann gibt es zu jedem offenen $\omega\relkom\Omega$ 
endlich viele globale Schnitte, die $\Rscript_{s}$ in $\omega$ erzeugen ([Grauert
Fritzsche] ).

Wir konstruieren eine Folge $(\overline{g}_j)_j$ in $\zzc$ wie folgt:

$\ffc{\overline{g}}{m_0}$ seien die in $\zzc$ kanonisch eingebetteten globalen
Schnitte, die $\Rscript_0$ in $\stackrel{\circ}{K}_0$ erzeugen.
%
%
%
$\overline{g}_{m_{s-1}+1},\ldots,\overline{g}_{m_{s}}$ seien die in $\zzc$
kanonisch eingebetteten globalen Schnitte, die $\Rscript_s$ in 
$\stackrel{\circ}{K}_s$ erzeugen für $s\ge 1$. Ohne Einschränkung können wir annehmen,
daß $m_t$ streng monoton in $t$ wächst.

Für jedes $\lambda\in \abbkern p_U$ liegt $Q_s(\lambda)$ in dem Kern der Abbildung 
\[   \pi : H_U(\oka^{b_s}) \rightarrow H_U(\F_s)\komma \ffc{c}{b_s} \mapsto 
                 \sum_j c_j f_j \punkt  \]
Nach Lemma \auslabelnsatz{7.2}    ist 
\[   H_U(\oka^{m_s}) \rightarrow \abbkern \pi \komma  \ffc{c}{m_s} \mapsto 
                 \sum_j c_j \overline{g}_j          \]
eine stetige lineare surjektive Abbildung.

Insgesamt folgt:

Ist $U\subset K_s$ offen holomorph konvex, dann ist die Abbildung
\formelohne{  
    \tilde{q}_U: H_U(\oka^{m_s})\times \zzc &\rightarrow& \zzc   \\
    \left( (\ffc{\beta}{m_s}),(\alpha_{i,t})_{t\ge s\komma 1\le i\le a_{t} - b_{t-1}}
            \right) &\mapsto& \ffa{j}{m_s}\beta_j \overline{g}_j + \sum_{t\ge s}
                                 \sum_{i=1}^{a_{t}-b_{t-1}} \alpha_{i,t} C_i^t     }
stetig linear mit $\bild \tilde{q}_U =\abbkern p_U$.

{\sc Bemerkung:} Ist $\lambda\in \abbkern p_U\cap F_r$ für ein $r\ge s$, dann gilt
$Q_{t+1}(\lambda)=Q_t(\lambda) = \lambda$ für $t\ge r$. Daher gilt 
\formel{{\einlabeln{7.0.5}}   \lambda=\sum_{t=s-1}^{\infty} x_t
         &=& \sum_{t=s-1}^{r-1} x_t = \ffa{j}{m_s}\beta_j 
    \overline{g}_j + \sum_{t=s}^{r-1}   \sum_{i=1}^{a_{t}-b_{t-1}} \alpha_{i,t} C_i^{t-1} \punkt 
}

Die erzeugenden Schnitte $C_i^t$ aus der Abbildung $\tilde{q}_U$ haben den Nachteil, 
daß sie lediglich jeweils auf einer Umgebung
von $K_s$ definiert sind für $t\ge s$.

Wir zeigen im nächsten Schritt, daß es globale Schnitte $G_i^t$, $i=1,\ldots,a_{t+1}
 - b_{t}$, $t\ge 0$ gibt, so daß   
\formelohne{  
    \tilde{\tilde{q}}_U: H_U(\oka^{m_s})\times \zzc &\rightarrow& \zzc   \\
    \left( (\ffc{\beta}{m_s}),(\alpha_{i,t})_{t\ge s\komma 1\le i\le a_{t+1} - b_{t}}
            \right) &\mapsto& \ffa{j}{m_s}\beta_j \overline{g}_j + \sum_{t\ge s}
                                 \sum_{i=1}^{a_{t+1}-b_{t}} \alpha_{i,t} G_i^t     }
stetig linear ist mit $\bild \tilde{\tilde{q}}_U= \abbkern p_U$ für $U\relkom K_s$ 
offen, holomorph konvex.

Hierzu benötigen wir die folgenden zwei Lemmata:

\begin{satz}{Lemma}\einlabelnsatz{7.5}
Sei $E$ ein Fr\'echetraum mit einem Fundamentalsystem von stetigen Halbnormen $(\norm{p}{\leer{0.5}})_p$
und $\Phi:E\rightarrow E$ stetig linear. Ferner gelte für die Identität $Id$ auf $E$, daß
$\norm{p}{ (Id -\Phi)(x)} \le \theta_p\norm{p}{x}$ mit $\theta_p<1$ für 
jedes $p\in \N$. 

Dann ist $\Phi $ invertierbar in L($E$) und es gilt 
\[   \norm{p}{\Phi^{-1}x}\le \frac{1}{1-\theta_p} \norm{p}{x}  \]
\end{satz}

Es sei \[ \Psi_n:=\sum_{k=0}^{n} (Id -\Phi)^k \mbox{ mit } (Id-\Phi)^0:=Id \]
Dann gilt für $m\ge n$ und für alle $p\in \N$ und $x\in E$
\[   \norm{p}{\Psi_m(x)-\Psi_n(x)} = \norm{p}{\sum_{k=n+1}^m (Id -\Phi)^k(x)}
               \le \sum_{k=n+1}^m \theta_p^k \norm{p}{x}  \]
Hieraus folgt insbesondere $\norm{p}{\Psi_m(x)}\le \frac{1}{1-\theta_p}\norm{p}{x} $.

Nach dem Satz von Banach-Steinhaus (siehe z.B. \cite{Koethe}, §39, Nr. 5) gibt es ein 
$\Psi\in {\rm L}(E)$ mit $\Psi_n\rightarrow \Psi$ punktweise. 
Es folgt daher wie beim von Neumannschen Lemma für jedes $x\in E$

\[   \Psi \circ \Phi(x) = \lim_n \Psi_n\circ \Phi(x) =x  \]
Genauso folgert man wegen der Stetigkeit von $\Phi$, daß $\Phi\circ\Psi(x)=x$.

Wir setzen also $\Phi^{-1}:=\Psi$. Für $p\in \N$ gilt 
\[   \norm{p}{\Phi^{-1}} = \lim_n\norm{p}{\Psi_n(x)} \le 
      \frac{1}{1-\theta_p} \norm{p}{x}  \leer{1} \Box  \]

\begin{satz}{Lemma}\einlabeln{7.6}
Sei $(E,(\norm{p}{\leer{0.5}})_p)$ eine Fr\'echetalgebra mit Eins und 
$A=(a_{i,j})_{i,j}$ eine 
unendliche Matrix mit Einträgen aus $E$, so daß $A$ eine  Abbildung
aus L$(l_1(E))$ repräsentiert. Ferner gebe es eine aufsteigende
Folge positiver ganzer Zahlen $(n_m)_m\in \N_0$ mit $n_0=0$, so daß
\[   a_{i,j}=0 \mbox{ für } i>n_{m(j)}\komma
                \mbox{ falls } n_{m(j)-1}<j\le n_{m(j)} \]

Ferner gelte für die Einheitsmatrix $I$ bestehend aus dem Einselement
der Algebra auf der Diagonalen, daß für jedes $v\in l_1(E)$
\[   \norm{p}{ (I -A)(v)} \le \theta_p\norm{p}{v} \mbox{ mit }
 \theta_p<1 \mbox{ für  jedes } p\in \N \]

Dann ist $A$ invertierbar auf $l_1(E)$ und für $A^{-1}=(b_{i,j})_{i,j}$
gilt   

\[   b_{i,j}=0 \mbox{ für } i>n_{m(j)}\komma 
         \mbox{ falls } n_{m(j)-1}<j\le n_{m(j)} \]
\end{satz}

Sei $c_{i,j}^k$ das $i,j$-te Element von $(I-A)^k$, dann zeigt man
wie folgt induktiv über $k$, daß 
\[   c_{i,j}^k=0 \mbox{ für } i>n_{m(j)}\komma 
       \mbox{ falls } n_{m(j)-1}<j\le n_{m(j )}\]
Der Induktionsanfang $k=1$ ist nach der Voraussetzung an $A$ klar. Für
den Induktionsschritt $k\rightarrow k+1$ sei zu festem $j$ ein Index $m(j)$ 
so gewählt, daß
$n_{m(j)-1}<j\le n_{m(j)}$ ist.  Für ein festes $i$ gibt es ein $l\in \N$, so daß
$n_{l-1}<i\le n_l$. Dann gilt
\formelohne{  c_{i,j}^{k+1} &=& \sum_{\nu=1}^{\infty} 
              c_{i,\nu}^k c^1_{\nu,j} \\
       &=&    \sum_{\nu=n_{l-1} +1}^{n_{m(j)}} c_{i,\nu}^k c^1_{\nu,j}  }
Damit ist $ c_{i,j}^{k+1}=0$, falls $l-1\ge m(j)$, was der Fall ist, wenn $i>n_{m(j)}$.

Ist $b_{i,j}^n$ das $i,j$-te Element von $A_n:=\sum_{k=0}^n (I-A)^k$, dann folgt
\[   b^n_{i,j}=0 \mbox{ für } i>n_{m(j)}\komma 
         \mbox{ falls } n_{m(j)-1}<j\le n_{m(j)} \]
Sei $e_j$ der $j$-te Einheitsvektor in $l_1(E)$. Dann ist das 
$i$-te Element des Vektors $A_n e_j$ gerade $b_{i,j}^n$. 
Da $(A_n e_j)_n$ eine Cauchyfolge in 
$l_1(E)$ ist, die gegen $A^{-1} e_j$ konvergiert, verschwindet für $i>n_{m(j)}$ auch das 
$i$-te Element des Vektors $A^{-1} e_j = b_{i,j}$ $\Box$

Wir kommen zurück zum Beweis der Proposition \auslabelnsatz{7.3}:

{\sc Bemerkung:} Ist $d\in \Gamma(\omega,\chi)$ ein Schnitt auf einer Umgebung $\omega$
von $K_s$ in eine kohärente analytische Garbe $\chi$, dann gibt es
 eine Folge $(d_j)_j\subset\Gamma(\Omega,\chi)$, so daß $\norm{K_s}{d-d_j}
\rightarrow 0$, wenn $j\rightarrow \infty$, da $K_s$ holomorph
konvex ist (\cite{Hoer}, Theorem 7.2.7).

Ist  $(D_k^t)_{k=1}^{k_t}$ für $t\in \N_0$ eine Folge von Schnitten in $\chi$, so daß
$\ffc{D^t}{k_t}$ jeweils $\chi$ um $K_t$ erzeugt, und ist für eine  offene Umgebung
$\tilde{\omega}$ von $K_t$
\[     \norm{K_t}{x}:= \inf_{x=\ffa{k}{k_t} h_k D_k^t}  \ffa{k}{k_t} \ffb{K_t}{h_k} 
       \mbox{ für } x\in \Gamma(\tilde{\omega},\chi)\komma  \]
dann gibt es zu vorgegebenem $\eps >0$ ein
$J$, so daß $\norm{K_t}{d-d_j}<\eps$ für $j\ge J$, $0\le t \le s$. Soweit die 
Bemerkung.

Sei $\rho<1$ fest. $C_i^t$ ist ein Schnitt in $\Rscript_{t+1}$ auf einer Umgebung von 
$K_t$. Wenden wir \auslabeln{7.0.5} auf $U=\stackrel{\circ}{K}_\tau$, $\tau\le t$ an,
so wird $\Rscript_{t+1}$ in $\stackrel{\circ}{K}_\tau$ durch
\[   M_\tau:=\{ \ffc{\overline{g}}{m_\tau},(C_i^l)_{l=\tau,\ldots,t,i=1,\ldots,a_{l+1}-b_l} \} \]
erzeugt. Nach Theorem 7.2.9 in \cite{Hoer}, einer unmittelbaren Folgerung aus
Theorem A von Cartan, wird  $\Rscript_{t+1}$ auch in pseudokonvexen Umgebungen 
 von $K_{\tau-1}$, die in $\stackrel{\circ}{K}_\tau$ liegen,   durch $M_\tau$ erzeugt.
Sei für  einen Schnitt $x$ in $\Rscript_{t+1}$ um $K_\tau$ 
\[    \norm{s,\tau}{x} := \inf\{\sum_{j=1}^{m_{s+1}} \ffb{K_s}{\beta_j} + 
                         \sum_{l=s+1}^t \sum_{i=1}^{a_{l+1}-b_{l}} 
                                               \ffb{K_s}{\alpha_{i,l}}\doppelpunkt
x=\sum_{j=1}^{m_{s+1}} \beta_j \overline{g}_j + 
                    \sum_{l=s+1}^t \sum_{i=1}^{a_{l+1}-b_{l}} \alpha_{i,l} C_i^l \}
                                                                 \punkt  \]
Dann gibt es entsprechend der obigen Bemerkung jeweils
 ein $G_i^t\in \Gamma(\Omega,\Rscript_{t+1})$, so daß
\formel{   \norm{s}{C_i^t - G_i^t} < \frac{\rho}{2}\leer{1} \mbox{ für } 
              0\le s\le t-1\komma   i=1,\ldots,a_{t+1} - b_{t}     \einlabeln{8.0.0.1}}
Auf diese Weise konstruieren wir für jedes $t\ge 0$ jeweils $G_i^t$ 
für $i=1,\ldots,a_{t+1} - b_{t}$.

Ist $s$  fest, dann gibt es somit für jedes $\tau\ge s$ jeweils um $K_{s-1}$ eine Darstellung
\[   G_j^\tau - C_j^\tau = \ffa{l}{m_s} \beta_l^{\tau,j} \overline{g}_l + 
                     \sum_{t=s}^\tau  \sum_{i=1}^{a_{t+1}-b_{t}} \alpha_{t,i}^{\tau,j} C_i^t  \]
mit 
\[   \ffa{l}{m_s}\ffb{K_{s-1}}{\beta_l^{\tau,j}} + \sum_{t=s}^\tau  \sum_{i=1}^{a_{t+1}-b_{t}} 
                 \ffb{K_{s-1}}{\alpha_{t,i}^{\tau,j}} < \rho    \]


Sei eine Abbildung $H:\zzc^2\rightarrow \zzc^2$ durch $H(\mu^1,\mu^2):=(H_1(\mu^1,\mu^2),
H_2(\mu^1,\mu^2)) $ definiert, wobei wir
\[       \left[ H_1(\mu^1,\mu^2)\right]_l := \left\{ \begin{array}{c@{\leer{0.5} ; \leer{0.5}}l} 
           \mu_l^1 + \sum_{\tau=s}^{\infty} \sum_{j=1}^{a_{\tau+1}-b_\tau}
                                           \mu_{b_\tau+j}^2 \beta_l^{\tau,j}   & l=1,\ldots, m_s \\
                 0                                                       & \mbox{sonst}
             \end{array} \right.           \]
setzen 
und $ H_2(\mu^1,\mu^2)= A\mu^2$ mit einer Matrix $A=(a_{l,n})_{l,n}$. Hierbei seien die 
Matrixeinträge wie folgt definiert:
\[    a_{l,n} := \left\{  \begin{array}{c@{\leer{0.5} ; \leer{0.5}}l}
                \alpha_{t,i}^{\tau,j} & l=b_t +i\komma s\le t\le \tau\komma i=1,
                           \ldots,a_{t+1} - b_{t}\komma l\not= b_\tau+j\komma n=b_\tau +j  \\
             1+ \alpha_{t,i}^{\tau,j} & l=n=b_\tau +j\komma t=\tau\komma i=j   \\
                        0          & \mbox{sonst}   \end{array} \right.  \punkt         \]
$H$ hängt natürlich von $s$ ab. Wir schreiben aus Vereinfachungsgründen $H$ statt $H_s$.

Wir zeigen \be
\item $\norm{p}{H_1(\mu^1,\mu^2)}\le \norm{p}{\mu^1} + \norm{p}{\mu^2}$ für alle $p\in \N$.
\item $A$ ist eine Abbildung in L$(l_1(H_U(\oka)))$. Es gilt $\norm{p}{A\lambda} \le
      2\norm{p}{ \lambda}$ für alle $p\in \N$, $\lambda\in l_1(H_U(\oka))$.
      $A$ ist eingeschränkt als Abbildung auf 
\[   N_s := \{ \mu\in \ffd{U}{\oka} \doppelpunkt \mu_1=\ldots  \mu_{b_s}=0\komma
                  \mu_{a_\tau +j}=0\komma \tau\ge 0\komma j=1,\ldots, b_{\tau} - a_\tau  \}  \]
      in sich invertierbar.

 
\item $H$ ist eine stetige lineare Abbildung und $\norm{p}{H(\mu^1,\mu^2)}\le 3(\norm{p}{
      \mu^1} + \norm{p}{\mu^2})$ für $p\in \N$.
\ee
\be
\item[ad 1)] 
\formelohne{  \norm{p}{ H_1(\mu^1,\mu^2)} &=& \sum_{l=1}^{m_s} \ffb{A_p}{\mu_l^1 +
                  \sum_{\tau\ge s} \sum_{j=1}^{a_{\tau+1}-b_\tau}
                                           \mu_{b_\tau+j}^2 \beta_l^{\tau,j}   }  \\
           &\le&  \sum_{l=1}^{m_s} \ffb{A_p}{\mu_l^1} + \sum_{\tau\ge s} 
                     \sum_{j=1}^{a_{\tau+1}-b_\tau}  \ffb{A_p}{\mu_{b_\tau+j}^2} 
                   \sum_{l=1}^{m_s} \ffb{K_{s-1}}{\beta_l^{\tau,j}}    \\
           &\le&  \norm{p}{\mu^1} + \norm{p}{\mu^2}   \punkt   }
\item[ad 2)] Es ist $A\mu\in N_s$ für
   $\mu\in N_s$.  Für $\lambda\in N_s$ gilt
\formelohne{ \norm{p}{(Id_{N_s} - A)\lambda} &=& \sum_{t=s}^{\infty} \sum_{i=1}^{a_{t+1}-b_{t}}
             \ffb{A_p}{ \sum_{\tau=t}^{\infty} \left( 
              \sum_{j=1\atop {i\not=j\komma \!falls \atop t=\tau}}^{a_{\tau+1}-b_{\tau}}
             -\alpha_{t,i}^{\tau,j} \lambda_{b_{\tau}+j}\right)
               + (1-(1-\alpha_{t,i}^{t,i}))\lambda_{b_t+i} }\\
       &\le&  \sum_{t=s}^{\infty} \sum_{i=1}^{a_{t+1}-b_{t}} \sum_{\tau=t}^{\infty}
              \sum_{j=1}^{a_{\tau+1}-b_{\tau}} \ffb{A_p}{\alpha_{t,i}^{\tau,j} } \betrag{
              \lambda_{b_{\tau}+j} }  \\
       &=&    \sum_{\tau=s}^{\infty} \sum_{j=1}^{a_{\tau+1}-b_{\tau}} \left( \sum_{t=s}^{\tau}
              \sum_{i=1}^{a_{t+1}-b_{t}} \ffb{A_p}{\alpha_{t,i}^{\tau,j} }\right) 
              \ffb{A_p}{\lambda_{b_{\tau}+j} } \\
       &<&    \rho \cdot \norm{p}{\lambda}     }
Hieraus folgt einerseits $\norm{p}{A\lambda} = \norm{p}{Id_{N_s}\lambda-(Id_{N_s}-A)\lambda}
\le \norm{p}{\lambda} + \rho \norm{p}{\lambda} \le 2\norm{p}{\lambda}$. Andererseits folgt
aus Lemma \auslabelnsatz{7.6}, daß $A\mid_{N_s}$ invertierbar ist. 

Wir bezeichnen mit $A^{-1}$ die kanonische Fortsetzung von $(A\mid_{N_s})^{-1}$ auf
$l_1(H_U(\oka))$. Für $A^{-1}$ gilt $\norm{p}{A^{-1}\lambda}\le \frac{1}{1-\rho}\norm{p}{
\lambda}$ für jedes $p\in \N$ und jedes $\lambda\in l_1(H_U(\oka))$.

\item[ad 3)] 3) folgt direkt aus 1) und 2).
\ee

Wenn wir den Definitionsraum von $\tilde{q}_U$ und $\tilde{\tilde{q}}_U$ jeweils
kanonisch in  $(l_1(H_U(\oka)))^2$ einbetten, dann
 gilt $\tilde{q}_U(H(\mu^1,\mu^2)) = \tilde{\tilde{q}}_U
(\mu^1,\mu^2)$.

Berechne für $\lambda\in \zzc$ mit $\lambda_{\tau,j}:=[\lambda]_{b_\tau + j}$, $\tau\ge s$,
$j=1,\ldots,a_{\tau+1} - b_{\tau}$:
\formelohne{  
           [A \lambda ]_l
         &=& \left\{  \begin{array}{c@{\leer{0.5} ; \leer{0.5}}l}
                  \sum\limits_{\tau>t} \sumlim_{j=1}^{a_{\tau+1}-b_\tau} a_{l,b_\tau+j}
                                    \lambda_{\tau,j}
               & l=b_t +i\komma t\ge s\komma i=1,\ldots,a_{t+1} - b_{t} \\
              0 & \mbox{sonst} \end{array}
              \right.   \\
         &=&  \left\{  \begin{array}{c@{\leer{0.5} ; \leer{0.5}}l}
              \begin{array}{l}
              \sumlim_{\tau>t} \sumlim_{j=1}^{a_{\tau+1}-b_\tau} \alpha_{t,i}^{\tau,j}
                     \lambda_{\tau,j} \\ 
               +\sumlim_{j=1,i\not=j}^{a_{t+1}-b_t} \alpha_{t,i}^{t,j}\lambda_{t,i} + 
                        (1+\alpha_{t,i}^{t,i})\lambda_{t,i} 
               \end{array}
                 &   l=b_t +i\komma t\ge s\komma i=1,\ldots,a_{t+1} - b_{t} \\
               0 & \mbox{sonst} \end{array} 
              \right.   \\
         &=&   \left\{  \begin{array}{c@{\leer{0.5} ; \leer{0.5}}l}
              \left(\sumlim_{\tau\ge t} \sumlim_{j=1}^{a_{\tau+1}-b_\tau} \alpha_{t,i}^{\tau,j}
                     \lambda_{\tau,j}\right) + \lambda_{t,i}       
                 &   l=b_t +i\komma t\ge s\komma i=1,\ldots,a_{t+1} - b_{t} \\
               0 & \mbox{sonst} \end{array} 
              \right.   }

Daher ist mit $\mu^2_{\tau,j}:=[\mu^2]_{b_\tau +j}$
\formelohne{ \tilde{q}_U(H(\mu^1,\mu^2)) 
      &=& \sum_{l=1}^{m_s} [H_1(\mu^1,\mu^2)]_l \overline{g}_l + 
           \sum_{t\ge s} \hspace{-0.75em}\sum_{i=1}^{a_{t+1}-b_{t}}\hspace{-0.75em} [A\mu^2]_{l=b_t+i} C_i^t \\
      &=& \sum_{l=1}^{m_s}(\mu_l^1 + \sum_{\tau \ge s} \hspace{-0.75em}
                       \sum_{j=1}^{a_{\tau+1}-b_\tau}\hspace{-1em}
            \mu^2_{\tau,j} \beta_l^{\tau,j}) \overline{g}_l 
          + \sum_{t\ge s }\hspace{-0.75em}\sum_{i=1}^{a_{t+1}-b_{t}}
            (\sum_{\tau\ge t} \hspace{-0.75em}\sum_{j=1}^{a_{\tau+1}-b_\tau}\hspace{-1em} \alpha_{t,i}^{\tau,j}
                                          \mu^2_{\tau,j} + \mu_{t,i}^2) C_i^t  \\
      &=& \sum_{l=1}^{m_s} \mu_l^1 \overline{g}_l +
          \sum_{\tau \ge s} \hspace{-0.75em}\sum_{j=1}^{a_{\tau+1}-b_\tau}\hspace{-1em}\mu^2_{\tau,j}
          \sum_{l=1}^{m_s} \beta_l^{\tau,j} \overline{g}_l
          + \sum_{\tau \ge s} \hspace{-0.75em}\sum_{j=1}^{a_{\tau+1}-b_\tau}\hspace{-1em}\mu^2_{\tau,j}
          \sum_{t=s}^\tau \hspace{-0.4em}\sum_{i=1}^{a_{t+1}-b_{t}}\hspace{-1em}  \alpha_{t,i}^{\tau,j} C_i^t 
          + \sum_{\tau \ge s} \hspace{-0.75em}\sum_{j=1}^{a_{\tau+1}-b_\tau}\hspace{-1em}\mu^2_{\tau,j} C_j^\tau  \\
      &=& \sum_{l=1}^{m_s} \mu_l^1 \overline{g}_l +
          \sum_{\tau \ge s} \hspace{-0.75em}\sum_{j=1}^{a_{\tau+1}-b_\tau}\hspace{-1em}\mu^2_{\tau,j} G_j^\tau \\
      &=& \tilde{\tilde{q}}_U(\mu^1,\mu^2)   }

Also ist $\tilde{\tilde{q}}_U: H_U(\oka^{m_s})\times \zzc \rightarrow \zzc$
stetig linear. Bleibt zu zeigen, daß $\bild \tilde{\tilde{q}}_U = \abbkern p_U$:

Sei $x\in \abbkern p_U=\bild \tilde{q}_U$. Finde also $\eta=(\ffc{\eta}{m_s})$,
$\xi=(\xi_j)_j\in \zzc$, so daß 
\[   x=\tilde{q}_U(\eta,\xi)=\ffa{l}{m_s} \eta_l\overline{g}_l + \sum_{t\ge s}\sum_{j=1}^{a_{\tau+1}-b_\tau}
                   \xi_{b_t+i} C_i^t       \]
Ohne Einschränkung ist also $\xi_j=0$ für $j\le b_s$ und $j=a_t+i$, $t\ge s$,
$i=1,\ldots,a_{t+1} - b_{t}$, d.h. $\xi \in N_s$.

Setze $\lambda:= A^{-1}(\xi)$ und 
\[   \beta_l:= \eta_l - \sum_{\tau=s}^{\infty} \sum_{j=1}^{a_{\tau+1}-b_\tau} 
               \lambda_{b_\tau+j} \beta_l^{\tau,j}\komma l=1,\ldots, m_s\komma
               \beta:=(\ffc{\beta}{m_s})  \]
Dann ist  $H_1(\beta,\lambda)=\eta$ und $H_2(\beta,\lambda)=A\lambda= \xi $.
Also ist 
\[   \tilde{\tilde{q}}_U(\beta,\lambda)=\tilde{q}_U(H(\beta,\lambda))=
          \tilde{q}_U(\eta,\xi) =x\punkt \] 
Ist umgekehrt $x\in \bild \tilde{\tilde{q}}_U$, dann ist $x\in \bild \tilde{q}_U=\abbkern p_U$.

Im nächsten Schritt konstruieren wir aus den Folgen $(\overline{g}_j)_j$ und
$(G_i^t)_{i,t}$ eine Folge $(g_l)_l$ aus $\bigoplus\limits_{l\in \N} H_{\Omega}(\oka)$,
die $\abbkern p_U$ in der gewünschten
Weise erzeugt und die Eigenschaften \auslabeln{7.0.2} und \auslabeln{7.0.3}
erfüllt.

Wir setzen unter Verwendung der obigen Indexfolge $(m_t)_{t\ge0}$ 
\[ \begin{array}{l@{\leer{0.5}:=\leer{0.5}}l@{\leer{0.2},\leer{0.5}}
                                           l@{\leer{0.2},\leer{0.5}}l}
       g_{m_t+n_{t+1}+i} &  \overline{g}_{m_t+i} & i=1,\ldots,m_{t+1}-m_t & t\ge -1 \\
       g_{m_{t+1}+n_{t+1}+i} & G_i^{t+1} & i=1,\ldots,n_{t+2}-n_{t+1} & t\ge -1\komma
   \end{array}   \]
 
wobei $m_{-1}:=0$ und $n_t:=\sum_{\tau=1}^t a_\tau-b_{\tau-1}$ für $t\ge 1$, sowie
$n_0:=0$.

Für jedes $l\in \N$ ist damit $g_l$ wohldefiniert, denn, da $m_t$ strikt
monoton wachsend ist, gibt es zu jedem $l\in \N$ ein $t\ge -1$, so daß
$m_t+n_{t+1} < l \le m_{t+1} + n_{t+2}$.
Ist also im einen Fall $m_{t} + n_{t+1}<l\le m_{t+1} + n_{t+1}$,
 dann setze $i:=l-(m_{t} + n_{t+1})$.
Ist im anderen Fall  $m_{t+1} + n_{t+1}<l\le m_{t+1} + n_{t+2}$, 
dann setze $i:=l-(m_{t+1} + n_{t+1})$.

Wir wollen Eigenschaft \auslabeln{7.0.2} zeigen. Sei dazu $l\in \N$ und es gelte 
$m_{t} + n_{t+1}<l\le m_{t+1} + n_{t+1}$ für ein $t\ge s\ge -1$, dann gilt das folgende

\begin{satz}{Lemma}\einlabelnsatz{8.0.1} Es gibt für $U\subset K_{s-1}$ Funktionen $(\mu_k)_k$ aus 
$H_U(\oka)$, so daß
\[   g_l=\sum_{k=1}^{m_t+n_{t+1}} \mu_k g_k  \]
\end{satz} 
Es ist $g_l=\overline{g}_{m_t+i}$ für ein $i\in \{1,\ldots,m_{t+1}-m_t\}$
Also ist $g_l\in \Gamma(\Omega,\Rscript_{t+1})$. Ist $F_t:=\{\lambda\in \zzc
\doppelpunkt \lambda_k=0 \mbox{ für } k\ge b_{t-1} +1\}$, dann gilt natürlich auch
$g_l\in \abbkern p_U\cap F_{t+2}$. Für Elemente aus $\abbkern p_U\cap F_{t+2}$ wurde gezeigt,
daß die folgende Darstellung auf $U$ existiert:
\[   g_l= \ffa{k}{m_{s-1}} \tilde{\beta}_k \overline{g}_k + \sum_{\tau=s-1}^{t+1}
          \sum_{j=1}^{a_{\tau}-b_{\tau-1}} \alpha_{j,\tau} C_j^{\tau-1}  \]
Sei $\alpha=(\alpha_l)_l$ eine Folge mit 
$\alpha_{b_{\tau-1}+j}=\alpha_{j,\tau}$ für $\tau=s-1,\ldots,t+1$, $j=1,\ldots,
\alpha_{\tau} - b_{\tau-1}$ und $\alpha_l=0$ sonst, dann ist $\alpha\in N_{s-1}$. 
Setze nun $\lambda:=A^{-1}\alpha$ und

\[   \beta_k:= \tilde{\beta}_k - \sum_{\tau=s-1}^{t+1} \sum_{j=1}^{a_{\tau+1}-b_\tau}
        \lambda_{b_\tau +j} \beta_l^{\tau,j}\komma k=1,\ldots,m_s\komma \beta_k=0 \mbox{ sonst}. \]
Hieraus folgt $\tilde{\beta}=H_1(\beta,\lambda)$ und $\alpha=H_2(\beta,\lambda)$. Daher gilt
\[   g_l= \tilde{q}_U(H(\beta,\lambda) = \tilde{\tilde{q}}_U(\beta,\lambda) =
          \ffa{k}{m_{s-1}} \beta_k \overline{g}_k + \sum_{\tau=s-1}^{\infty} 
          \sum_{j=1}^{a_{\tau}-b_{\tau-1}} \lambda_{b_{\tau-1} +j} G_j^{\tau-1}   \]

Da $\alpha= \sum_{\tau=s-1}^{t} \sum_{j=1}^{a_{\tau+1}-b_\tau} \alpha_{b_\tau+j} e_{b_\tau+j}$
ist, folgt mit Blick auf Lemma \auslabelnsatz{7.6} $\lambda=A^{-1} \alpha \in F_{t+1}$.
Damit haben wir
\[   g_l= \ffa{k}{m_{s-1}} \beta_k \overline{g}_k + \sum_{\tau=s-1}^{t+1}
          \sum_{j=1}^{a_{\tau+1}-b_\tau} \lambda_{b_\tau +j} G_j^\tau 
    =\sum_{k=1}^{m_{s-1}+n_{s-1}} \mu_k g_k + \sum_{k=m_{s-1}+n_{s-1} +1}^{m_t+n_{t+1}}
            \mu_k g_k \komma   \]
denn $\{ \overline{g}_k\doppelpunkt k=1,\ldots,m_{s-1}\}\subset 
\{ g_k\doppelpunkt k=1,\ldots, m_{s-1}+
n_{s-1} \}$ und $\{ G_j^\tau\doppelpunkt s-1\le\tau\le t+1\komma j=1,\ldots, a_{\tau+1}-b_\tau\} 
\subset \{ g_k\doppelpunkt k=m_{s-1}+n_{s-1}+1,\ldots,m_{t+1}+n_{t+2} \}$.

Hieraus folgt das Lemma \Ende

Für $t\ge 0$ setzen wir $u_t:=m_{t+1} + n_{t+1}$ und $v_t:=m_{t+1} + n_{t+2}$.
Dann ist $u_t\le v_t\le u_{t+1}$ und $v_{t+1}-v_t\ge m_{t+2}-m_{t+1}>0$ alle
$t\ge 0$. Sei $l\in \N$ mit $m_t+n_{t+1}<l\le m_{t+1} + n_{t+1}$, dann ist
$l=v_{t-1} +i$ mit $i\in \{1,\ldots,m_{t+1}-m_t\} = \{ 1,\ldots,u_t-v_{t-1} \}$.

Nach Lemma \auslabelnsatz{8.0.1} gibt es Funktionen $\mu_{t-1,i,k}$ aus $\zzb$ für
$K_{s-2}\subset U\subset K_{s-1}$ und $s\le t$, so daß 
\[    g_{v_{t-1}+i} = \ffa{k}{v_{t-1}}  \mu_{t-1,i,k} g_k  \punkt \]
Wir können, ohne die Erzeugungseigenschaft der Folge $(g_l)_l$ zu  verändern, die
Folge $(g_l)_l$ bezüglich festem $\rho<1$  wie folgt induktiv skalieren:

Sei $\ffc{\tilde{g}}{v_s}$ bereits skaliert, dann gibt es Funktionen $\tilde{\mu}_{s,i,k}
\in H_U(\oka)$, $K_{s-1}\subset U\subset K_s$, so daß 
\[   g_{v_s+i}=\sum_{k=1}^{v_s} \tilde{\mu}_{s,i,k} \tilde{g}_k \mbox{ auf } U \punkt \]

Sei \formel{{\einlabeln{8.0.1.1}}
    S_{s,i} &:=& \ffa{k}{v_s} \ffb{K_s}{\tilde{\mu}_{t,i,k}} \leer{2} \mbox{für } 
                      i=1,\ldots,u_{s+1}-v_s \mbox{ und} \\
    S_s     &:=& \max_{i=1,\ldots,u_{s+1}-v_s} S_{s,i}\cdot\frac{1}{\rho}
       }
Für $l=v_s +i$ sei $\tilde{g}_l := \frac{g_l}{S_s}$, $i=1,\ldots,u_{s+1}-v_s$. 
Für $l=u_{s+1}+i$, $i=1,\ldots,v_{s+1}-u_{s+1}$ sei $\tilde{g}_l=g_l$.
Dann gilt für die auf diese Weise induktiv skalierte Folge $(\tilde{g}_l)_l$
die Eigenschaft \auslabeln{7.0.2}: Ist $l=v_t +i$, $i\in \{ 1,\ldots,u_{t+1}-v_{t} \}$,
$t\ge s$, dann gibt es in einer Umgebung von $K_{t-1}$ holomorphe Funktionen $\mu_{t,i,k}$,
$k=1,\ldots, v_t$, so daß
\[   \tilde{g}_{v_t+i} = \ffa{k}{v_t} \mu_{t,i,k} \tilde{g}_k   \]
und
\[    \ffa{k}{v_t} \ffb{K_{s-1}}{\mu_{t,i,k}} 
\le \ffa{k}{v_t} \ffb{K_{t-1}}{\mu_{t,i,k}}  \le\rho < 1 \foralls t\ge s \]

Wir benennen von hier an $(\tilde{g}_l)_l$ in $(g_l)_l$ um und definieren für ein
$U\relkom \Omega$ holomorph konvex
\[   q_U :\zzc \to \zzc\komma (\mu_l)_l \mapsto \sum_{l} \mu_l g_l   \]

Wir zeigen jetzt, daß $q_U$ stetig ist:

Da $U\subset K_{s-1}$ gilt für alle $t\ge s$ und $1\le i\le u_t- v_t$
\[   g_{v_t +i} = \sum_{k=1}^{v_t} \mu_{t,i,k} g_k \leer{0,5} \mbox{ mit }
     \sum_{k=1}^{v_t} \ffb{U}{\mu_{t,i,k}} < 1 \]

Sei $(A_p)_{p\in \N}$ eine kompakte Ausschöpfung von $U$, dann folgt für die 
$p$-te Halbnorm des zugehörigen Halbnormensystems in $l_1(H_U(\oka))$ induktiv:

\[   \norm{p}{g_{v_t+i}} \le C_{p,s}:= \max \{ \norm{p}{g_k} \doppelpunkt
      1\le k\le v_s \}   \]
für jedes $t\ge s$ und $1\le i\le u_t- v_t$.

Für $t\ge s$ und $1\le i\le v_{t+1} - u_{t+1}$ gilt $g_{u_{t+1}+i} = G_i^{t+2}$.
Um $K_{s-1}$ sind  $C_i^{t+1}$ definiert für jedes $t\ge s$, $1\le i\le 
v_{t+1}-u_{t+1}=a_{t+3}-b_{t+2}$, und es gilt $\norm{p}{C_i^{t+1}} \le 2$ für
alle $p$.

Wegen \auslabeln{8.0.0.1} in Verbindung mit Lemma \auslabelnsatz{7.0.1} gibt es
$C'_{p,s}$, so daß 

\[   \norm{p}{G_i^{t+1}} \le C_{p,s}' \]
für alle $p$ und alle $t\ge s$, $1\le i\le v_{t+1} - u_{t+1}$ gilt.

Daher gibt es zu jedem $p\in \N$ eine Konstante $C_p$, so daß
\formel{ \norm{p}{q_U(\mu)} \le C_p \norm{p}{\mu}
\leer{0,8}\mbox{ für jedes } \mu\in l_1(H_U(\oka))  \einlabeln{8.0.1.2} }

{\sc Bemerkung:} Es ist $g_l\in \Gamma(\Omega,\Rscript_{t+1})$, für $u_t+1\le l\le u_{t+1}$.

Zum Schluß bleibt nur noch Eigenschaft \auslabeln{7.0.3} zu zeigen:

Sei $\mu=(\mu_l)_l\in \abbkern q_U$ für ein $U\subset K_{s-1}$ holomorph konvex.
Sei $t\ge s$ und $\mu_{v_r+i}=0$ für alle $r\ge t$ und $i=1,\ldots,u_{r+1}-v_r$.
Also ist 
\formelohne{
 0=\sum_{l\ge 1} \mu_l g_l &=&\sum_{l\le u_t} \mu_l g_l + \sum_{r\ge t}\leer{0.8}(\leer{-0.8}
               \sum_{u_r<l\le v_r}\leer{-0.5}\mu_l g_l +\leer{-1} \sum_{v_r<l\le u_{r+1}}
                                                   \leer{-1} \mu_l g_l)\\
         &=&  \sum_{l\le u_t} \mu_l g_l + \sum_{r\ge t} \sum_{i=1}^{a_{r+2}-b_{r+1}}
                   \mu_{u_r+i} G_i^{r+1}   \komma  }
denn $v_r-u_r=n_{r+2}-n_{r+1}= a_{r+2}-b_{r+1}$. Da $\sum_{l\le u_t} \mu_l g_l\in 
\Gamma(\Omega,\Rscript_{t+1})$, folgt nach Indexverschiebung 
\[   0=\sum_{r\ge t+1}\sum_{i=1}^{a_{r+1}-b_{r}}\tilde{\mu}_{r,i}[G_i^r]_k 
     \leer{2}\mbox{ mit }\tilde{\mu}_{r,i}=\mu_{u_{r-1}+i}\mbox{ für } k\ge b_{t+1}+1  \]
Also auch nach Umrechnung
\[   0=\sum_{r\ge t+1}\sum_{i=1}^{a_{r+1}-b_{r}} [A\tilde{\mu}]_{b_r+i}[C_i^r]_k 
           \leer{2}    \mbox{ für } k\ge b_t+1 \komma \]
wobei $\tilde{\mu}$ eine Folge ist mit $\tilde{\mu}_l=\tilde{\mu}_{r,i} $ für
$l=b_r+i$, $i=1,\ldots,a_{r+1} - b_{r}$ und $r\ge t+1$ und $\tilde{\mu}_l=0$ in allen
anderen Fällen.

Setze $\lambda=(\lambda_l)_l$ und $\lambda_l:= [A\tilde{\mu}]_{br+i}$ für $l=b_r+i$,
$i=1,\ldots,a_{r+1} - b_{r}$ und $r\ge t+1$ und sonst gleich Null. Da $[C_j^\tau]_{b_\tau +j}=-1$
für $\tau\ge t+1$, $j=1,\ldots,a_{\tau+1} - b_{\tau}$, folgt
\[   \ffb{A_p}{\lambda_{b_\tau+j }} = \ffb{A_p}{\sum_{r\ge t+1}
        \sum^{a_{r+1}-b_{r}}_{{i=1\atop i\not= j,\mbox{{\tiny falls }} r=\tau}} 
          \lambda_{b_r+i} [C_i^r]_{b_\tau+j }  }  \] 
für alle $p$, wobei $(A_p)_p$ die oben fest gewählte Ausschöpfung von $U$ ist.

Also gilt
\formelohne{
   \sum_{\tau\ge t+1} \sum_{j=1}^{a_{\tau+1}-b_\tau} \ffb{A_p}{\lambda_{b_\tau+j}}
       &=& \sum_{\tau\ge t+1} \sum_{j=1}^{a_{\tau+1}-b_\tau} \ffb{A_p}{
           \sum_{r\ge t+1} \sum^{a_{r+1}-b_{r}}_{{i=1\atop i\not= j,\mbox{{\tiny falls }} r=\tau}}
               \lambda_{b_r+i} [C_i^r]_{b_\tau+j }  }  \\
       &\le& \sum_{r\ge t+1} \sum_{i=1}^{a_{r+1}-b_r} \ffb{A_p}{\lambda_{b_r+i}}
             \sum_{\tau\ge t+1} 
             \sum^{a_{\tau+1}-b_{\tau}}_{{j=1\atop j\not= i,\mbox{{\tiny falls }} r=\tau}} 
             \ffb{A_p}{[C_i^r]_{b_\tau+j}}  \\
       &=& \sum_{r\ge t+1} \sum_{i=1}^{a_{r+1}-b_r} \ffb{A_p}{\lambda_{b_r+i}}
             \sum_{\tau= t+1}^{r-1} 
             \sum^{a_{\tau+1}-b_{\tau}}_{j=1} 
             \ffb{A_p}{c_{r,i,b_\tau+j}}\\
       &\le& \sum_{r\ge t+1} \sum_{i=1}^{a_{r+1}-b_{r}} \ffb{A_p}{\lambda_{b_r+i}} \cdot \theta }
für ein $\theta<1$ nach Eigenschaft \auslabeln{7.0.2} der Voraussetzung.

Hieraus folgt $\lambda_{b_r+i}=0$ für $r\ge t+1$, $i=1,\ldots,a_{r+1} - b_{r}$.

Da $\tilde{\mu}\in N_{t+1}$ folgt $\tilde{\mu}=A^{-1}\lambda =0$. Also ist $\tilde{\mu}_{r,i}
=0$ für $r\ge t+1$, $i=1,\ldots,a_{r+1} - b_{r}$. Dies bedeutet nach Rückverschiebung der 
Indizes  $\mu_{u_r+i}=0$
für $r\ge t$, $i=1,\ldots,v_{r} - u_{r}$. Insgesamt ist also $\mu_l=0$ für
alle $l\ge u_t +1$.

Hieraus folgt Eigenschaft \auslabeln{7.0.3} und die Proposition
\auslabelnsatz{7.3} ist bewiesen \Ende

Ist $U\relkom \Omega$ holomorph konvex. $(A_p)_p$ eine holomorph konvexe kompakte Ausschöpfung
von $U$. Dann sei 
\formelohne{
   B_p &:=& \{ f\in \zzc \doppelpunkt \norm{p}{f} \le 1 \} \\
   B_p^0 &:=& \{ f\in \zza\doppelpunkt \norm{p}{f}\le 1 \}   }
Wir kommen nun zum Hauptergebnis dieses Abschnittes.

\begin{satz}{Satz}\einlabelnsatz{8.1} Sei $\F$ ein kohärente analytische Garbe über $\Omega$.
$U\relkom\Omega$ holomorph konvex. Dann gibt es eine exakte Sequenz
\formel{  
   \ldots \zzc\stackrel{p_2^U}{\to} \zzc \stackrel{p_1^U}{\to} \zzc 
      \stackrel{p_0^U}{\to} \zza \to 0 \punkt   \einlabeln{8.0.3}}
Hierbei können die Abbildungen $p_i^U$, $i=0,1,\ldots$, so gewählt werden, daß
$p_i^U=\sum\limits_l \lambda_l f_l^i$ mit von $U$ unabhängigen $f_l^i\in l_1(H_\Omega(\oka))$
für alle $i\ge 1$ bzw. $f_l^i\in H_\Omega(\F)$ für $i=0$.

Es gilt ferner
\formel{  B_{p+1}\cap\abbkern p_{i-1}^U \subset p_i^U(C\cdot B_p) \leer{2}\mbox{für } i\ge 1 }
\einlabeln{8.0.4}
und einer von $p$ und $i$ abhängigen Konstante $C$, sowie 
\formel{  B_p^0\subset p_0^U(\theta\cdot B_p) \leer{2}\mbox{für jedes }\theta>1  
\einlabeln{8.0.5} \punkt }
\end{satz}

Sei $(K_t)_t$ eine holomorph konvexe kompakte Ausschöpfung von $\Omega$. Sei $0<\rho<1$ fest.
Ferner gelten die Bezeichnungen wie in der Proposition \auslabelnsatz{7.3}.

Wir zeigen die Exaktheit von \auslabeln{8.0.3} über die induktive Konstruktion 
der $(f_l^i)_l$:

Im Beweis werden wir $i'$ statt $i$ verwenden.

1) Induktionsanfang $i'=0$:

Ist $t\in \N_0$, dann gibt es nach Theorem A von Cartan $\ffc{\tilde{f}}{a_t} \in 
\Gamma(\Omega,\F)$, die $\F$ um $K_t$ erzeugen. Ohne Einschränkung ist $(a_t)_{t\ge 0}$ strikt
monoton steigend. (Ist z.B. $\F$ endlich erzeugt, dann werden die Erzeuger dupliziert.)

Setze ferner $f_1^0:=\tilde{f}_1,\ldots, f_{a_0}^0:=\tilde{f}_{a_0}$. Es seien für $t>0$ die Schnitte
$\ffc{f^0}{a_t}$ durch Skalierung mit Konstanten aus $\ffc{\tilde{f}}{a_t}$ hervorgegangen.
Betrachten wir $\tilde{f}_{a_t+i}$ in einer Umgebung von $K_t$ nahe genug an $K_t$, dann
lassen sich um $K_t$ holomorphe Funktionen $\lambda_{t,i,j}$, $j=1,\ldots,a_t$ finden, so daß 
\[      \tilde{f}_{a_t+i}= \ffa{j}{a_t} \lambda_{t,i,j} f_j^0 \leer{2}\mbox{gilt.}   \]
Wir setzen $C_{t,i}:=\ffa{j}{a_t} \ffb{K_t}{\lambda_{t,i,j}}$ und $f_{a_t+i}^0 := \tilde{f}_{a_t+i}\cdot
\frac{\rho}{C_{t,i}}$ für $i=1\ldots,a_{t+1}-a_t$.

Dann erzeugen $\ffc{f^0}{a_t}$ die Garbe $\F$ um $K_t$. Für jedes $U\relkom K_s$ holomorph konvex erhalten
wir eine Abbildung 
\[     p_0^U : \zzc \to H_U(\F)\komma (\lambda_l)_l \mapsto \sum_l \lambda_l f_l^0\komma  \]
die surjektiv und stetig linear ist. Denn es gilt für $U\relkom K_s $
und eine kompakte Ausschöpfung $(A_p)_{p\in \N}$ von $U$, sowie 
für $t\ge s$ 
\formelohne{ \ffb{A_p}{f_{a_t+i}^0} &=& \ffb{A_p}{\ffa{j}{a_t} \lambda_{t,i,j} f_j^0}\frac{\rho}{C_{t,i}} \\
                   &\le& \max_{j=1}^{a_t} \ffb{A_p}{f_j^0}  \ffb{A_p}{\lambda_{t,i,j}}\frac{\rho}{C_{t,i}} \\
                   &<& \max_{j=1}^{a_t} \ffb{A_p}{f_j^0}   }
Also induktiv 
\[   \ffb{A_p}{f_l^0} \le \max_{l=1}^{a_s} \ffb{A_p}{f_l^0} =: C_{p,s}  \leer{2} \mbox{für alle }
             l\ge a_s\punkt \]
Also gilt 
\formel{{ \einlabeln{8.0.6}}   \ffb{A_p}{\sum_l \lambda_l f_l^0} 
          &\le& \sum_l \ffb{A_p}{\lambda_l}\ffb{A_p}{f_l^0} \\
                &\le& C_{p,s} \sum_l \ffb{A_p}{\lambda_l} \\
                &=&  C_{p,s} \norm{p}{(\lambda_l)_l}  }
Wir bemerken, daß die $C_{p,s}$ eine gleichmäßige obere Schranke in $p$ besitzen.

$(f_l^0)_l$ erfüllt die Eigenschaften \auslabeln{7.0.2} und \auslabeln{7.0.3} aus der 
Proposition mit $b_t:=a_t$ alle $t$. Eigenschaft \auslabeln{7.0.2} ist wegen der obigen Skalierung
erfüllt und Eigenschaft \auslabeln{7.0.3} gilt automatisch wegen $a_t=b_t$  für alle $t$.

2) Die Folgen $(f_l^{i'})_l$ für $i'\ge 1$ erhalten wir induktiv durch Anwenden von 
Proposition \auslabelnsatz{7.3}. 

Bezogen auf die Folge $(f_l^{i'-1})_l$ ist damit 
\formelohne{ f_l^{i'} &=& \frac{\overline{g}_{l-n_{t+1}}}{S_{t-1}} \leer{0,5}
                         \mbox{ für } v_{t-1} +1\le l\le u_t \mbox{ und} \\
             f_l^{i'} &=& G_{l-n_t}^{t+1} \leer{0,5}
                         \mbox{ für } u_{t-1} \le l\le v_t \mbox{, wobei} }
die $S_t$ für $t\ge 1$ die Skalierungskonstanten aus \auslabeln{8.0.1.1} sind.

3) Es bleibt \auslabeln{8.0.4} und \auslabeln{8.0.5} zu zeigen. 
Sei $f\in H_U(\F)$ für ein $U\relkom \Omega$ holomorph konvex und sei $\norm{p}{f}\le 1$. 
Es sei $s$ so gewählt, daß $U\subset K_s $. Da $\ffc{f^0}{a_s}$ globale erzeugende Schnitte von
$\F$ über $U$ sind, kann Lemma \auslabelnsatz{7.2} angewandt werden.

Es gibt daher für $\theta>1$ Funktionen $\ffc{\lambda}{a_s}\in H_U(\oka)$ mit $\ffa{l}{a_s} \lambda_l f_l^0 
=f$ und $\ffa{l}{a_s} \ffb{A_p}{\lambda_l}\le \theta$. Damit ist $f\in p_0^U(\theta B_p)$.

Sei jetzt $i'\ge 1$, $f\in \abbkern p_{i'-1}^U$ und $\norm{p+1}{f}\le 1$ mit $f=(f_k)_k$.
Da $U\subset K_s$ können wir mit in $U$ holomorphen Funktionen
$\beta=(\ffc{\beta}{m_s})$ und $\alpha=(\alpha_{i,t})_{t\ge s\komma i=1,
\ldots,a_{t+1}-b_t}$ 
\[   f= \ffa{l}{m_s} \beta_l \overline{g}_l + \sum_{t\ge s}\ffa{i}{a_{t+1}-b_t} \alpha_{i,t}C_i^t
           \leer{2}\mbox{schreiben, wobei} \]
$(a_t)_t$, $(b_t)_t$, sowie $(C_i^t)_{t,i}$ etc. von der Folge $(f_l^{i'-1})_l$ herrühren gemäß
Proposition \auslabelnsatz{7.3}.

Da $\overline{g}_l\in \Gamma(\Omega,\Rscript_s)$ für $l=1,\ldots,m_s$ ist, 
gilt $[\overline{g}_l]_k=0$ für $k\ge b_s+1$.
Also gilt 
\[   f_k= \sum_{t\ge s}\ffa{i}{a_{t+1}-b_t} \alpha_{i,t}[C_i^t]_k \leer{2}\mbox{für } 
               k\ge b_s+1 \]

Wie im Beweis der Proposition \auslabelnsatz{7.3} gezeigt wurde, gilt
\formelohne{
    1\ge \sum_{k\ge b_s+1} \ffb{A_{p+1}}{f_k} &=& \sum_{k\ge b_s+1}\ffb{A_{p+1}}{
           \sum_{t\ge s}\ffa{i}{a_{t+1}-b_t} \alpha_{i,t}[C_i^t]_k}   \\
     &\ge& (1-\rho) \sum_{t\ge s}\ffa{i}{a_{t+1}-b_t} \ffb{A_{p+1}}{\alpha_{i,t}}\komma   }
wobei $\rho<1$ gemäß \auslabeln{7.0.2} gewählt ist.
Sei $\tilde{f}:= f- \sum_{t\ge s}\ffa{i}{a_{t+1}-b_t} \alpha_{i,t} C_i^t$. Dann ist 
\formelohne{  \norm{p+1}{\tilde{f}} &\le& \norm{p+1}{f} + \sum_{t\ge s}\ffa{i}{a_{t+1}-b_t}
               \ffb{A_{p+1}}{\alpha_{i,t}} \norm{p+1}{C_i^t}  \\
                   &\le& 1+\frac{1}{1-\rho}\cdot 2   }

Da $\ffc{\overline{g}}{m_s}$ die Garbe $\Rscript(f^{i'-1}_1,\ldots,f^{i'-1}_{b_s})$ 
um $K_s$ erzeugt, kann in Anwendung
von Lemma \auslabelnsatz{7.2} und Lemma \auslabelnsatz{7.0.1}
 $\ffc{\gamma}{m_s}\in H_U(\oka)$ gefunden werden, so daß 
\[   \tilde{f}=\ffa{l}{m_s} \gamma_l\overline{g}_l  \leer{2}\mbox{und}\leer{2}
           \ffa{l}{m_s}\ffb{A_p}{\gamma_l}  \le C_{i,p}\komma  \]
wobei $C_{i,p}$ von $f$ unabhängig ist.
Wir setzen $\alpha:=(\alpha_l)_l$ mit $\alpha_{b_t+i}:=\alpha_{i,t}$ für $t\ge s$, $i=1,\ldots, a_{t+1}-b_t$
und $a_l=0$ sonst.
Sei $\lambda:=A^{-1}\alpha$. Dann gilt 
$\norm{p+1}{\lambda}\le \frac{1}{1-\rho}\norm{p+1}{\alpha} \le \frac{1}{(1-\rho)^2}$.

Wir bilden $\mu=(\mu_l)_l$ durch 
\[   \mu_l:= \left\{ \begin{array}{c@{\leer{0.5} ; \leer{0.5}}l} 
                 \gamma_{m_t+i}\cdot S_{t-1} &\mbox{für } l=v_{t-1}+i\komma t\le s-1\komma i=1,
                                                 \ldots,u_{t}-v_{t-1} \\
                 \lambda_{i,t} &\mbox{für } l=u_{t}+i \komma t\ge s-1 \komma i=1,
                                                 \ldots,v_{t}-u_{t}   \\
                 0         &\mbox{sonst} \end{array} \right.   \]
Dann folgt
\formelohne{   p_{i'}^U\mu  & = & \sum_l \mu_l f_l^{i'} \\
                         & = & \sum_{t\le s-1} \sum_{i=1}^{u_t-v_{t-1}} \gamma_{m_t+i}\cdot S_{t-1}
                                f_{v_{t-1}+i}^{i'}
                              + \sum_{t\ge s-1} \sum_{i=1}^{v_t-u_{t}} \lambda_{i,t} f_{u_{t}+i}^{i'}\\
                         & = & \sum_{t\le s-1} \sum_{i=1}^{u_t-v_{t-1}} \gamma_{m_t+i}
                                \overline{g}_{m_t+1}
                              + \sum_{t\ge s-1} \sum_{i=1}^{v_t-u_{t}} \lambda_{i,t} G_i^{t+1}\\
                         & = & \ffa{l}{m_s} \gamma_l\overline{g}_l  +  
                               \sum_{t\ge s-1}\ffa{i}{a_{t+2}-b_{t+1}} \lambda_{i,t} G_i^{t+1} \\
                         & = & \tilde{f} + \sum_{t\ge s}\ffa{i}{a_{t+1}-b_t} \alpha_{i,t} C_i^t \\
                         & = & f  }
und $\norm{p}{\mu} = \norm{p}{\gamma} + \norm{p}{\lambda} \le C_{i,p} + \frac{1}{(1-\rho)^2}=:
\tilde{C}_{i,p}$. 
Also ist $f\in p_i^U(\tilde{C}_{i,p} H_p)$\Ende

\newpage
\section[{\sc Theorem B mit Schranken}]{
Theorem B mit Schranken für lokal endlich erzeugte Untergarben von 
$\oka_{\Omega}^p$}

Wir wollen nun Satz \auslabelnsatz{10.3} anwenden auf Koketten Werten in lokal endlich
erzeugten Untergarben ${\cal F}$ von $\oka^p$ für ein beliebiges natürliches $p$.
$\oka$ sei hier die Garbe von Keimen holomorpher Funktionen auf einer Steinschen
Mannigfaltigkeit $\Omega$.   
Nach dem Satz von Oka (siehe z.B. \cite{Hoer}, Theorem 6.4.1)
 ist ${\cal F}$ eine kohärente analytische Garbe.

Sei $(U_i)_{i\in I}$ eine  Überdeckung von $\Omega$ mit relativ kompakten,
Steinschen Gebieten in $\Omega$, wobei
es eine Schranke $M$ geben soll, so daß $U_i$ nicht von mehr als $M$ verschiedenen
Mengen aus $(U_i)_{i\in I}$ geschnitten wird.

Es sei $(A_{i,n})_{i\in I,n\in \N_0}$ ein System von holomorph konvexen kompakten Mengen in 
$\Omega$, so daß für jedes $i$ die Folge $(A_{i,n})_n$ eine Ausschöpfung von $U_i$ bildet.
Ohne Einschränkung sei $\Omega\subset \bigcup_{i\in I} \stackrel{\circ}{A}_{i,0}$ und 
$A_{i_1,0}\cap\ldots \cap A_{i_\sss,0}=\emptyset$ genau dann, wenn $U_{i_1}\cap\ldots
\cap U_{i_\sss}=\emptyset$. Für $\alpha=(i_1,\ldots,i_\sss)$ sei $A_{\alpha,n}:=
A_{i_1,n}\cap\ldots \cap A_{i_\sss,n}$. 

Auf den Fr\'echeträumen der holomorphen Funktionen auf $U_{i_0,\ldots,i_\sss}:= U_{i_0}\cap
\ldots\cap U_{i_\sss}$ ist jeweils eine Gradierung durch 
\[   \betrag{f}_{i_0,\ldots,i_\sss,n} := \ffb{A_{i_0,\ldots,i_\sss,n}}{f(z)}   \]
gegeben. 

Es sei für $\alpha\in I^{\sss+1}$, $\sss\ge 0$, jeweils $F^k_\alpha:=  l_1(H(U_\alpha))$, falls
$U_\alpha\not=\emptyset$ für 
$k=1,2,\ldots$. Hiebei trägt $F^k_\alpha$ die folgende Gradierung:
\[    \norm{\alpha,n}{(f_l)_{l\in \N}} := \sum_{l=1}^{\infty} \betrag{f_l}_{\alpha,n}\komma 
             n\in \N \punkt \]
Für $k=0$ und $\alpha\in I^{\sss+1}$, $\sss\ge 0$ sei jeweils $F_\alpha^0=\Gamma(U_\alpha,\F)$
versehen mit der Gradierung, die sich durch die Suprema auf $(A_{\alpha,n})_n$ ergibt.

Im Fall  $U_\alpha=\emptyset$ sei $F^k_\alpha:=\{0\}$ für $k=0,1,2,\ldots$.

Sei also $S=(S_{k,\sss})_{k,\sss}$ mit $S_{k,\sss}:=\{F_\alpha^k, \alpha\in I^{\sss+1}\}$.

Für die Indexmengen $I_\alpha$, $\alpha$ Multiindex, gilt dann $I_\alpha=\{i\in I\doppelpunkt U_\alpha\cap U_i
\not= \emptyset\}$. $I_\alpha$ ist also für jeden Multiindex $\alpha$ durch $M$ beschränkt
und es gilt $I_{\alpha,j}\subset I_\alpha$ für jeden Index $j\in I$ und jeden Multiindex 
$\alpha$.

F"ur $\alpha\in I^{\sss +1}$,  $\sss\ge 0$ mit $U_\alpha\not=\emptyset$ und $j\in \{0,\ldots,\sss\}$ 
sowie $k\ge 1$ sei 
\[   \iota_{\alpha_j}^{k,\alpha}: l_1(H(U_{\alpha_j})) \to l_1(H(U_{\alpha}))   \]

diejenige Abbildung, die komponentenweise $H(U_{\alpha_j})$ in $H(U_{\alpha})$ einbettet.
Diese Abbildung ist stetig, linear und injektiv und es gilt:
\[   \norm{\alpha,n}{\iota_{\alpha_j}^{k,\alpha}(f_l)_{l\in \N} } \le 
                   \norm{\alpha_j,n}{(f_l)_{l\in \N} }   \punkt\]
Ist $U_\alpha=\emptyset$, dann soll $\iota_{\alpha_j}^{k,\alpha}$ komponentenweise aus
der Nullabbildung bestehen.

Ist  $k=0$, dann  sei für $\alpha\in I^{\sss +1}$, $\sss\ge 0$ mit $U_\alpha\not=\emptyset$
und $j\in \{1,\ldots,\sss\}$
\[    \iota_{\alpha_j}^{0,\alpha}: H_{U_{\alpha_j}}(\F) \to H_{U_{\alpha}}(\F)  \]

die Einschränkungsabbildung von Schnitten auf $U_{\alpha_j}$ zu Schnitten auf $U_{\alpha}$.
Die Topologie auf $H_U(\F)$ für ein offenes $U\subset \Omega$ wurde in Kapitel \ref{Kap5}
eingeführt.
Diese Abbildung ist stetig linear mit 
\[   \norm{\alpha,n}{\iota_{\alpha_j}^{k,\alpha}(f_l)_{l\in \N} } \le 
                   \norm{\alpha_j,n}{(f_l)_{l\in \N} }   \]
für alle $n\in \N$. $\iota_{\alpha_j}^{k,\alpha}$ ist ebenfalls injektiv, da ja
$\F\subset \oka^p$ vorausgesetzt ist.

Ist $U_\alpha=\emptyset$, dann soll auch $\iota_{\alpha_j}^{0,\alpha}$ komponentenweise aus
der Nullabbildung bestehen.

Die Eigenschaften \auslabeln{10.1.1} sind erfüllt, da die $\iota_{\alpha_j}^{k,\alpha}$
jeweils Einschränkungsabbildungen sind.

Wir wollen zeigen, daß für $S$ die Eigenschaft (E1) gilt:

Wir setzen dazu $p_\alpha^k := p_k^{U_\alpha}$ für jedes $k\ge 0$ und $\alpha\in I^{\sss +1}$, $\sss\ge 0$,
 mit $ p^{U_\alpha}_k$ aus  Satz \auslabelnsatz{8.1}. Der Satz besagt gerade,
daß die Folge $(p^{U_\alpha}_k)_{k\in\N_0}$ mit absteigendem $k$ die in (E.1.1) 
geforderte exakte Sequenz bildet. 

Es gilt für $(\lambda_l)_l \in  l_1(H_{U_{\alpha_j}}(\oka))$ 
\[   p_\alpha^k\circ \iota_{\alpha_j}^{k+1,\alpha} ( (\lambda_l)_l) = 
                \sum_l \left[ \iota_{\alpha_j}^{k+1,\alpha} ( (\lambda_l)_l)\right]_l f_l^k\mid_{U_\alpha}
              =   \sum_l \lambda_l \mid_{U_\alpha} f_l^k\mid_{U_\alpha} 
              =   \iota_{\alpha_j}^{k+1,\alpha}\circ p_\alpha^k (( \lambda_l)_l)   \]
für $\alpha\in I^{\sss +1}$, $j=0,\ldots,\sss$, $k\ge 0$. Hierdurch wird (E.1.2) erfüllt. 

Die Eigenschaft (E.1.3) ist wegen \auslabeln{8.0.4} bzw. \auslabeln{8.0.5} erfüllt. 
(E.1.4) folgt aus 
  \auslabeln{8.0.1.2}.

Wir zeigen nun, daß auch Eigenschaft (E.2) für $S$ gilt.

Es sei für ein $\gamma>1$
\[ P_\gamma:= \{ \phi \in PSH(\Omega)\doppelpunkt \phi\ge 0\komma
        \sup_{U_i} \phi \le\gamma \inf_{U_i} \phi \mbox{ für jedes } i\in I\}\]
Wir setzen
\[   M_\gamma :=\{ [ (\sup_{U_\alpha} e^{-\phi})_{\alpha\in I^{\sss +1}} ,
                   (\sup_{U_\alpha} e^{-\phi})_{\alpha\in I^{\sss}}] \doppelpunkt
                     \phi\in P_\gamma \}   \]
Wir behaupten, daß zu jedem $\gamma>1$ die offensichtlich nicht leere
Menge $M_\gamma$ die Eigenschaften 1) bis 3) aus (E.2) besitzt:

{\sc dazu:} 

1) Sei $[C,C']\in M_\gamma$, dann gilt für alle $\alpha\in I^{\sss +1}$
beziehungsweise $\beta\in I^{\sss}$, daß $0<C_\alpha\le 1$ und $0<C'_\beta
\le 1$. Ist ferner $\alpha\in I^{\sss }$, dann gilt 
\[   C'_\alpha=\sup_{U_\alpha} e^{-\phi} \ge \sup_{U_\alpha\cap U_i} e^{-\phi} \]
für jedes $i\in I_\alpha$. Also gilt $C\prec C'$.

2) Ist $\theta>1$ und ist $[C,C']\in M_\gamma$, dann $[C^{\theta},{C'}^{\theta}]
\in M_\gamma$, denn $\theta\phi\in P_\gamma$.

3) Sei $[C, C']=[ (\sup_{U_\alpha} e^{-\phi})_{\alpha\in I^{\sss +1}} ,
                   (\sup_{U_\alpha} e^{-\phi})_{\alpha\in I^{\sss}}]\in M_\gamma$. 
Dann gilt für $C'':=(\sup_{U_\alpha} e^{-\phi})_{\alpha\in I^{\sss +2}}$, daß
$[C'',C]\in M_\gamma$ ist.

Es sei $C_1:=(C_{1,\alpha})_{\alpha\in I^{\sss +1}}$ und $n\in \N$. Es sei 
$[C,C']\in M_\gamma$ und $c\in {\cal C}_n(S_{k,\sss},C\cdot C_1)$ mit 
$\delta^\sss c=0$.

Für festes $\alpha\in I^{\sss +1}$ und $k\ge 1$ ist $c_\alpha=
(c_{\alpha,l})_{l\in \N}\in l_1(H(U_{\alpha}))$. Wir setzen $c^l:=
(c_{\alpha,l})_{\alpha\in I^{\sss +1}}$. Jedes $c^l$ ist eine Kokette der Länge
$\sss$ bezüglich der Überdeckung $(U_i)_{i\in I}$. $\delta^\sss c=0$ 
bedeutet gerade, daß $\delta^\sss c^l=0$ für jedes $l\in \N$.

Es gilt für jedes $n\in \N$ 
\formelohne{
        \infty > \dreinorm{2}{n,C\cdot C_1}{c}
             &=& \left(\sum_{\alpha\in I^{\sss +1}}
             C_\alpha C_{1,\alpha}   \sum_{l\in\N} \ffb{A_{\alpha,n}}{c_{\alpha,l}} \right)^2 \\ 
             &\ge& \sum_{\alpha\in I^{\sss +1}}
             C_\alpha^2 C_{1,\alpha}^2   \left( \sum_{l\in\N} \ffb{A_{\alpha,n}}
                                          {c_{\alpha,l}}\right)^2  \\
             &\ge& \sum_{\alpha\in I^{\sss +1}}
           \sum_{l\in\N}   C_\alpha^2 C_{1,\alpha}^2    \ffb{A_{\alpha,n}}
                                          {c_{\alpha,l}}^2  \punkt  }
Es sei  $\phi\in P_\gamma$ die zu $C$ gehörige Funktion. Wir wählen eine 
plurisubharmonische Funktion auf $\Omega$, so daß für jedes $\alpha\in I^{\sss +1}$ auf
$U_{\alpha}$ 
\[   e^{\psi_1}\ge C_{1,\alpha}^{-2}\cdot D_\alpha  \]
gilt, wobei $D_\alpha:=\int_{U_{\alpha}}  d\Omega(z)$.

Also gilt
\formelohne{
     \norm{2\phi+\psi_1,n}{c^l}^2 &:= &  \sum_{\alpha\in I^{\sss +1}}
                          \int_{A_{\alpha,n}^{\circ}} \betrag{c_{\alpha,l}}^2
                                       e^{-2\phi-\psi_1} d\Omega(z)  \\
                               &\le&   \sum_{\alpha\in I^{\sss +1}}
                          D_\alpha \ffb{A_{\alpha,n}}{c_{\alpha,l}}^2  
                                 \sup_{A_{\alpha,n}}    e^{-2\phi-\psi_1} \\
                               &\le&   \sum_{\alpha\in I^{\sss +1}}
                       C_\alpha^2 C_{1,\alpha}^2 \ffb{A_{\alpha,n}}{c_{\alpha,l}}^2<\infty}

Wir wenden nun das in Satz  \auslabelnsatz{3.6} gewonnene Theorem B mit 
Schranken für die Garbe $\oka=\H_0$ an auf  $n\in \N$ 
und die plurisubharmonische
Funktion $2\phi+ \psi_1$. Danach gibt es eine von $n$ unabhängige plurisubharmonische Funktion
$\psi_2$, so daß zu jeder holomorphen Kokette $c^l$ eine holomorphe Kokette ${c'}^l$ existiert mit 
$\delta^{\sss-1} {c'}^l = c^l$ und 
\[     \norm{2\phi +\psi_1 +\psi_2,n-1}{{c'}^l} \le \norm{2\phi+\psi_1,n}{c^l} \punkt\]

Bezüglich des Maßes $e^{-\psi_1}d\Omega$ und einer jeweils festen offenen Umgebung
$\tilde{U}_\alpha$ von $A_{\alpha,n-2}$ mit $\tilde{U}_\alpha\subset\subset 
A_{\alpha,n-1}^{\circ}$, $\alpha\in I^{\sss }$ gibt es Konstanten $K_{\alpha,n}$,
$\alpha\in I^{\sss }$, so daß 
\[   \ffb{\tilde{U}_\alpha}{f(z)}^2 \le K_{\alpha,n} \int_{A_{\alpha,n-1}^{\circ}}
                 \betrag{f(z)}^2 e^{-\psi_1} d\Omega(z)   \]
für jedes $f\in H(A_{\alpha,n-1}^{\circ})\cap L_2(A_{\alpha,n-1}^{\circ}, 
e^{-\psi_1}d\Omega)$.

Es gilt daher
\formelohne{
    \norm{2\phi+\psi_1,n}{c^l}^2 &\ge& \norm{2\phi+\psi_1+\psi_2,n-1}{{c'}^l}^2
                = \sum_{\alpha\in I^{\sss +1}} \int_{A_{\alpha,n-1}^{\circ}}
                  \betrag{c'_{\alpha,l}}^2 e^{-2\phi-\psi_1-\psi_2} d\Omega \\
                              &\ge& \sum_{\alpha\in I^{\sss}} K_{\alpha,n-1}^{-1}
                  \inf_{A_{\alpha,n-1}} e^{-\psi_2} \inf_{A_{\alpha,n-1}} e^{-2\phi}
                            \ffb{A_{\alpha,n-2}}{{c'}_{\alpha,l}}^2  }
                          
Wir setzen $C_{2,\alpha}:=K^{-\frac{1}{2}}_{\alpha,n-1} \inf_{U_\alpha} e^{-
\frac{1}{2}\psi_2} n(\alpha)^{-1}$ mit einer Funktion $n: I^{\sss} \to \R^{+}$,
so daß \[ \sum_{\alpha\in I^{\sss}} n(\alpha)^{-2} \le 1\punkt \]

Es gilt für $\alpha\in I^{\sss }$
\[   {C'}_\alpha^{2\gamma}  = \sup_{U_\alpha} e^{-2\phi\gamma} = e^{-2(\inf_{U\alpha} \phi)
      \gamma} \le e^{-2\sup_{U_\alpha}\phi} = \inf_{U_\alpha}e^{-2\phi} \le \inf_{
      A_{\alpha,n-1}}e^{-2\phi} \punkt \]
Insgesamt folgt: \formelohne{
     \dreinorm{2}{n-2,{C'}^\gamma C_2}{c'}
   &=& \left( \sum_{\alpha\in I^{\sss }} {C'}_\alpha^\gamma C_{2,\alpha} 
               \sum_{l\in \N} \ffb{A_{\alpha,n-2}}{{c'}_{\alpha,l}}\right)^2  \\
   &\le& \left[ \sum_{l\in \N} \left( \sum_{\alpha\in I^{\sss }} n(\alpha)^{-2}
         \right)^{\frac{1}{2}} \left( \sum_{\alpha\in I^{\sss }} K_{\alpha,n-1}^{-1}
         \inf_{U_\alpha} e^{-\psi_2} {C'}_\alpha^{2\gamma} 
            \ffb{A_{\alpha,n-2}}{{c}_{\alpha,l}}^{2}\right)^{\frac{1}{2}} \right]^2 \\
   &\le& \left[ \sum_{l\in \N} (
                 \norm{2\phi+\psi_1,n}{c^l}^2)^{\frac{1}{2}}\right]^2 \\
   &\le& \left[ \sum_{l\in \N} \left( \sum_{\alpha\in I^{\sss+1 }} 
         {C}_\alpha^2 C_{1,\alpha}^2 \ffb{A_{\alpha,n}}{{c}_{\alpha,l}}^2
                \right)^{\frac{1}{2}}  \right]^2  \\
   &\le& \left[ \sum_{l\in \N} \sum_{\alpha\in I^{\sss+1 }} 
         {C}_\alpha  C_{1,\alpha} \ffb{A_{\alpha,n}}{{c}_{\alpha,l}} \right]^2  \\
   &=&   \dreinorm{2}{n,{C} C_1}{c}             }
            
Damit ist nach Umnummerierung Eigenschaft (E.2) gezeigt.

Wir kommen nun zum Hauptergebnis dieses Abschnittes:

Sei $\Omega$ eine Steinsche Mannigfaltigkeit und $\F$ eine lokal
endlich erzeugte Untergarbe von $\oka^p$ auf $\Omega$ für
irgendein natürliches $p$. Sei ferner $(U_i)_{i\in I}$ eine abzählbare Steinsche
Überdeckung von $\Omega$, wobei höchstens $M$ paarweise verschiedene Überdeckungsmengen eine feste Überdeckungsmenge schneiden.

Für jedes $\alpha\in I^{\sss +1}$, $\sss\ge 0$ wird durch erzeugende Schnitte
$\ffc{F}{q}\in \Gamma(\Omega,\F)$ eine Fr\'echetraum\-topo\-logie auf $\Gamma(U_\alpha,\F)$ induziert. Die
einzelnen Stufen werden durch eine kompakte holomorph konvexe Aus\-schöpf\-ung
$(A_n)_n$ von $U_\alpha$ gegeben:
\[   \norm{\alpha,n}{f} := \inf_{f=\sum_{j=1}^q d_j F_j} \sum_{j=1}^q \ffb{
          A_n}{d_j}  \]
für $f\in \Gamma(U_\alpha,\F)$. Es gilt der folgende 

\begin{satzohnebeweis}{Satz}\einlabelnsatz{11.1} 
Sei für jedes $\alpha\in I^{\sss +1}$, $\sss\ge 0$, ein zu $(\normleer_{\alpha,n})_n$
äquivalentes Halbnormensystem $(\betragleer_{\alpha,n})_n$ gegeben auf
 $\Gamma(U_\alpha,\F)$ und sei $m\in \N$ fest gewählt. 

Dann gibt es zu jedem $\sss\ge 0$ ein Konstantensystem 
$(D_\alpha)_{\alpha\in I^{\sss}}$, $0< D_\alpha \le 1$, so daß für jede auf 
$\Omega$ plurisubharmonische Funktion $\phi\ge 0$, für die es ein $\gamma>1$ gibt mit 
\[    \sup_{U_i} \phi \le \gamma\inf_{U_i} \phi \mbox{ für alle } i\in I\komma\]
               
das folgende gilt:

Zu jeder Kokette $c=(c_\alpha)_{\alpha\in I^{\sss +1}}$, $\sss\ge 0$, mit 
$c_\alpha\in \Gamma(U_\alpha,\F)$, $\delta^\sss c=0$ und
\formel{ 
     \sum_{\alpha\in I^{\sss +1}} \sup_{U_\alpha} e^{-\phi} \sup_n
         \betrag{c_\alpha}_{n,\alpha} < \infty   \einlabeln{11.1.1} }
gibt es eine Kokette $c':=({c'}_\alpha)_{\alpha\in I^{\sss }}$
mit ${c'}_\alpha \in \Gamma(U_\alpha,\F)$, $\delta^{\sss-1} c'=c$ und 
\[   \sum_{\alpha\in I^{\sss }} D_\alpha \sup_{U_\alpha}{e^{-\gamma^M \phi}}
           \betrag{{c'}_\alpha}_{\alpha,m}  \le  \sum_{\alpha\in I^{\sss +1}}
           \sup_{U_\alpha}{e^{-\phi}}  \sup_n  \betrag{c_\alpha}_{n,\alpha} \]
\end{satzohnebeweis}
{\sc Bemerkung:} Wegen Lemma \auslabelnsatz{7.0.1} sind die Supremumsnormen
in $H(U_\alpha)$ auf einer jeweiligen kompakten Ausschöpfung von $U_\alpha$
ein äquivalentes Halbnormensystem zu $(\normleer_{\alpha,n})_n$.

{\sc Beweis des Satzes:} Für das in Kapitel \ref{Kap2} definierte System $S$ von Fr\'echeträumen 
ist oben die Eigenschaft
(E) nachgewiesen worden. Wir wollen Satz \auslabelnsatz{10.3} für $S$ anwenden bezüglich der
Überdeckung $(U_i)_i$. Wir übernehmen die Bezeichnungen aus Kapitel \ref{Kap2}

Wir können   ohne Einschränkung der 
Allgemeinheit annehmen, daß es zu einem festen $m$ ein
 $n\in \N$ und Konstanten $(K_{1,\alpha})_{\alpha\in I^{\sss +1}}$ gibt, so daß
\[    K_{1,\alpha} \betrag{f}_{\alpha,m} \le \norm{\alpha,n-2(M+1)}{f}  \]
für jedes $f\in \Gamma(U_\alpha,\F)$ und jedes $\alpha\in I^{\sss+1}$.

Ferner gibt es $K_{2,\alpha}$, so daß
\[   K_{2,\alpha} \norm{\alpha,n}{f} \le \sup_k \betrag{f}_{\alpha,k} \]
für jedes $f\in \Gamma(U_\alpha,\F)$, für das die rechte Seite existiert.

Da $\phi\in P_\gamma$ ist $[C:=(\sup_{U_\alpha} e^{-\phi})_{\alpha\in I^{\sss +1}},
C':=(\sup_{U_\alpha} e^{-\phi})_{\alpha\in I^{\sss }} ] \in M_\gamma$.
Es sei $C_1:=(K_{2,\alpha})_{\alpha\in I^{\sss +1}}$.

Sei also $c=(c_\alpha)_{\alpha\in I^{\sss +1}}$ mit $\delta^\sss c=0$ und 
\auslabeln{11.1.1}, dann folgt $c\in {\cal C}_n(S_{0,\sss},C\cdot C_1)$,
denn
\[   \dreinorm{}{n,C\cdot C_1}{c} = \sum_{\alpha\in I^{\sss +1}} C_\alpha
           C_{1,\alpha} \norm{\alpha,n}{c_\alpha}  \le  \sum_{\alpha\in I^{\sss +1}}
            \sup_{U_\alpha} e^{-\phi} \sup_k\betrag{c_\alpha}_{\alpha,k} <\infty \]

Wir wenden Satz \auslabelnsatz{10.3} an und erhalten 
$c'\in {\cal C}_{n-2(M+1-\sss)}(S_{0,\sss-1},{C'}^{\gamma^{M-\sss}}\cdot C_2)$
mit $\delta^{\sss-1}c'=c$ und 
\[   \dreinorm{}{n-2(M+1-\sss),{C'}^{\gamma^{M-\sss}}\cdot C_2}{c'} \le
     \dreinorm{}{n,C\cdot C_1}{c} \punkt  \]
Hierbei ist das Konstantensystem $C_2$ unabhängig von $C$, $C'$ und $c$.
Wir setzen $D_\alpha:= C_{2,\alpha}\cdot K_{1,\alpha}$ für $\alpha\in I^{\sss }$ und
erhalten 
\formelohne{  
     \sum_{\alpha\in I^{\sss }} D_\alpha \sup_{U_\alpha} e^{-\gamma^M\phi} 
     \betrag{c_\alpha}_{\alpha,m}
  &\le& \sum_{\alpha\in I^{\sss }} C_{2,\alpha} {C'}^{\gamma^{M-\sss}} 
        \norm{\alpha,n-2(M+1-\sss)}{{c'}_\alpha}  \\
  &\le& \sum_{\alpha\in I^{\sss +1}} C_{1,\alpha} C_\alpha \norm{\alpha,n}{c_\alpha} \\
  &\le& \sum_{\alpha\in I^{\sss +1}} \sup_{U_\alpha} e^{-\phi} \sup_k 
          \betrag{c_\alpha}_{\alpha,k}  \punkt    }
Ende des Beweises $\Box$

\newpage
\section[{\sc Konstruktion einer Überdeckung des $\R^N$}]{Konstruktion 
einer Überdeckung des $\R^N$}

Wir zeigen in diesem Kapitel, daß es spezielle Überdeckungen des
$\R^N$ mit achsenparallelen Kuben gibt:

\begin{satzohnebeweis}{Proposition}\einlabelnsatz{14.1}
Gegeben sei eine stetige Funktion $\phi:\R^N \to \R^{+}$, dann gibt es 
eine abzählbare Überdeckung $(K_i)_{i\in I}$ mit achsenparallelen offenen Kuben 
der Seitenlänge $s_i$, so daß folgendes erfüllt ist:
\be
\item[1)] $\phi(x)> s_i$ für jedes $x \in K_i$ und $i\in I$.
\item[2)] Für jedes $i\in I$ ist die Anzahl der Elemente der Menge
   $\{  j\in I\doppelpunkt K_i\cap K_j\not=\emptyset \}$ durch die 
   Zahl $4^N - 2^N$ beschränkt.
\ee
\end{satzohnebeweis}

Sei $T\subset \{ f:\R^{+}_0 \to \R^{+}_0 $ beschränkt $\}$ durch folgendes
definiert:

Ist $f\in T$, dann gibt es eine Folge $0=t_0 < t_1 < t_2 <\ldots$ mit $t_k \nearrow \infty$,
so daß
\be
\item[1)] entweder
\ee
   \formel{ t_{k+1}-t_k = t_k - t_{k-1} \mbox{ oder } t_{k+1}-t_k = 
           \frac{1}{2}(t_k - t_{k-1})  \mbox{ gilt }\foralls k\in \N\komma 
          \einlabeln{14.1.0}}
\be
\item[2)] $f(x)=t_{k+1}-t_k$ für jedes $x\in [t_k,t_{k+1}[ $ und
          $k=0,1,\ldots$.
\ee

\begin{satz}{Lemma}\einlabelnsatz{14.2} Sei $(f_k)_k$ eine Folge
wie oben, dann gibt es
zu jedem $k\in\N$  ein $m\in \N$ und $n\in \N_0$, so daß
$t_k= m(t_k-t_{k-1})$ und $t_k-t_{k-1}= 2^{-n} (t_1-t_0)$

\end{satz}
Wir führen eine Induktion über $k$ durch. Der Induktionsanfang $k=1$ ist klar.
Sei die Behauptung für $k>1$ richtig. Dann gilt $m\,t_{k-1}=(m-1)t_k$. Ferner
gibt es ein $n\in \N$, so daß $t_k-t_{k-1}= 2^{-n} (t_1-t_0)$.

Ist im 1. Fall von \auslabeln{14.1.0} $t_{k+1}-t_k = t_k - t_{k-1}$, dann 
folgt zum einen $t_{k+1}-t_k=2^{-n} (t_1-t_0)$, zum anderen folgt
\[  mt_{k+1} = 2mt_k - mt_{k-1} = (2m - m+1)t_k  \]
also gilt
\[  (m+1)t_{k+1} = (m+1)t_k + t_{k+1} \punkt  \]

Im 2. Fall ist  $t_{k+1}-t_k = \frac{1}{2}(t_k - t_{k-1})$. Dann folgt zum
einen $t_{k+1}-t_k=2^{-(n+1)} (t_1-t_0)$, zum anderen folgt
\[ 2m t_{k+1} = 3mt_k - mt_{k-1} = (3m-m+1)t_k \]
also gilt
\[ (2m+1)t_{k+1} = (2m+1)t_k + t_{k+1}  \mbox{ \Ende} \]

\begin{satz}{Lemma}
Falls es ein $n_0\in \N$ gibt, so daß $t_1-t_0 = 2^{-n_0}$, dann gibt es zu 
jedem $l\in \N$ ein $k\in \N$, so daß $t_k=l$.
\end{satz}
Da $t_k\nearrow \infty$, gibt es ein $k\in \N$, so daß $t_{k-1}<l\le t_k$.
Nach Lemma \auslabelnsatz{14.2} gibt es $n,m\in \N$, so daß
$t_k=2^{-n}m$, sowie $t_k-t_{k-1} = 2^{-n}$. Hieraus folgt
$t_{k-1}=(m-1)2^{-n}$. Aus
\[   (m-1)2^{-n} < l\cdot 2^n 2^{-n} \le m \cdot2^{-n}   \]
folgt
\[ m-1<l \cdot 2^n \le m  \]
Also gilt $m=l\cdot 2^n$ und damit $t_k=2^{-n}l\, 2^n  = l $ \Ende

\begin{satz}{Lemma}\einlabelnsatz{14.4}
Sei $\phi :\R_0^{+} \to \R^{+}$ eine stetige Funktion, dann gibt es ein 
$f\in T$, mit $f<\phi $.
\end{satz}
Sei $\tilde{a}_k := \inf \{ \phi(x)\doppelpunkt x\in [k,k+1] \}$ und 
$a_k := \min \{ \tilde{a}_l \doppelpunkt l\le k\} \cup \{1\}$ jeweils für 
$k=0,1,2,\ldots$. Dann ist $(a_k)_k$ eine monoton fallende Folge positiver
Zahlen mit $a_0 \le 1$.

Wir konstruieren induktiv eine Folge $(t_k)_k$ mit \auslabeln{14.1.0}
und $t_k - t_{k-1} <a_l$, falls $t_k\in ]l,l+1]$.

Sei die Folge $(t_k)_k$ bereits konstruiert bis zu einem $t_K = l$ und 
sei $t_K - t_{K-1}<a_l$. 
Sei für $j\ge 1$ mit $t_{K+j} \le l+1$ 
\formel{
   t_{K+j} := \left\{ \begin{array}{l@{\mbox{ falls }}l}
   t_{K+j-1} + \frac{1}{2} (t_{K+j-1}-t_{K+j-2})\komma & t_{K+j-1}-t_{K+j-2} \ge a_{l+1} \\
   t_{K+j-1} + (t_{K+j-1}-t_{K+j-2})\komma & t_{K+j-1}-t_{K+j-2} < a_{l+1}
   \end{array} \right. \punkt  \einlabeln{14.4.1}}

Dann ist \auslabeln{14.1.0} erfüllt. Es ist klar, daß dann
\[  t_k - t_{k-1} \le t_K - t_{K-1} < a_l \foralls k\ge K \mbox{ mit }
t_k \le l+1 \mbox{ gilt.} \]

Es gibt ein $J$, so daß $t_{K+J}=l+1$. Wir zeigen, daß 
$t_{K+J}-t_{K+J-1} < a_{l+1}$ gilt:

Angenommen, es sei $t_{K+J}-t_{K+J-1} \ge a_{l+1}$, dann gilt für 
$t_{K+J-i}-t_{K+J-i-1} \ge a_{l+1}$ für jedes $1\le i\le J$.
Aus der Konstruktion in \auslabeln{14.4.1} folgt daher für
$1\le i\le J$
\[   t_{K+J-i}-t_{K+J-i-1} = 2^{-J+i} ( t_K-t_{K-1})  \]
Andererseits folgt wegen
\formelohne{
   1=t_{K+J} - t_K &=& \sum_{i=0}^{J-1} (t_{K+J-i}-t_{K+J-i-1}) \\
                   &=& (t_K-t_{K-1})\sum_{i=0}^{J-1} 2^{-J+i} \\
                   &<&  t_K-t_{K-1} <a_l \le 1  }
ein Widerspruch.

Der Induktionsschritt der Konstruktion ist damit gezeigt. Für das
Intervall $[0,1]$ setzen wir $t_0=0$ und $t_1:=2^{-n_0}$, wobei
$2^{-n_0}<a_0$ gelten soll. Wir benutzen \auslabeln{14.4.1}, um
die weiteren $t_k$ bis zur Stelle $1$ zu konstruieren.

Sei $f(x):=t_k-t_{k-1}$ für $x\in [t_{k-1},t_k[$. $f$ ist in $T$
und $f(x)<\phi(x)$ für jedes $x\in \R_0^{+}$ \Ende

{\sc Beweis der Proposition:}
Wir setzen 
\formelohne{
    A_t     &:=& \{ x\in \R^N \doppelpunkt \betrag{x_j}\le t \foralls 1\le j\le N \}\\
    R_t     &:=& \partial A_t  \\
    \psi(t) &:=& \inf_{x\in R_t} \phi(x)  \punkt }
Da $\psi$ stetig ist, können wir Lemma \auslabelnsatz{14.4} anwenden.
Wir erhalten eine Folge $0=t_0<t_1<t_2<\ldots$ mit den genannten Eigenschaften.
Nach Lemma \auslabelnsatz{14.2} gibt es zu jedem $k\ge 1$ ein $m\in \N$,
so daß $t_k=m\cdot(t_k-t_{k-1})$. Daher läßt sich $A_{t_k}$ durch
exakt $(2m)^N$ achsenparallele abgeschlossene Kuben der Seitenlänge
$(t_k-t_{k-1})$ überdecken. Die Menge dieser Kuben sei mit $U_k$
bezeichnet.

Für $K\in U_k$ gilt
\[   t_k-t_{k-1} < \phi(x) \foralls x \in K  \komma \]
denn wenn $t_k$ in $]l,l+1]$ liegt, dann folgt 
\[   t_k-t_{k-1} < a_l \le \inf_{t\in [0,l+1]}\psi(t) = 
     \inf_{x\in A_{l+1}}\phi(x) \le\inf_{x\in K}\phi(x) \]

Sei $W_k:=\{K\in U_k\doppelpunkt \stackrel{\circ}{K}\cap A_{t_{k-1}} =
\emptyset \}$. Dann sei $(L_i)_i$ eine Abzählung von $\bigcup_{k} W_k$.
Es ist leicht zu sehen, daß $\bigcup_i L_i =\R^N$.

Ist $l_i$ die Seitenlänge von $L_i$, dann gilt ferner, daß $L_i\cap
L_j\not=\emptyset$ entweder $l_i=l_j$ oder $l_i=2l_j$ oder $l_j=2l_i$ impliziert,
denn ist $L_i\in W_k$, dann ist $L_j\in W_{k-1}\cup W_k\cup W_{k+1}$, wenn
$L_i\cap L_j\not=\emptyset$ ist.

Hieraus folgt wiederum, daß die Anzahl der Elemente in $J_i:=
\{ j\in I\doppelpunkt L_i\cap L_j\not=\emptyset \}$ durch
$4^N-2^N$ beschränkt ist.

Die gesuchte Überdeckung $(K_i)_{i\in I}$ erhalten wir durch
jeweilige Vergrößerung der  $L_i$ zu einem Kubus $K_i$ der
Seitenlänge $s_i$ mit 
\be 
\item[1)] $\stackrel{\circ}{K}_i \supset L_i$,
\item[2)] $l_i<s_i< \inf_{x\in K_i}\phi(x)$ und 
\item[3)] $K_i\cap K_j\not=\emptyset$ nur dann, wenn $L_i\cap L_j\not=\emptyset$,
\ee
für alle $i,j\in I$ \Ende

\newpage
\section[{\sc Ein Fortsetzungssatz für holomorphe Funktionen}]{Ein 
Fortsetzungssatz für holomorphe Funktionen auf Steinschen
Mannigfaltigkeiten}

Als Anwendung von Satz \auslabelnsatz{11.1} zeigen wir den folgenden 
Fortsetzungssatz für holomorphe Funktionen 
auf Steinschen Untermannigfaltigkeiten des $\C^N$.

Sei $V$ eine abgeschlossene und damit Steinsche Untermannigfaltigkeit des $\C^N$. 
 Dann gilt der folgende 

\begin{satzohnebeweis}{Satz}\einlabelnsatz{12.1}
Zu jedem Kompaktum $K\subset \C^N$ und jedem $\gamma>1$ gibt es eine Konstante 
$B_{K,\gamma}$, so daß  für alle $f\in H(V)$,  und jede stetige Funktion $C_f:
\C^N\to \R_0^+$ das folgende äquivalent ist:
\be
\item[i)] Es gibt eine plurisubharmonische Funktion $\phi>0$ auf dem $\C^N$, so daß
\be
   \item[a)] $\log\betrag{f(z)} \le \phi(z)$ für jedes $z\in V$;
   \item[b)] $\sup_K\phi(z)\le \sup_{K\cap V} \log\betrag{f(z)} + \sup_K C_f(z)$
             für jedes $K\subset \C^N$ kompakt mit $K\cap V\not= \emptyset$ und 
             $K\supset\not= K_f:=\{z\in V\doppelpunkt \betrag{f(z)}\le 1\} $.
\ee
\item[ii)] Für jedes $\gamma>1$ gilt: Zu jedem $m\in \N$ und jedem Kompaktum 
            $K_0\subset \C^N$ mit $K_0\cap V\not=\emptyset$ und 
           $K_0\supset\not= K_f$ gibt es eine ganze Funktion $F_{m,K_0}$, so daß 
\be
   \item[a)]$ F_{m,K_0}(z)=f^m(z)$ für jedes $z\in V$;
   \item[b)] $ \sup_K \betrag{F_{m,K_0}(z)} \le e^{B_{K,\gamma} + m\gamma\sup_K C_f } 
             \sup_{K\cap V} \betrag{f^m(z)}^\gamma $ für jede kompakte Menge
             $K_0\supseteq K\supset\not= K_f$ mit $K\cap V\not= \emptyset$.
\ee
\item[iii)] Für jedes $\gamma>1$ gilt: Zu jedem $m\in \N$ gibt es eine ganze Funktion
$F_m$, so daß
\be
   \item[a)]$ F_{m}(z)=f^m(z)$ für jedes $z\in V$;
   \item[b)] $ \sup_K \betrag{F_{m}(z)} \le e^{B_{K,\gamma} + m\gamma\sup_K C_f } 
             \sup_{K\cap V} \betrag{f^m(z)}^\gamma $ für jede kompakte Menge
             $ K\supset\not= K_f$ mit $K\cap V\not= \emptyset$.
\ee
\ee
\end{satzohnebeweis}

Bevor wir mit dem Beweis des Satzes beginnen, benötigen wir zwei Lemmata:

\begin{satz}{Lemma}\einlabelnsatz{12.2}
Sei $X$ ein lokalkompakter, $\sss$-kompakter topologischer Raum. Sei 
$(K_n)_{n\in \N_0}$ eine kompakte Ausschöpfung von $X$ mit
$\stackrel{\circ}{K}_{n+1}\supset K_n$ für jedes $n\in \N_0$
 und $(\delta_n)_{n\in \N}$
eine monoton fallende Folge positiver Zahlen. Dann gibt es eine stetige positive
Funktion $\psi$ auf $X$, so daß für jedes $x\in X$ ein $n\in \N$ existiert mit
$\psi(x)\le \delta_n$ und $x\in K_n$. 
\end{satz}

Nach dem Lemma von Urysohn finden wir Funktionen $\phi_n:X\to[0,1]$,  $n\in \N$, so
daß $\phi_n(x)=1$ für $x\in K_{n-1}$ und $\phi_n(x)=0$ für $x\in X\backslash K_n$ (siehe
\cite{M-V}, Lemma 4.18). Sei $\psi_n:=\max(\delta_1\phi_1, \ldots,$ $ \delta_n\phi_n)$.
Dann gilt $\psi_n(x)\le \delta_j$ für jedes $x\in X\backslash K_{j-1}$, $j=1,\ldots,n$, und 
$\psi_n(x)\ge \delta_n$ für jedes $x\in K_{n-1}$. Es gilt für jedes $m\ge n$ und
jedes $x\in K_{n-1}$ daher $\psi_n(x)=\psi_m(x)$. Es wird also durch 
$\psi(x):=\lim_n\psi_n(x)$ auf $X$ eine stetige positive Funktion definiert mit
$\psi(x)=\psi_{n+1}(x)\le \delta_n$ für jedes $x\in K_n\backslash K_{n-1}$ \Ende







\begin{satz}{Lemma}\einlabelnsatz{12.3}
Sei $\gamma>1$, $\psi>0$ eine stetige reellwertige Funktion auf einer offenen Menge $\Omega
\subset \C^N$. Dann gibt es eine stetige reellwertige Funktion $\xi>0$ auf $\Omega$,
so daß für jedes relativ kompakte, offene $U\subset \Omega$ 
mit $\sup \{\betrag{x-y}\doppelpunkt
x,y\in U\} <\inf_U\xi(z)$ folgt
\[   \sup_U \psi \le \gamma \inf_U\psi  \]

\end{satz}

Zu jedem $K\subset\Omega$ kompakt wählen wir $r_K>0$, so daß $\betrag{\log\psi(x)-
\log\psi(y)} <\log \gamma$ für $x,y\in K$ mit $\betrag{x-y}<r_K$. 
Hieraus folgt $\frac{\psi(x)}{\psi(y)} < \gamma$ für $\betrag{x-y}<r_K$.

Sei $(K_n)_{n\in \N_0}$ eine kompakte Ausschöpfung von $\Omega$. Ohne
Einschränkung ist $r_n:=r_{K_n}$, $n\in \N$, eine monoton fallende Folge. Dann gibt
es nach Lemma \auslabelnsatz{12.2} eine stetige positive Funktion $\xi$ auf 
$\Omega$, so daß zu jedem $x\in \Omega$ ein $n\in \N$ existiert mit 
$\xi(x)\le r_n$ und $x\in K_n$. 

Sei $U\subset \Omega$ relativ kompakt mit $\diam(U):=\sup \{\betrag{x-y}\doppelpunkt
x,y\in U\} <\inf_U\xi(z)$. Sei $n\in \N$ so gewählt, daß $U\subset K_n$ und es gibt
ein $x\in U$ mit $x\not\in K_{n-1}$, dann folgt $\xi(x)\le r_n$. Also gilt
für alle $z,y\in U$, daß $\betrag{z-y}<\xi(x)\le r_{K_n}$. 
Daher folgt $\sup_U\psi \le \gamma \inf_u\psi$ \Ende



{\sc \bf Beweis von Satz \auslabelnsatz{12.1}:} Wir zeigen zunächst die Äquivalenz
von $ii)$ und $iii)$, wobei $iii)\Rightarrow ii)$ klar ist.

$ii)\Rightarrow iii)$: Sei $(K_n)_{n\in \N}$ eine kompakte Ausschöpfung des $\C^N$ mit
$K_1 \supset\not= K_f$ und $K_1\cap V\not= \emptyset$, dann ist
$(\betrag{F_{m,K_n}})_n$ eine  gleichmäßig auf Kompakta beschränkte Folge ganzer
Funktionen. Nach dem Satz von Stieltjes-Vitali in mehreren komplexen Veränderlichen
(siehe z.B. \cite{Hoer}, Corollary 2.2.5) gibt es eine ganze Funktion $F_m$, deren
Betrag ein Häufungspunkt der obigen Folge ist. Damit erfüllt $F_m$ die Aussagen
iii), a) und iii), b).

$i)\Rightarrow ii)$: Wir nehmen zunächst an, daß $\phi$ stetig ist und wenden Lemma 
\auslabelnsatz{12.3} auf $\phi$ an.

Es gibt eine offene Umgebung $W$ von $V$ im $\C^N$ und eine holomorphe 
Retraktionsabbildung $\pi:W\to V$ mit $\pi(z)=z$ für jedes $z\in V$ 
(siehe z.B. \cite{G-R}, Chap. VIII, Sec. C, Theorem 8).
 Allgemeiner als hier benötigt zeigt Siu in \cite{Siu}, Corollary 1, 
daß für jede komplexe Untermannigfaltigkeit $V$ einer Steinschen Mannigfaltigkeit
$\Omega$ ein holomorpher Retrakt von einer in $\Omega$ offenen Umgebung von $V$  auf 
$V$ existiert.

Sei $\gamma >1$. Wir wenden Lemma \auslabelnsatz{12.3} auf $\phi$, $\gamma$ und 
$W$ an und erhalten eine stetige positive Funktion $\xi$ auf $W$ mit den entsprechenden
Eigenschaften.

Sei für ein $U\subset W$ offen und relativ kompakt in $\C^N$ und $x\in U$ jeweils
$B_{x,U}$ die kleinste Kugel um $\pi(x)$, die $U$ enthält.

Es sei $K\subset V $ kompakt, $x_0\in K$ beliebig.
Da $\xi $ stetig und positiv ist auf $W$, gibt es daher ein $\eps_K$, so daß
für jedes $y\in K\cup(B_{2\eps_K}(x_0)\cap V)$ gilt
$B_{2\eps_K}(y) \subset W$ und 
\[   \inf\{ \xi(x)\doppelpunkt x\in B_{2\eps_K}(y)\} \ge \frac{1}{2}\xi(y) \ge
     4\eps_K  \punkt \]

Andererseits gibt es ein $\delta_K$, so daß für jedes $x_0\in K$ gilt
\[   \betrag{\pi(x)-x}< \eps_K \foralls x\in B_{\delta_K}(x_0) \und \]
\[  \delta_K \le \eps_K  \mbox{ und  damit }\delta_K \le \xi(x_0)  \punkt \]

Sei $x_0\in K$ und $\eta $ eine stetige positive Funktion auf dem $\C^N$ mit
$\eta(x_0) \le \delta_K$. Dann gilt für jede offene zusammenhängende Umgebung
$U$ von $x_0$ in $W$ mit $\diam(U)<\inf_U \eta$, daß 
$U\subset B_{\delta_K}(x_0) $. Damit ist $\betrag{\pi(x)-x}<\eps_K$ für
jedes $x\in U$. Sei $y\in U$, dann ist 
\[  \diam(B_{y,U}) \le 2(\betrag{\pi(x)-x} + \diam(U)) \le 2(\eps_K + \delta_K)
    \le 4 \eps_K  \punkt \]

Da $\betrag{\pi(y)-x_0}<2\eps_K$, gilt $4\eps_K\le \frac{\xi(\pi(y))}{2}$. Also
folgt 
\[  \diam(B_{y,U})\le \inf\{ \xi(x)\doppelpunkt x\in B_{y,U} \}  \punkt \]

Sei nun $(K_n)_{n\in \N_0}$ eine kompakte Ausschöpfung von $V$ und $\delta_n:=
\delta_{K_n}$, dann gibt es nach Lemma \auslabelnsatz{12.2} eine 
stetige positive Funktion $\eta$ auf $V$, so daß zu jedem $x_0\in V$ ein
$n\in \N$ existiert mit $\eta(x_0)\le \delta_n$ und $x_0\in K_n$. Sei
$\tilde{\eta}$ eine positive stetige Funktion auf dem $\C^N$ mit
$\tilde{\eta}(x)=\eta(x)$ für jedes $x\in V$ nach dem Theorem von Tietze-Urysohn
(siehe z.B. \cite{Eng}, Theorem 2.1.8). Dann hat $\tilde{\eta}$ die Eigenschaft,  
 daß für eine offene zusammenhängende Menge $U\subset\C^N$, $U\cap V\not=\emptyset$
mit 
\formel{     \sup\{\betrag{x-y}\doppelpunkt x,y\in U\}< \inf_U\tilde{\eta}(z) 
\einlabeln{12.1.1}}
 schon
folgt $U\subset W$ und $ \sup\{ \betrag{z-y}\doppelpunkt z,y\in B_{x,U}\}
<\inf_{B_{x,U}}\xi(z)$ für
jedes $x\in U$.

Wir wenden Satz \auslabelnsatz{14.1} auf $\psi:=\frac{1}{\sqrt{2N}} \tilde{\eta} $ an und 
erhalten eine Überdeckung des $\C^N=\R^{2N}$ mit achsenparallen Kuben $(U_i)_{i\in I}$.

Ist $M:= 16^N-4^N$, dann haben höchstens M paarweise verschiedene Überdeckungsmengen
aus $(U_i)_{i\in I}$ einen nichtleeren Schnitt mit einer festen Überdeckungsmenge aus
$(U_i)_{i\in I}$.

Es gilt $\sup\{ \betrag{x-y}\doppelpunkt x,y\in U_i\}<\inf_{U_i}\tilde{\eta}$ für jedes $i\in I$.
Hieraus folgt:
\be
\item[1)] Ist $U_i\cap V\not=\emptyset$, dann folgt $U_i\subset W$ für alle
        $i\in I$.

\item[2)] Ist $x\in U_i$ für ein $U_i$ mit $U_i\cap V\not=\emptyset$, dann folgt
         $\sup\{\betrag{z-y}\doppelpunkt z,y \in B_{x,U_i} \} < \inf_{B_{x,U_i}} \xi$.
         Mit Lemma \auslabelnsatz{12.3} folgt hieraus $\sup_{B_{x,U_i}} \phi\le
         \gamma \inf_{B_{x,U_i}} \phi$ für jedes $x\in U_i$.

\ee

Sei $U_i\cap V\not=\emptyset$ und  $x\in U_i$, dann folgt
\formelohne{ 
   \log\betrag{f(\pi(x))}  &\le& \phi(\pi(x)) \\
                           &\le& \sup_{B_{x,U_i}}\phi  \\
                           &\le& \gamma \inf_{B_{x,U_i}}\phi \\
                           &\le& \gamma \inf_{U_i} \phi }

Also folgt
\formel{ \sup_{U_i} \log\betrag{f(\pi(x))} \le \gamma \inf_{U_i}\phi }

Wir setzen
\[   f_i(x) := \left\{ \begin{array}{cl} 
               f(\pi(x)) & \komma \mbox{falls $ x\in U_i$ und $U_i\cap V\not=\emptyset$} \\
               0         & \komma \mbox{sonst}
               \end{array} \right.  \punkt \]
Sei ferner $f_{ij}(x):= f_i(x) - f_j(x)$ für alle $x\in U_{ij}:=U_i\cap U_j$,
$i,j\in I$, $U_i\cap U_j\not=\emptyset$.

Sei $J:=\{(i,j)\in I^2\doppelpunkt U_i\cap V\not=\emptyset \mbox{ und }
              U_j\cap V=\emptyset\mbox{ oder umgekehrt } \}$. 

Es gilt für $(i,j)\in J$
\[   \sup_{U_{ij}} \betrag{f_{ij}} \le \left\{
\begin{array}{cl} \sup_{U_i} \betrag{f(\pi(x))}\mbox{, falls } & U_i\cap V\not=\emptyset \\
                  \sup_{U_j} \betrag{f(\pi(x))}\mbox{, falls } & U_j\cap V\not=\emptyset
\end{array} \right. \punkt \]

Ist $(i,j)\in I^2\backslash J$, dann folgt $f_{ij}\equiv 0$.

Sei $h:I\to \R^{+}$ eine Funktion mit $\sum_{i\in I} h(i)\le 1$. Wir setzen $h_{ij}:=
\min(h(i),h(j))$ für $(i,j)\in I^2$.

Es folgt
\formelohne{
   \sum_{(i,j)\in I^2 \atop i\not=j} h_{ij} \sup_{U_{ij}} e^{-\gamma\phi}
                         \sup_{U_{ij}}\betrag{f_{ij}} 
   &=&   \sum_{(i,j)\in J} h_{ij} \sup_{U_{ij}} e^{-\gamma\phi} 
                         \sup_{U_{ij}}\betrag{f_{ij}} \\
   &\le& \sum_{i\in I \atop U_i\cap V\not=\emptyset} h(i) \sup_{U_i} e^{-\gamma\phi}
                         \sup_{U_i} \betrag{f(\pi(z))}\\
   &\le& \sum_{i\in I} h(i) \sup_{U_i} e^{-\gamma\phi} e^{\gamma \inf_{U_i}\phi} \\
   &\le& 1   }

Sei $\J_V$ die Garbe von Keimen holomorpher Funktionen, die auf $V$ verschwinden.
$\J_V$ ist eine analytische kohärente Garbe
(siehe z.B. \cite{Hoer}, Theorem 6.5.2 i.V.m. Theorem 7.1.5).

Wir wenden nun Satz \auslabelnsatz{11.1} auf die Kokette $(f_{ij})_{(i,j)\in I^2}$ an.
Es gilt $f_{ij} \in \Gamma(U_{ij}, \J_V)$ und $\delta^1((f_{ij})_{ij}) =0$, denn
für $i,j,k\in I$ gilt $f_{ij} + f_{jk} - f_{ik} =0$ auf $U_i\cap U_j\cap U_k$.
\auslabeln{11.1.1} ist mit $\gamma \phi$ erfüllt. Sei $V_i\subset\subset U_i$, 
$i\in I$  jeweils eine offene holomorph konvexe Teilmenge und es gelte 
$\bigcup_{i\in I} V_i=\C^N$.

Satz \auslabelnsatz{11.1} liefert dann eine Kokette $(d_i)_{i\in I}$ mit $d_i \in \Gamma(
U_i, \J_V)$, $i\in I$ und $\delta^0((d_i)_i) = (f_{ij})_{ij}$ also
$d_i - d_j=f_{ij}$ auf $U_{ij}$ für alle $(i,j)\in I^2$.

Ferner liefert der Satz Konstanten $(D_i)_{i\in I}$ mit $0<D_i\le 1$, die unabhängig
von $\phi$ und $(f_{ij})_{ij}$ sind, so daß die folgende Abschätzung gilt:

\[   \sum_{i\in I} D_i \sup_{U_i} e^{-\gamma^{M+1} \phi} \sup_{V_i} \betrag{d_i}
     \le \sum_{(i,j)\in I^2 \atop i\not=j} h_{ij} \sup_{U_{ij}} e^{-\gamma\phi}
                         \sup_{U_{ij}}\betrag{f_{ij}} \le 1   \]

Sei ohne Einschränkung $D_i\le h_i$ für alle $i\in I$.

Wir erhalten zu $f$ eine Fortsetzung $F$ auf dem $\C^N$, indem wir 
\[ F(z) := f_i(z) - d_i(z)  \] für jedes $z\in U_i$, $i\in I$, setzen.
$F$ ist wohldefiniert und ganz, denn es gilt $f_i(z) - d_i(z) = f_j(z) - d_j(z)$ 
auf jeder nichtleeren Überlappungsmenge $U_{ij}$.

Wir untersuchen nun $\sup_K\betrag{F(z)}$ für ein Kompaktum $K\subset \C^N$:

Sei $I_K:=\{ i\in I\doppelpunkt V_i\cap K\not=\emptyset \}$. $I_K$ ist endlich, 
denn angenommen, es gebe eine Folge $(i_j)_j$ paarweise verschiedener Indizes in $I_K$, dann
sei $(z_j)_j$ eine Folge von Punkten aus $K$ mit $z_j\in V_{i_j}$. Sei
$z\in K$ ein Häufungspunkt von $(z_j)_j$, dann gibt es ein $l\in I_K$, so daß
$z\in V_l$. $V_l$ hat dann allerdings einen nichtleeren Schnitt mit mehr als endlich 
vielen $V_i$, $i\in I_K$, was ein Widerspruch zu der endlichen Überlappungsordnung
von $(V_i)_i$ ist.

Es gilt
\formelohne{
  \sup_K \betrag{F(z)}
  &\le& \sup_K \betrag{F(z)} \sup_{K} e^{-\gamma^{M+1}\phi} \sup_{K} e^{\gamma^{M+1}\phi}
               (\min_{i\in I_K} D_i)  (\min_{i\in I_K} D_i)^{-1}  \\
  &\le& \left( \sup_{K} e^{\gamma^{M+1}\phi} \right) (\min_{i\in I_K} D_i)^{-1}
           \sum_{i\in I_K} D_i \sup_{V_i} \betrag{F(z)} \sup_{V_i} e^{-\gamma^{M+1}\phi}\punkt }

Wegen
\formelohne{
   \sum_{i\in I_K} D_i \sup_{V_i} \betrag{f_i(z)-d_i(z)} \sup_{V_i} e^{-\gamma^{M+1}\phi}
   &\le& 1 + \sum_{i\in I_K} h_i \sup_{V_i} \betrag{f_i(z)}\sup_{V_i} e^{-\gamma^{M+1}\phi} \\
   &\le& 1 + \sum_{i\in I_K \atop V_i\cap V\not= \emptyset} h_i 
              e^{\gamma^{M+1}\inf_{U_i}\phi} e^{-\gamma^{M+1}\inf_{U_i}\phi} \\
   &\le& 2 }
erhalten wir mit $B_{K,\gamma}:=\log(2(\min_{i\in I_K} D_i)^{-1})$
\formel{
   \sup_K \betrag{F(z)} \le e^{B_{K,\gamma}}  \sup_{K} e^{\gamma^{M+1}\phi} \le 
          e^{B_{K,\gamma} + \gamma^{M+1} \sup_K C_f} \sup_{K\cap V}\betrag{f(z)}^{\gamma^{M+1}}\punkt }
Hierbei hängt die Konstante $B_{K,\gamma}$ nur von $K$, $\gamma$ und der gewählten Überdeckung 
$(U_i)_i$ ab. $M$ ist nur von $N$ abhängig.

Wir betrachten nun die Potenzen $f^m$, $m\in \N$.

Es gilt für jedes Kompaktum $K\supset\not= K_f$ mit $K\cap V\not=\emptyset$
\formelohne{
   \log\betrag{f^m(z)}   &\le& m\cdot \phi \mbox{ für jedes } z\in V \\
   \sup_K m\cdot \phi(z) &\le& \sup_{K\cap V} \log \betrag{f^m(z)} + m\sup_K C_f(z)\punkt }

Nach dem bisher Bewiesenen erhalten wir $F_m\in H(\C^N)$, $m\in \N$, mit 
\formel{\einlabeln{12.3.10}
  F_m                   &=&    f^m(z) \mbox{ für jedes $z\in V$ und } \\
  \sup_K\betrag{F_m(z)} &\le&  e^{B_{K,\gamma} + \gamma^{M+1}m \sup_K C_f} 
                           \sup_{K\cap V}\betrag{f^m(z)}^{\gamma^{M+1}}  
   } 
für jedes Kompaktum $K\supset \not= K_f$ mit 
$K\cap V \not= \emptyset$.

Sei jetzt $\phi$ eine beliebige plurisubharmonische Funktion, die i) erfüllt.
Dann wählen wir eine fallende Folge plurisubharmonischer $C^{\infty}$-Funktionen 
$(\phi_k)_{k\in \N}$,
die punktweise gegen $\phi$ strebt (siehe z.B. \cite{Klim}, Theorem 2.9.2). 

Es gibt  also nach Lemma \auslabelnsatz{12.2} 
stetige positive Funktionen $(G_k)_{k\in \N}$ auf dem $\C^N$, so daß
\[  \sup_K \phi_k \le \sup_{K\cap V}\log\betrag{f(z)} + \sup_K(C_f + G_k)  \]
für jedes Kompaktum $K\supset\not= K_f$ mit $K\cap V\not=\emptyset$, wobei
$\sup_K G_k \to 0$, falls $k\to\infty$. Wir erhalten also \auslabeln{12.3.10}
auch mit $C_f+G_k$ an Stelle von $C_f$.

Wenn wir abhängig von $K$  ein genügend großes $k$ wählen, erhalten
wir in Anwendung von \auslabeln{12.3.10} ein $F_{m,k}$, 
wobei $m\sup_K G_k$ durch ein weiteres
$\gamma$ abgefangen werden kann. 
Es gibt also für jedes $m\in \N$ ein $F_{m,K}\in H(\C^N)$ mit  
\formelohne{
  F_{m,K}                   &=&    f^m(z) \mbox{ für jedes $z\in V$ und } \\
  \sup_L\betrag{F_{m,K}(z)} &\le&  e^{B_{L,\gamma} + \gamma^{M+2}m \sup_L C_f} 
                           \sup_{L\cap V}\betrag{f^m(z)}^{\gamma^{M+2}}  }
 für jedes 
kompakte $L$ mit $K\supset L\supset\not= K_f$ und $L\cap V\not=\emptyset$.
                   
Da wir den Beweis für jedes $\gamma>1$ durchlaufen lassen können,
folgt $ii)$. 

$ii)\Rightarrow i)$: Wir wählen eine kompakte Ausschöpfung $(K_n)_{n\in \N}$ mit 
$K_1\supset K_f$ und $K_1\cap V\not=\emptyset$. Es sei $\phi_{m,n}:=
\log \betrag{F_{m,K_n}}$ für $n,m\in \N$, wobei das zugehörige $\gamma$ in ii) jeweils 
$1+\frac{1}{n} $ sein soll.

Für ein $n\in \N$ ist die Funktion
\[   \phi_n(z):= \overline{\lim_{\zeta \to z}}\, \overline{\lim_{m}} \,
                 \frac{\phi_{m,n}(\zeta)}{m}  \]
plurisubharmonisch, denn die $(\phi_{n,m}/m)_m$ sind in $m$ lokal gleichmäßig
nach oben  beschränkte plurisubharmonische Funktionen 
(vergleiche Theorem 2.9.17 in \cite{Klim}).

Sei $K_n\supset K$ und $\stackrel{\circ}{K}\supset L$ für 
kompaktes $L$ und $L\supset\not= K_f$. Ist $z$ im Inneren 
von $K$, dann folgt
\[   \phi_n(z) \le  \overline{\lim_{m}} \, \sup_{\zeta\in K}
     \betrag{\phi_{m,n}(\zeta)/m} \le (1+\frac{1}{n}) \sup_K C_f + (1+\frac{1}{n})
     \sup_{K\cap V}\betrag{f}  \]
Wegen der Stetigkeit von $C_f$ und $\betrag{f}$ folgt
\[   \inf\{ (1+\frac{1}{n}) \sup_K C_f + (1+\frac{1}{n})\sup_{K\cap V}\betrag{f}
                    \doppelpunkt K\supset 
     L  \} = (1+\frac{1}{n}) \sup_{L} C_f + (1+\frac{1}{n})
                     \sup_{L\cap V}\betrag{f} \]
und damit
\[   \sup_{L} \phi_{n} \le (1+\frac{1}{n}) \sup_{L} C_f + 
                                  (1+\frac{1}{n})\sup_{L\cap V}\betrag{f} \punkt \]
Mit der gleichen Argumentation ist 
\[  \phi(z) := \overline{\lim_{\zeta \to z}}\, \overline{\lim_{n}} \, 
                 \phi_{n}(\zeta) \]
plurisubharmonisch und  es gilt für jedes  $K\subset \C^N$ kompakt mit 
$K\cap V\not= \emptyset$ und $K\supset\not= K_f$ 
\[  \sup_K\phi \le \sup_K C_f + \sup_{K\cap V} \log\betrag{f}  \punkt \]

Ferner gilt für $z\in V$
\[   \phi_n(z) \ge \overline{\lim_{\zeta \to z\atop \zeta\in V}} \,\overline{\lim_{m}}
     \frac{m\log\betrag{f(\zeta)}}{m} = \log\betrag{f(z)} \punkt \]
Hieraus folgt $\phi(z)\ge \log\betrag{f(z)}$ für jedes
$z\in V$ \Ende

Für Anwendungen ziehen wir aus dem Beweis des Satzes \auslabelnsatz{12.1}
folgendes Ergebnis:

\begin{satz}{Satz}\einlabelnsatz{12.4}
Seien $\phi_1>0$ und $\phi_2>0$ stetige plurisubharmonische Funktionen auf dem 
$\C^N$. Sei ferner $\gamma>1$ fest. Dann gibt es zu jedem Kompaktum $K$ eine Konstante
$B_K$, so daß zu jedem $f\in H(V)$ und jedem $\alpha >0$ mit 
\formel{
   \log \betrag{f(z)} \le \phi_1(z) +\alpha \phi_2(z) \mbox{ für jedes }z\in V   
\einlabeln{12.4.1}}

eine ganze Funktion $F$ existiert mit 
\be
\item[a)] $F(z)=f(z)$ für jedes $z\in V$,
\item[b)] $\sup_K \betrag{ F(z)} \le e^{B_K} e^{\gamma \sup_K \phi_1+ \alpha\phi_2}$ 
          für jedes Kompaktum $K$.
\ee
Hierbei hängt $B_K$ zwar von $\phi_1$ und $\phi_2$ ab, aber nicht von $\alpha$.
\end{satz}

Es gibt nach doppelter Anwendung von Lemma \auslabelnsatz{12.3} eine Funktion $\xi$,
so daß für jedes $U\subset \Omega$ mit 
\[ \sup \{ \betrag{x-y} \doppelpunkt x,y \in U\} < \inf_U \xi(z)  \]
folgt
\[  \sup_U \phi_1 \le \gamma \inf \phi_1 \mbox{ und } 
                     \sup_U \phi_2 \le \gamma \inf \phi_2  \punkt \] 
Für solche $U$ gilt dann aber 
\formel{ 
     \sup_U (\phi_1 + \alpha\phi_2 )
   &\le& \sup_U \phi_1 + \sup_U \alpha\phi_2 \\
   &\le& \gamma \inf_U \phi_1 + \alpha\gamma\inf_U \phi_2  \\
   &\le& \gamma \inf_U (\phi_1 + \alpha\phi_2)   }
für jedes $\alpha\in \R_0^{+}$.

In dem Beweis zu Satz \auslabelnsatz{12.1} kann daher für jedes $\alpha>0$ die gleiche
Überdeckung $(U_i)_i$ gewählt werden.

Ist $f\in H(V)$ mit \auslabeln{12.4.1} gegeben, dann existiert nach \auslabeln{12.3.10}
für jedes $\gamma>1$ eine
ganze Funktion $F$ mit 

\[   \sup_K \betrag{F(z)} \le e^{B_{K,\gamma}} 
                         \sup_K e^{\gamma^{M+1} (\phi_1+\alpha\phi_2)}, \]
wobei $B_{K,\gamma}$ unabhängig von $f$ und $\alpha$ ist \Ende

\newpage
\section[{\sc Existenz eines linear zahmen 
Fortsetzungsoperators}]{Ein Kriterium für die Existenz eines linear zahmen 
Fortsetzungsoperators}

Sei $V$ eine in den $\C^N$ eingebettete zusammenhängende abgeschlossene
Untermannigfaltigkeit. $\phi$ sei eine stetige plurisubharmonische 
Ausschöpfungsfunktion des $\C^N$ mit $\inf_{z\in V} \phi(z)<1$.

Sei 
\[  K_n:= \overline{\{z\in \C^N\doppelpunkt \phi(z) < n\} } 
             \mbox{ für }n\in \N \]
und 
\formelohne{
    \norm{n}{F} &:=& \sup_{z\in K_n}\betrag{F(z)}\mbox{ für } 
                        F\in H(\C^N) \mbox{ bzw. }  \\
    \betrag{f}_n &:=& \sup_{z\in K_n\cap V}\betrag{f(z)}\mbox{ für } 
                        f\in H(V)   }
Wir fixieren die Halbnormensysteme $(\normleer_n)_n$ und 
$(\betragleer_n)_n$ und behandeln $(H(\C^N),\normleer_n)$ bzw.
$(H(V),\betragleer_n)$ von nun an als gradierte Fr\'echeträume.

\begin{satzohnebeweis}{Definition}
Seien $(X,\normleer_n)$ und $(Y,\betragleer_n)$ gradierte Fr\'echeträume,
dann bezeichnen wir einen Operator $T:(X,\normleer_n)\to (Y,\betragleer_n)$
als linear zahm, falls es Konstanten $\alpha\in \N$ und $\beta\in \N_0$ gibt,
so daß zu jedem $k\in \N$ ein $C_k>0$ existiert mit
\[   \betrag{T(x)}_k \le C_k \norm{\alpha k + \beta}{x} \foralls x\in X \punkt\]

Kann $\alpha=1$ gewählt werden, so heißt der Operator zahm.

Zwei gradierte Fr\'echeträume $X,Y$ heißen linear zahm bzw. zahm isomorph
zueinander, falls es einen Isomorhismus $\Phi:X\to Y$ gibt, so daß
$\Phi$ und $\Phi^{-1}$ linear zahm bzw. zahm sind. Zwei Gradierungen
$(\betragleer_n)_n$ und $(\normleer_n)_n$ auf $X$ heißen linear zahm
bzw. zahm äquivalent zueinander, falls die Identität $(X,\betragleer_n)\to
(X,\normleer_n)$ ein linear zahmer bzw. zahmer Isomorphismus ist.

\end{satzohnebeweis}

Die Restriktionsabbildung $R:(H(\C^N),\normleer_n)\to (H(V),\betragleer_n)$
ist ein zahmer Operator.

Sei $0\le a_0\le a_1\le \ldots$ eine wachsende Folge von Zahlen mit
$a_k \sim k^{\frac{1}{d}}$ für ein festes $d\in \N$, d.h. es
gibt eine Konstante $C>0$, so daß 
\[   \frac{1}{C} k^{\frac{1}{d}} \le a_k \le C k^{\frac{1}{d}} 
                   \foralls k\in \N \punkt \]

Der Potenzreihenraum unendlichen Typs, definiert durch 
\[   \Lambda_{\infty}(a_k) := \{ x\in \C^{\N_0}\doppelpunkt \norm{n}{x}^\infty := 
        \sum_{k=0}^{\infty} \betrag{x_k} e^{na_k} < \infty \}  \]
ist dann nuklear und zahm isomorph zu $\Lambda_\infty(k^{\frac{1}{d}})$ (vergleiche
hierzu \cite{M-V}, Satz 29.6).

{\sc Bemerkung:} Ist $V$ eine Steinsche Mannigfaltigkeit der Dimension $d$
und ist $H(V)$ isomorph zu einem Potenzreihenraum unendlichen Typs 
$\Lambda_{\infty}(a_k) $, dann ist $a_k \sim k^{\frac{1}{d}}$ (siehe hierzu
z.B. \cite{Ayt2}, Proposition I.1).

\begin{satz}{Satz}\einlabelnsatz{13.1a}
Zu jedem linear zahmen Operator $T: (\Lambda_\infty(a_k),\normleer_n^\infty)\to
(H(V),\betragleer_n)$ gibt es ein linear zahmes Lifting $L:
(\Lambda_\infty(a_k),\normleer_n^\infty) \to (H(\C^N),\normleer_n)$ mit
$T=R\circ L$
\end{satz}

Sei $f_k := T e_k$, wobei $e_k$ der $k$-te Einheitsvektor in $\Lambda_{\infty}(a_k)$
sei, $k\in \N_0$. Dann folgt 
\[   \sup_{K_n\cap V}\betrag{f_k}\le C_n e^{(\alpha n + \beta)a_k} \]
für jedes $n\in \N$ und jedes $k\in \N_0$.

Wir wählen eine stetige plurisubharmonische Funktion $\phi_1$ auf dem
$\C^\N$, so daß für jedes $n\ge 2$ gilt
\[   \inf_{z\notin K_{n-1}} \phi_1(z) \ge \log C_n \und \phi_1(z) \ge \log C_1
              \foralls z\in \C^N \punkt \]
Ferner sei eine plurisubharmonische Funktion $\phi_2$ für jedes $z\in \C^N$
durch 
\[   \phi_2(z) := \max(\alpha(\phi(z) + 1) + \beta, \alpha+\beta )\]
definiert. Für diese gilt
\[   \sup_{K_n} \phi_2 \le \alpha(n+1)+\beta \foralls n\in \N \punkt \]

Da wir Satz \auslabelnsatz{12.4} anwenden wollen bei festem $\gamma>1$, 
müssen wir \auslabeln{12.4.1} zeigen:

Sei $z\in (K_n\backslash K_{n-1})\cap V$, $n\ge 2$, dann folgt für jedes
$k\in \N_0$

\[   \log\betrag{f_k(z)} \le \sup_{K_n\cap V}\log\betrag{f_k} \le \log C_n
                                    + (\alpha n + \beta)a_k  \punkt\]

Andererseits ist für jedes $z\notin K_{n-1}\cap V$
\[  \phi_1(z) \ge \log C_n \und \phi_2(z) \ge \alpha((n-1)+1) +\beta \]

Ist $z\in K_1\cap V$, dann folgt
\[   \log\betrag{f_k(z)} \le \sup_{K_1\cap V}\log\betrag{f_k(z)} \le
     \log C_1 + (\alpha + \beta)a_k \le \phi_1(z) + a_k \phi_2(z) \]

Hieraus folgt insgesamt 
\[   \log \betrag{f_k(z)} \le \phi_1(z) + a_k \phi_2(z) \foralls z\in V \]

 Satz \auslabelnsatz{12.4} liefert also für jedes $k\in \N_0$ eine ganze Funktion $F_k$ mit
\be
\item[a)] $F_k(z) = f_k(z)$ für jedes $z\in V$, $k\in \N_0$,
\item[b)] $\sup_{K_n} \betrag{F_k(z)}\le e^{B_{K_n} + \gamma \sup_{K_n}( \phi_1
                                            + a_k\phi_2}) $ 
          für jedes $n\in \N$ und $k\in \N_0$.
\ee

Dies wenden wir für $\gamma=2$ an und erhalten eine Folge ganzer Funktionen
$(F_k)_k$ mit a) und b).

Wir setzen für eine Folge $(\lambda_k)_k \in \Lambda_\infty(a_k)$ 
\[   L((\lambda_k)_k):= \sum_{k=0}^\infty \lambda_k F_k \punkt \]

Durch die folgende Abschätzung wird die Wohldefiniertheit und 
linear Zahmheit gleichzeitig gezeigt. Die Linearität von $L$ ist offensichtlich.
Für jedes $n\in \N$ gilt 

\formelohne{
\norm{n}{L((\lambda_k)_k)}
               &=& \sup_{K_n} \betrag{\sum_{k=0}^\infty \lambda_k F_k} \\
               &\le& \sum_{k=0}^\infty \betrag{\lambda_k} 
                     e^{B_{K_n} + 2 \sup_{K_n} (\phi_1 + a_k\phi_2)} \\
               &\le& e^{B_{K_n}} e^{2\sup_{K_n} \phi_1}
                     \sum_{k=0}^\infty \betrag{\lambda_k} e^{2(\alpha(n+1) + \beta) a_k} \\
               &=&   D_n \norm{2\alpha n + 2(\alpha + \beta)}{(\lambda_k)_k}\punkt  }
   
Der Satz ist bewiesen \Ende

Wir zeigen nun das angekündigte Ergebnis: 

\begin{satz}{Korollar}\einlabelnsatz{13.2}
Ist $(H(V),\betragleer_n)$ linear zahm isomorph zu einem Potenzreihenraum unendlichen
Typs $\Lambda_\infty(a_k)$, dann gibt es einen linear zahmen Ausdehnungsoperator
$(H(V),\betragleer_n)\to (H(\C^N),\normleer_n)$.

\end{satz}
Sei $T: (\Lambda_\infty(a_k),\normleer_n^\infty)\to
(H(V),\betragleer_n)$ ein  linear zahmer Isomorphismus.

 Wir verwenden das im Beweis zu 
Satz \auslabelnsatz{13.1a} konstruierte linear zahme Lifting $L:
(\Lambda_\infty(a_k),\normleer_n^\infty) \to (H(\C^N),\normleer_n)$ mit
$T=R\circ L$. Wir setzen nun $E:=L\circ T^{-1}$. E ist daher linear zahm und 
es gilt 
\[   R\circ E( \sum_{k=0}^\infty \lambda_k f_k)=R\circ L((\lambda_k)_k)=
     R(\sum_{k=0}^\infty \lambda_k F_k)=\sum_{k=0}^\infty \lambda_k f_k\punkt \]
Also gilt $R\circ E=id_{H(V)}$ \Ende

\newpage
\section[{\sc Konvexität plurisubharmonischer Funktionen}]{Konvexität
plurisubharmonischer Funktionen auf abgeschlos\-senen
Untermannigfaltigkeiten des $\C^N$} \label{Kap10}

In diesem Kapitel wollen wir die Konvexität plurisubharmonischer Funktionen
auf Steinschen Untermannigfaltigkeiten $V$ des $\C^N$ untersuchen und zwar relativ
zu den Kugeln im $\C^N$ mit Radius $r$ und Mittelpunkt im Ursprung.


Wir werden die Algebraizität von $V$ in einen Zusammenhang mit der Qualität der 
Konvexität plurisubharmonischer Funktionen auf $V$  bringen.

Ausgehend von der Existenz eines linear zahmen Ausdehnungsoperators 
verwenden wir die Argumentation von A. Aytuna in \cite{Ayt3}, um eine
Konvexitätsbedingung für plurisubharmonische Funktionen auf $V$ herzuleiten.
Während A. Aytuna direkt die zur Algebraizität von $V$ äquivalente 
Bedingung von A. Sadullaev (siehe \cite{S1}, Theorem 2.2, bzw. Kapitel 1) zeigt,
werden wir hier ferner zeigen, daß die Bedingung von A. Sadullaev äquivalent zu der
angesprochenen Konvexitätsbedingung für plurisubharmonische Funktionen auf $V$
ist. 

Wir benötigen zunächst das folgende Lemma (siehe \cite{Rock}, Corollary 17.1.5),
das direkt aus dem 
Satz über konvexe Hüllen von Caratheodory  folgt (siehe \cite{Rock}, Theorem 17.1).

Hier können wir uns auf die Situation in der reellen Ebene beschränken:

\begin{satzohnebeweis}{Lemma}\einlabelnsatz{16.1} 
Sei $f:\R^+_0\to \R$ eine monoton wachsende Funktion und sei $E:=\{
(x,f(x))\doppelpunkt x\in \R^+_0 \}$. $ch(E)$ sei die konvexe Hülle
von $E$. Sei ferner $h(x):= \inf\{ y\doppelpunkt (x,y)\in ch(E)\}$.
Dann ist $h$ eine konvexe Funktion und es gilt 
\[   h(x)=\inf\{ \lambda f(x_1) + (1-\lambda)f(x_2)\doppelpunkt x=\lambda x_1
+(1-\lambda)x_2\komma \lambda\in [0,1]\komma x_1,x_2\in \R^+_0 \}   \]
\end{satzohnebeweis}

Hieraus folgt eine Aussage über eine schwache Form von Konvexität:

\begin{satz}{Proposition}\einlabelnsatz{16.2}
Sei $f:\R^+_0\to \R$ eine monoton wachsende Funktion und seien
$a,b\in \R^+_0$, $a\ge 1$. Dann sind äquivalent:

\begin{enumerate}
\item[1)] $f(\lambda x_1+(1-\lambda)x_2)\le \lambda f(ax_1+b) + (1-\lambda)f(ax_2 + b)$ für alle
$\lambda \in [0,1]$, $x_1,x_2\in \R^+_0$.
\item[2)] Es gibt ein konvexe Funktion $h$, so daß 
\formelohne{   h(x)&\le& f(x) \mbox{ für alle } x\ge 0 \mbox{ und } \\
       f(x)&\le& h(ax+b) \mbox{ für alle } x\ge 0\punkt }
\end{enumerate}
\end{satz}
 2)$\Rightarrow$1): Sei $x_1,x_2\in \R^+_0$ und $\lambda\in [0,1]$, dann gilt
\formelohne{ f(\lambda x_1+(1-\lambda)x_2) 
     &\le& h(a(\lambda x_1+(1-\lambda)x_2)+b) \\
     &\le& \lambda h(ax_1+b) + (1-\lambda)h(ax_2+b) \\
     &\le& \lambda f(ax_1+b) + (1-\lambda)f(ax_2+b) \punkt }
 1)$\Rightarrow$2):
Wir wählen zu $f$ die Funktion $h$ aus Lemma \auslabelnsatz{16.1}. $h$ ist
konvex und liegt unterhalb von $f$. Sei $x\in \R^+_0$ fest. Ist $\eps>0$, dann 
können wir nach Lemma \auslabelnsatz{16.1} $x_1,x_2\in \R^+_0$ wählen, so daß
$ax+b=\lambda x_1+(1-\lambda)x_2$ und $\lambda f(x_1) + (1-\lambda)f(x_2)\le h(ax+b)+\eps$.
Wegen 1) ist 
\formelohne{ \lambda f(x_1) + (1-\lambda)f(x_2)
   &\ge & f(\lambda(\frac{x_1-b}{a}) + (1-\lambda) (\frac{x_2-b}{a})) \\
   &= & f(x) \punkt } 
Für $\eps\searrow 0$ erhalten wir $f(x)\le h(ax+b)$. Da $x\in \R^+_0$ 
beliebig gewählt war, folgt 2) \Ende

Sei $\Phi\in PSH(\C^N)$, dann sei $M_{\Phi}(r):=\sup \{ \Phi(z)\doppelpunkt
\betrag{z}=e^r \}$ für $r\in \R^+_0$. $M_{\Phi}(r)$ ist eine konvexe 
wachsende Funktion nach dem Hadamarschen Drei-Kreise-Satz für plurisubharmonische
Funktionen (siehe z.B. \cite{Hoer1}, Theorem 4.1.13).

Sei $V$ eine abgeschlossene Untermannigfaltigkeit des $\C^N$ und sei
durch die Spuren der Kugeln $B_r:=\{ z\in \C^N\doppelpunkt \betrag{z}<e^r\}$, 
$r\in \R^+_0$,
auf $V$ eine Gradierung $\norm{n}{f}:=\sup_{B_n\cap V} \betrag{f}$, $n\in \N$,
$f\in H(V)$, für $H(V)$ gegeben.

Es gilt der folgende Satz:

\begin{satz}{Satz} \einlabelnsatz{16.3}
Sei $E:H(V) \to H(\C^N)$ ein linear zahmer Ausdehnungsoperator, d.h. es
gebe $a,b\in \N_0$, $a\ge 1$, so daß für alle $n\in \N$ ein $C_n>0$ existiert mit
\formel{  \sup_{B_n} \betrag{E(f)(z)} \le C_n \norm{an+b}{f} \foralls f\in H(V)
        \einlabeln{16.3.1}\punkt }
Dann gibt es zu jedem $\phi\in PSH(V)$ ein $\Phi\in PSH(\C^N)$ mit 
\formelohne{  \sup_{B_r\cap V} \phi(z) &\le& M_{\Phi}(r) \mbox{ und } \\
      M_{\Phi}(r) &\le& \sup_{B_{ar+a+b}\cap V} \phi(z) \foralls r\ge 0 \punkt }
\end{satz}

Jedes $\psi\in PSH(\C^N)$ besitzt nach einem Ergebnis von Bremermann (\cite{Brem},
Theorem 2)
eine Hartogsdarstellung auf dem $\C^N$ 
\[    \psi(z)= \overline{\lim_n}\,  \frac{\log\betrag{f_n(z)}}{c_n} \mbox{ mit }
      f_n \in H(\C^N) \foralls n\in \N \punkt \]
Ohne Einschränkung gilt $c_n\nearrow \infty$. Der Fortsetzungssatz für 
plurisubharmonische Funktionen auf einer abgeschlossenen komplexen Untermannigfaltigkeit
des $\C^N$ von Sadullaev (siehe \cite{S-C}, Theorem 1.1) zeigt, daß auch
jedes $\phi\in PSH(V)$ eine Hartogsdarstellung
\formel{    \phi(z)= \overline{\lim_n}\, \frac{\log\betrag{f_n(z)}}{c_n} \mbox{ mit }
      f_n \in H(V) \foralls n\in \N \einlabeln{16.3.2} }
besitzt. Sei also $\phi\in PSH(V)$ mit einer Darstellung wie in \auslabeln{16.3.2}
gegeben. 

Wir zeigen, daß $\{ \frac{\log\betrag{E(f_n)}}{c_n} \doppelpunkt n\in \N \} $ eine
lokal gleichmäßig beschränkte Familie ist.

Sei $D_r:=\sup_{B_r\cap V} \phi(z)$, $r\ge 0$. Nach dem Hartogslemma
(siehe z.B. \cite{S1}, Theorem 1.4) gibt es zu jedem $\eps>0$ und 
jedem $r'<r$ ein $n_0$, so daß für alle $n\ge n_0$
\[   \sup_{B_{r'}\cap V} \frac{\log\betrag{f_n(z)}}{c_n} \le D_r + \eps\punkt \]

Wir setzen $b':=a+b$. Für jedes $r\in \R^+_0$ liefert \auslabeln{16.3.1} für
jedes $f\in H(V)$ die Abschätzung
\[   \sup_{B_r} \betrag{E(f)(z)} \le C_r \norm{ar+b'}{f} \komma \]
wobei $C_r=C_m$ für $m-1\le r\le m$ ist.

Wir erhalten zu jedem $\eps>0$ und jedem $r'<r$ ein $n_0$, so daß für
$n\ge n_0$ gilt
\formelohne{  \sup_{B_{r'}}\frac{\log\betrag{E(f_n)(z)}}{c_n}
     &\le& \frac{\log C_{r'}}{c_n} + \sup_{ B_{ar'+b'}\cap V}
                        \frac{\log\betrag{f_n(z)}}{c_n}   \\
     &\le& \frac{\log C_{r'}}{c_n} + D_{ar+b'} + \eps \punkt }

Also ist $\{ \frac{\log\betrag{E(f_n)}}{c_n} \doppelpunkt n\in \N \} $  eine
lokal gleichmäßig beschränkte Familie. Daher ist 
\[   \Phi(z):= \overline{\lim_{\zeta\to z}}\, \overline{\lim_n}\, \frac{\log\betrag{E(f_n)(\zeta)}}{c_n}\]
eine plurisubharmonische Funktion auf dem $\C^N$ (siehe z.B. \cite{Klim}, Theorem
2.9.17).

Für $r'<r''<r$ und jedes $\eps >0$ gilt dann
\formelohne{  M_{\Phi}(r')  
   &\le& \sup \{ \Phi(z)\doppelpunkt \betrag{z} < e^{r''} \}  \\
   &\le& \sup_{B_{r''}} \overline{\lim_n} \frac{\log\betrag{E(f_n)(z)}}{c_n} \\
   &\le& \overline{\lim_n}\, ( \frac{\log C_{r''}}{c_n} + D_{ar+b'} + \eps) \\
   &=&   D_{ar+b'} + \eps \punkt}

Da $M_{\Phi}(r)$ stetig ist in $r$, folgt insgesamt $M_{\Phi}(r)\le 
\sup_{B_{ar+a+b}\cap V} \phi(z)$.

Wegen des Maximumprinzips für plurisubharmonische Funktionen auf $V$ gilt
umgekehrt auch $\sup_{B_r\cap V} \phi(z) \le M_{\Phi}(r)$ für $r>0$ \Ende 

Man beachte, daß wegen der Regularisierung zwar $\phi(z)\le\Phi(z)$ für jedes
$z\in V$, aber nicht notwendig $\phi(z)=\Phi(z)$ gilt.

Sei für ein $\phi\in PSH(V)$ eine Funktion $m_{\phi} :\R^+_0\to \R$ durch
$ m_{\phi}(r):= \sup_{B_r\cap V} \phi(z)$ definiert, wobei wie oben 
$B_r:=\{ z\in \C^N\doppelpunkt \betrag{z}<e^r\}$ sei.

Kombiniert mit Proposition \auslabelnsatz{16.2} ergibt sich

\begin{satzohnebeweis}{Satz} \einlabelnsatz{16.4}
Unter den Voraussetzungen von Satz \auslabelnsatz{16.3} gilt  
für jedes $\phi\in PSH(V)$ und jedes $\lambda\in [0,1]$, sowie
jedes $r_1,r_2\in \R^+_0$
\formel{   m_{\phi}(\lambda r_1 + (1-\lambda)r_2) \le \lambda m_{\phi}(ar_1 + b') + (1-\lambda)
     m_{\phi}(ar_2+b')  \komma \einlabeln{16.4.1} }
wobei $a,b$ aus Satz \auslabelnsatz{16.3} sind und $b'=a+b$.
\end{satzohnebeweis}

Wir zeigen nun, daß aus \auslabeln{16.4.1} bereits folgt, daß $V$ algebraisch ist.

\begin{satz}{Satz} \einlabelnsatz{16.5}
Sei $V$ eine abgeschlossene Untermannigfaltigkeit des $\C^N$ und für
jede plurisubharmonische Funktion auf $V$ gelte \auslabeln{16.4.1}, dann
ist $V$ algebraisch.
\end{satz}

Wir zeigen, daß die Siciaksche Extremalfunktion $S(z,K)$ für mindestens ein
Kompaktum $K\subset V$ lokal beschränkt ist auf $V$.
Dann können wir ein Ergebnis von Sadullaev anwenden, nach dem dann und nur dann
$V$ algebraisch ist. Sadullaevs Ergebnis ist allgemeiner als hier
benötigt, es erfaßt  nämlich auch  
Varietäten, die singuläre Punkte besitzen (\cite{S1}, Theorem 2.2).

$S(z,K)$ ist für ein Kompaktum $K\subset V$ wie folgt definiert:

\[  S(z,K):= \overline{\lim_{\zeta\to z}}\, \sup\{ u(\zeta)\doppelpunkt u\le 0 \mbox{ auf }
      K \komma u\in L\} \komma  \]
wobei
\[  L:=\{ u\in PSH(\C^N)\doppelpunkt \mbox{ es gibt } \alpha=\alpha(u)\doppelpunkt
       u(z)\le \log(1+\betrag{z}) + \alpha \foralls z\in \C^N \} \punkt \]

Sei $u\in L$, $u\le 0$ auf $\overline{B_0\cap V}$. Sei $\phi:=u\mid_{\textstyle V}$. 
Aus Proposition \auslabelnsatz{16.2} folgt die Existenz einer
konvexen Funktion $h$, so daß $h(r)\le m_{\phi}(r)$ und $m_{\phi}(r)\le 
h(ar+b)$ für alle $r\ge 0$.
Also ist
\formelohne{ h(r) &\le& \sup\{ \log(1+\betrag{z}) + \alpha(u) \doppelpunkt 
    \betrag{z}< e^r \}  \\
   &\le& \log(1+e^r) + \alpha(u)  \\
   &\le& r + \alpha(u) + \log 2 }
für $r\ge 0$.

Damit ist $h$ eine lineare Funktion mit Steigung kleiner
oder gleich 1 auf $\R^+_0$ und wegen
$h(0)\le m_{\phi}(0)\le 0$ liegt $h$ unterhalb der Identität $\R^+_0 \to \R^+_0$.

Ist $z\in V$, dann wählen wir $r\ge 0$, so daß $z\in B_r\cap V$.
Also gilt 
\[  u(\zeta) = \phi(\zeta) \le m_{\phi}(r) \le h(ar +b) \le ar +b \]
für alle $\zeta\in V$ nahe $z$. Damit ist $S(z,\overline{B_0\cap V})$ lokal
beschränkt in $V$ \Ende


\newpage \addcontentsline{toc}{section}{\sc Literatur}
\begin{thebibliography}{99}

\bibitem{Ayt1} Aytuna, A., {\it On the linear topological structure of 
spaces of analytic functions}, Do\v{g}a Turkish J. Math {\bf 10} (1986), 46-49.

\bibitem{Ayt2} Aytuna, A., {\it Stein spaces for which O(M) is isomorphic
to a power series space}, Advances in the theory of Fr\'echet spaces (ed: 
T. Terzio\v{g}lu), 115-154, Kluwer Academic Publishers, 1989.


\bibitem{Ayt3} Aytuna, A., {\it Linear tame extension operators from closed
subvarieties of $\C^d$}, Proc. Amer. Math. Soc. {\bf 123} (1995), 759-763.

\bibitem{A-K-T} Aytuna, A., Krone, J., Terzio\v{g}lu, T., {\it Complemented infinite
type power series subspaces of nuclear Fr\'echet spaces}, Math. Ann. {\bf 283}
(1989), 193-202.

\bibitem{Brem} Bremermann, H.J., {\it On the conjecture of equivalence of
plurisubharmonic functions and Hartogs functions}, Math. Ann. {\bf 131}, 76-86.

\bibitem{S-C} Cegrell, U., Sadullaev, A., {\it Approximation of plurisubharmonic
functions and the Dirichlet problem for the complex Monge-Amp\`ere operator},
Math. Scand. {\bf 71} (1992), 62-68. 

\bibitem{Eng} Engelking, R., {\it General Topology}, Heldermann, Berlin, 1998.

\bibitem{Hoer} Hörmander, L., {\it An introduction to complex analysis in
several variables}, North Holland, Amsterdam, 1990, 3rd edition.

\bibitem{Hoer1} Hörmander, L., {\it Notions of convexity}, Birkhäuser,
Boston-Basel-Berlin, 1994.

\bibitem{G-V} Grauert, H., Fritzsche, K., {\it Einführung in die Funktionentheorie
mehrerer Veränderlicher}, Springer, Berlin-Heidelberg-New York, 1974.

\bibitem{G-R} Gunning, R., Rossi, H., {\it Analytic functions of several complex variables},
Prentice-Hall, 1965.

\bibitem{Klim} Klimek, M., {\it Pluripotential theory}, Clarendon Press, Oxford-New 
York-Tokyo, 1991, London Math. Soc. Monographs N.S. 6.

\bibitem{Koethe} Köthe, G., {\it Topological vector spaces II}, Grundlehren der
Math. Wiss. 257, Springer, Berlin-Heidelberg-New York, 1979.

\bibitem{Kunz} Kunz, E., {\it Introduction to Commutativ Algebra und Algebraic 
Geometry}, Birkhäuser, Boston-Basel-Stuttgart, 1985.

\bibitem{M-V} Meise, R., Vogt, D., {\it Einführung in die Funktionalanalysis}, Vieweg,
Braunschweig-Wiesbaden, 1992.

\bibitem{M-Venglisch} Meise, R., Vogt, D., {\it Introduction to Functional Analysis}, 
Clarendon Press, Oxford, 1997.

\bibitem{M-H} Mityagin, V.S., Khenkin, G.M., {\it Linear problems of complex analysis},
Russian Math. Surveys {\bf 26} (1971), 99-164.

\bibitem{Rock} Rockafellar, T., {\it Convex Analysis}, Princeton Univ. Press, 
1970.

\bibitem{S1} Sadullaev, A., {\it An estimate for polynomials on analytic sets}, 
Math. USSR Izvestiya {\bf 20} (1983), 493-502.

\bibitem{Siu} Siu, Y., {\it Every Stein subvariety admits a Stein neighborhood}, Inv. math.
{\bf 38} (1976), 89-100.

\bibitem{V1} Vogt, D., {\it On the functor Ext$^1(E,F)$ for Fr\'echet spaces}, Studia
 Math. {\bf 85} (1987), 163-197.

\bibitem{V2} Vogt, D., {\it Power series spaces representations of 
nuclear Fr\'echet spaces}, Trans. Amer. Math. Soc. {\bf 319} (1990), 191-208.

\bibitem{Z} Zaharyuta, V.P., {\it Isomorpisms of spaces of analytic functions},
Sov. Math. Dokl. {\bf 22} (1980), 631-634.

 


\end{thebibliography}
\end{document}